\documentclass[11pt,reqno]{amsproc}
\usepackage[utf8]{inputenc}
\usepackage[T1]{fontenc}
\usepackage{lmodern}
\usepackage{amsmath, amssymb, amsthm}
\usepackage{hyperref}
\usepackage{geometry}
\usepackage{tikz}
\usetikzlibrary{arrows.meta}
\usepackage{pgfplots}
\pgfplotsset{compat=1.18}
\usepgfplotslibrary{groupplots}
\title[Matrix Kawada--It\^o theorem and harmonic rigidity]{The Matrix Kawada--It\^o Theorem and Harmonic Rigidity on Locally Compact Groups}
\author{Ali Jabbari}
\address{Department of Applied Mathematics and Informatics, Kyrgyz-Turkish Manas University, Bishkek, Kyrgyzstan}
\email{ali.jabbari@manas.edu.kg, jabbari\_al@yahoo.com}
\thanks{Corresponding author: Ali Jabbari, Department of Applied
Mathematics and Informatics, Kyrgyz-Turkish Manas University, Bishkek,
Kyrgyzstan. E-mail: \texttt{ali.jabbari@manas.edu.kg}.}

\theoremstyle{plain}
\newtheorem{theorem}{Theorem}[section]
\newtheorem{proposition}[theorem]{Proposition}
\newtheorem{lemma}[theorem]{Lemma}
\newtheorem{corollary}[theorem]{Corollary}
\theoremstyle{definition}
\newtheorem{definition}[theorem]{Definition}
\newtheorem{example}[theorem]{Example}

\theoremstyle{remark}
\newtheorem{remark}[theorem]{Remark}

\newcommand{\C}{\mathbb{C}}
\newcommand{\R}{\mathbb{R}}

\newcommand{\tr}{\operatorname{Tr}}
\newcommand{\Fix}{\operatorname{Fix}}
\newcommand{\Ran}{\operatorname{Ran}}
\newcommand{\Lp}{L^p(G,M_n)}
\newcommand{\hsp}{H_S^p(G, M_n)}

\begin{document}
\begin{abstract}
A positive $M_n$-valued measure on a locally compact group need not have
positive convolution products: matrix multiplication preserves neither
self-adjointness nor positivity.  We prove that this obstruction disappears
completely at an idempotent.  If $\mu\in M(G,M_n)$ is positive and
$\mu*\mu=\mu$, then all its matrix values commute and there is a unique
finite decomposition
\[
  \mu=\sum_{j=1}^r m_{K_j}\otimes E_j,
\]
where the $K_j$ are distinct compact subgroups and the $E_j$ are
pairwise orthogonal nonzero projections.  The result holds for arbitrary
locally compact groups, without countability assumptions or a priori
commutativity.  The proof combines projection-valued dilations,
rigidity of contractive idempotent Fourier coefficients, and
faithfulness of the regular representation on finite measures over
$G\times G$.  For normalized positive matrix convolution semigroups,
we further identify the Ces\`aro harmonic projection with convolution
by such an idempotent, even if mass escapes to infinity.  Its surviving
matrix directions are exactly those uniformly confined to a compact
subset of the group, yielding an explicit direct Haar-channel
decomposition for all $1<p<\infty$ (orthogonal in $L^2$).  A bounded Zeno product
formula and independent density, atomic, and subgroup-support
classifications illustrate the rigidity.
\end{abstract}
\maketitle

\section{Introduction}\label{sec:intro}

The Kawada--It\^o theorem identifies the idempotent probability measures
on a locally compact group with normalized Haar measures of compact
subgroups.  In the scalar problem the positivity assumption interacts
well with convolution, and the support of an idempotent is forced to be
a compact group.  For matrix-valued measures even the most elementary
part of that argument becomes problematic.  If $A,B\in M_n$ are positive,
the product $AB$ need not be self-adjoint; consequently convolution of
positive matrix measures need not be positive.  This obstruction is
present before one considers Fourier analysis or the topology of the
underlying group.

We establish an exact rigidity principle: \emph{the assumption that a
positive matrix measure is convolution-idempotent forces all its values
to commute}.  Once the coefficient algebra becomes commutative, the
classical theorem applies simultaneously to its finitely many
spectral channels.  Neither commutativity nor a simultaneous choice of
eigenvectors is assumed in advance.  The group need not be abelian,
compact, second countable, or unimodular.

\subsection*{The classification theorem}
Let $M(G,M_n)$ denote the finite regular matrix-valued Borel measures
on a locally compact group $G$.  Positivity is meant in the
\emph{measure} sense: $\mu(D)$ is a positive semidefinite matrix for
every Borel set $D$.  Convolution uses ordinary multiplication of the
matrix coefficients in the displayed order.  Our principal theorem is
the following (Theorem~\ref{thm:general-matrix-kawada}).

\textbf{Classification Theorem.}
\emph{If $\mu\ge0$ and $\mu*\mu=\mu$, there exist distinct compact
subgroups $K_1,\ldots,K_r\le G$ and pairwise orthogonal nonzero
projections $E_1,\ldots,E_r\in M_n$ such that}
\begin{equation}\label{eq:intro-classification}
   \mu=\sum_{j=1}^r m_{K_j}\otimes E_j,
   \qquad \sum_{j=1}^rE_j=\mu(G).
\end{equation}

\emph{The decomposition is unique up to permutation; $r\le n$ and the
zero measure corresponds to the empty family.  In particular,
$[\mu(D),\mu(D')]=0$ for all Borel sets $D,D'$.}
\medskip

The extension to measures with coefficients in any finite-dimensional
$C^*$-algebra follows blockwise
(Corollary~\ref{cor:finite-dimensional-cstar-kawada}).

The projections $E_j$ are intrinsic: their ranges are precisely the
subspaces on which the vector measure $D\mapsto\mu(D)u$ equals
$m_{K_j}(D)u$; the subgroup $K_j$ is its stabilizer under left
translation.  Thus the classification gives more than a basis in
which the measure happens to be diagonal: it describes all its
channels without a choice of basis.  The finite union of compact
subgroups is the support of $\mu$, while $\ker\mu(G)$ is a distinguished
inactive coefficient space.  The relevant subgroups need not be
nested, normal, or pairwise disjoint.

This statement should be distinguished from the classification of
\emph{complex} scalar convolution idempotents, for which signed
combinations and coset-ring phenomena occur
\cite{CohenIdempotent,HostIdempotents}.  The positivity hypothesis is
indispensable.  Even a scalar compression
$D\mapsto\langle\mu(D)v,v\rangle$ of a matrix idempotent is not, in
general, idempotent: it can be a nontrivial convex combination of
Haar measures of distinct subgroups.  Thus applying the classical
Kawada--It\^o theorem separately to all vector states is not a
valid proof of \eqref{eq:intro-classification}.

\subsection*{Why the matrix coefficients must commute}
The proof isolates the noncommutative step as
Lemma~\ref{lem:tensor-leg}.  First compress to the range of the
idempotent mass $\mu(G)$ and realize the resulting normalized
positive operator-valued measure as
$\mu(D)=V^*P(D)V$, with $P$ a projection-valued measure.
For each strongly continuous unitary representation $\pi$ of $G$,
integration of $\pi$ against the dilation yields a unitary
$\mathcal U_\pi$; its compression
\[
 A_\pi=(I\otimes V)^*\mathcal U_\pi(I\otimes V)
      =\sum_{a,b}\pi(\mu_{ab})\otimes e_{ab}
\]
is idempotent by $\mu*\mu=\mu$.  As a compression of a unitary it is
contractive, and hence is an \emph{orthogonal} projection.

For two representations $\pi,\rho$, let $P_1,P_2$ denote the
operators $A_\pi,A_\rho$ acting on the two representation legs and
a common coefficient leg.  The dilated unitaries commute on their
separate representation factors, but this alone does \emph{not}
allow one to factor arbitrary compressions.  Instead we prove two
specific range identities with the projection
$R=A_{\pi\otimes\rho}$:
\[
             RP_2=P_1P_2,\qquad RP_1=P_2P_1.
\]

The above equalities imply that $P_1P_2$ is a contractive idempotent;
therefore $P_1P_2=P_2P_1$.  When both representations are regular,
these two products are integrated forms of finite matrix-valued
Radon measures on $G\times G$ whose values on rectangles are
$\mu(D)\mu(D')$ and $\mu(D')\mu(D)$.  Faithfulness of the regular
convolution representation on \emph{finite measures} shows that the
measures agree and hence that the matrix values commute.  The
Radon-product construction and a cyclic spectral decomposition avoid
separability assumptions.  This argument is the heart of the paper;
none of the harmonic results below is used in its proof.

\subsection*{Harmonic projection without tightness}
The classification theorem also resolves a structural question for
positive matrix convolution semigroups.  Let $\{\sigma_t\}_{t>0}$ be
normalized by $\sigma_t(G)=I_n$, and assume the corresponding
operators $T_t$ form a strongly continuous semigroup on
$L^2(G,M_n)$, with Hilbert--Schmidt matrix norm.  Positivity gives a
uniform variation bound, so the mean ergodic theorem produces the
orthogonal projection $P_2$ onto the common fixed space.  An abstract
translation-invariant $L^2$ projection need not be convolution by a
finite measure.  Here, however, the positive measures
\[
      \mu_R=\frac1R\int_0^R\sigma_t\,dt
\]
converge weak-$*$ in $M(G,M_n)$ to a \emph{unique} positive
idempotent $\mu$, even if some directions lose their total mass in
the limit.  Moreover (Theorem~\ref{thm:harmonic-measure-representation}),
\begin{equation}\label{eq:intro-ergodic-representation}
           P_p=T_\mu\qquad(1<p<\infty),
\end{equation}
and the same measure represents the common fixed-space projection
for every such exponent.

There is a precise geometric interpretation of the mass which survives.
For a compact subset $C\subset G$, put
$V_C=\bigcap_{t>0}\ker\sigma_t(G\setminus C)$.
Then (Theorem~\ref{thm:matrix-harmonic-compact-confinement})
\begin{equation}\label{eq:intro-confinement}
   \bigcup_{C\Subset G}V_C=\operatorname{Ran}\mu(G).
\end{equation}

Thus a nonzero finite-$p$ harmonic space exists exactly when a nonzero matrix direction remains confined to one fixed compact set for \emph{all} times.  The semigroup itself need not have commuting coefficient matrices, and other directions may escape to infinity. The classification then yields a finite decomposition of every $L^p$ harmonic space into left-invariant compact-subgroup channels. This is formulated explicitly in Theorem~\ref{thm:geometric-harmonic-space}.

The same Haar-channel description identifies every positive matrix convolution measurement on $L^2$.  We include a short bounded-generator Zeno consequence, for which the measurement is completely classified; the more extensive theory of unbounded Hunt and Zeno generators is outside the scope of this article.

\subsection*{Relation to existing work and scope}
The closest classical antecedent is the Kawada--It\^o characterization
of positive scalar convolution idempotents as Haar probabilities of
compact subgroups \cite{Heyer,KawadaIto,Kelley,Wendel54}.
This is an implication from positivity and idempotence to a
\emph{single subgroup}.  Signed or complex scalar idempotents are a
different matter, connected to the idempotent theorems of
Cohen and Host \cite{CohenIdempotent,HostIdempotents};
commutation questions for contractive scalar idempotents have also
been studied in the measure algebra \cite{SpronkCommuting}.
Mustafayev's work on scalar Kawada--It\^o convergence removes a
metrizability requirement from the corresponding classical
Ces\`aro-limit setting \cite{Mustafayev2022}.
These results neither assume nor resolve the problem of
\emph{noncommuting positive matrix coefficients} under the ordered
matrix-valued convolution of the present paper.

There is a closer connection with the matrix harmonic analysis of
Chu \cite{Chu,ChuMatrixHarmonic} and the matrix-valued
Choquet--Deny theorem of Chu, Hilberdink and Howroyd
\cite{ChuHilberdinkHowroyd}.
In the latter work, a normalized positive matrix-valued measure on a
locally compact \emph{abelian} group is the input, and the question is
when its bounded matrix-valued harmonic functions are constant.
It does not assert that a positive convolution-idempotent
matrix measure has commuting values, nor does it classify such
idempotents on arbitrary locally compact groups.
The distinction is especially important: being the measure defining
a harmonic equation is not the same hypothesis as
$\mu*\mu=\mu$, and a classification of the latter is not obtained
by applying a scalar theorem to each matrix coefficient or vector
state.  Our applications concern $L^p$ mean-ergodic projections
and the compactly confined coefficient subspaces; they do not
claim to settle the general bounded matrix Choquet--Deny problem.

A second neighboring theory concerns idempotent states on  locally compact quantum groups and quantum semigroups \cite{FranzSkalski,KasprzakKhosraviSoltan,SalmiSkalski}. There the convolution is  defined via a coproduct on an operator algebra.  Here $\mu$ is a finite $M_n$-valued \emph{Radon measure on a classical group}, positivity is required for every matrix value $\mu(D)$, and  convolution multiplies coefficient matrices in a specified order.  Neither a quantum idempotent-state theorem nor a scalar measure-algebra classification directly gives the forced commutativity  and canonical orthogonal Haar channels in \eqref{eq:intro-classification}.  The related frameworks are important comparisons, but their objects and conclusions should not be identified with ours.

The proof of the general theorem proceeds in two stages.  The commuting-value case follows from simultaneous diagonalization and the scalar classification.  The genuinely new stage is the \emph{reduction of an arbitrary positive convolution idempotent to commuting values}, using the tensor-leg projection argument. The hypothesis of positive \emph{matrix measures} is crucial: products of noncommuting positive coefficients need not remain positive, and no positivity of the convolution square is available for a general non-idempotent measure.  In particular, our result is not a classification of arbitrary complex matrix-valued idempotents.

The next section fixes the analytic and measure conventions, and
Section~\ref{sec:idempotents} proves the unrestricted classification.
Section~\ref{sec:harmonic} identifies the mean-ergodic projections,
proves compact confinement, and describes the resulting harmonic
channels.  A short bounded-generator Zeno consequence follows.
The appendix gives complementary arguments for measures with
densities, atomic parts, and transverse subgroup carriers.  Spectral
gaps, nonunimodular modular growth, general Hunt--Zeno convergence,
and quantum-Markov extensions constitute separate problems and are
not needed for any theorem in the present article.

\section{Setting and conventions}\label{sec:setting}
% ------------------------------------------------------------------

Throughout, $G$ is a locally compact group with left Haar measure $dx$. We use the modular function $\Delta:G\longrightarrow(0,\infty)$ characterized by
$$\int_G \varphi(xa)\,dx   =\Delta(a)^{-1}\int_G \varphi(x)\,dx  \qquad (a\in G,\ \varphi\in C_c(G)).$$

With this convention, the inversion formula is
$$  \int_G\varphi(x^{-1})\,dx  =\int_G\varphi(x)\Delta(x^{-1})\,dx.$$

No unimodularity assumption is imposed. We write $M_n$ for the algebra of complex $n\times n$ matrices and $I_n$ for its identity. For $1\le p<\infty$, let $\Lp=L^p(G,M_n)$ denote the space of equivalence classes of measurable functions $f:G\longrightarrow M_n$ for which
$$ \|f\|_p  =\left(\int_G\|f(x)\|_2^p\,dx\right)^{1/p}<\infty,  \qquad  \|a\|_2=\bigl(\tr(a^*a)\bigr)^{1/2}.$$

When $p=\infty$, $L^\infty(G,M_n)$ is equipped with
$$  \|f\|_\infty  =\operatorname*{ess\,sup}_{x\in G}\|f(x)\|_2.$$

Thus the matrix norm used in these spaces is the Hilbert--Schmidt norm. Although every norm on $M_n$ gives an equivalent $L^p$ norm, assertions about contractivity depend on the chosen matrix  norm. The space $L^2(G,M_n)$ is a Hilbert space with inner product
$$ \langle f,g\rangle  =\int_G\tr\bigl(f(x)g(x)^*\bigr)\,dx,$$
which is linear in the first variable. Let $M(G,M_n)$ denote the space of regular $M_n$-valued Borel measures of finite total variation. We take total variation with respect to the operator norm on $M_n$, and denote it by $|\sigma|$. A matrix measure $\sigma$ is \emph{positive} if $\sigma(E)$ is positive semidefinite for every Borel set $E\subseteq G$; in particular, $\sigma(E)^*=\sigma(E)$. Write $\check\sigma(E)=\sigma(E^{-1})$. For $\sigma,\tau\in M(G,M_n)$ define the convolution by
$$  (\sigma*\tau)(E)  =\iint_{G\times G}\mathbf1_E(yz)\,d\sigma(y)\,d\tau(z),$$
where matrix factors are multiplied in the displayed order. Equivalently, writing $\sigma=(\sigma_{ij})$ and $\tau=(\tau_{ij})$,
$$  (\sigma*\tau)_{ij}  =\sum_{k=1}^n\sigma_{ik}*\tau_{kj},$$
where the convolutions on the right are scalar-measure convolutions. For
an arbitrary locally compact group the double integral is taken with
respect to the Radon product of the scalar entries
(Lemma~\ref{lem:radon-prerequisites}(b)); equivalently,
$\int_Gf\,d(\sigma*\tau)=\iint f(yz)\,d\sigma(y)\,d\tau(z)$ for
$f\in C_0(G)$ \cite[Ch.~2]{Folland}. For second-countable $G$ this is the
ordinary product measure. For $f\in\Lp$ (and likewise for $f\in L^\infty(G,M_n)$), define \emph{left convolution} by
$$  T_\sigma f(x)=(\sigma\star f)(x)  =\int_G d\sigma(y)\,f(y^{-1}x),$$ with the matrix measure acting by multiplication on the left. Because left translations preserve left Haar measure, the variation bound gives
$$  \|T_\sigma f\|_p\le |\sigma|(G)\,\|f\|_p,  \qquad 1\le p\le\infty.$$

The bounded operators so defined satisfy
\begin{equation}\label{eq:composition}
  T_\tau T_\sigma=T_{\tau*\sigma}.
\end{equation}

Each $T_\sigma$ commutes with right translations $\rho_a f(x)=f(xa)$ and with right multiplication by constant matrices. For non-unimodular $G$, the operators $\rho_a$ on $L^p$ need not be isometries; the commutation assertion does not require that property. Positivity of the matrix \emph{measure} $\sigma$ does not in general mean that $T_\sigma$ preserves the positive cone of
$L^\infty(G)\,\bar\otimes\,M_n$.

\begin{remark}\label{rem:why-left}
Left translations preserve left Haar measure. Right-convolution norms on
non-unimodular groups involve a modular weight, whereas the left
convolution operators used here are bounded by the total variation of
the defining matrix measure.\end{remark}

\begin{definition}
A \emph{matrix convolution semigroup} is a family $\{\sigma_t\}_{t>0}$ of positive measures in $M(G,M_n)$ such that
$$  \sigma_t(G)=I_n,\qquad  \sigma_{s+t}=\sigma_s*\sigma_t\quad(s,t>0).$$

Write $T_t=T_{\sigma_t}$ for $t>0$ and $T_0=I$. Whenever the generator on $L^p(G,M_n)$ is considered, we assume that $\{T_t\}_{t\ge0}$ is strongly continuous on that space and denote its generator by $A_p$; when $p$ is fixed, we write simply $A$. The \emph{harmonic space} on that space is
$$  \hsp=\Fix(T)  =\bigcap_{t>0}\ker(T_t-I)  =\ker A_p.$$

The semigroup is \emph{symmetric} if $\check\sigma_t=\sigma_t$ for every $t>0$. For positive matrix measures this implies that $T_t$ is self-adjoint on $L^2(G,M_n)$; the converse follows from  the injectivity of the regular convolution representation.
\end{definition}

By \eqref{eq:composition} and the measure semigroup law, $T_sT_t=T_{\sigma_s*\sigma_t}=T_{s+t}$. The equality $\ker A_p=\Fix(T)$ is the standard fixed-point characterization for a strongly  continuous semigroup \cite[Ch.~II]{EngelNagel}. Strong continuity on $L^\infty$ is not assumed by default.

\begin{lemma}\label{lem:trivial}
For every matrix convolution semigroup, every $t>0$, and every
$N\in\mathbb N$,
\[
  \bigl(T_{\sigma_{t/N}}\bigr)^N
  =T_{\sigma_{t/N}^{*N}}
  =T_{\sigma_t}.
\]
\end{lemma}

\begin{proof}
This follows immediately from the measure semigroup law and
\eqref{eq:composition}.
\end{proof}

\begin{remark}\label{rem:failure}
For a convolution semigroup, $(T_{t/N})^N=T_t$ for all $N$.
A genuine Zeno product includes a fixed projection $Q$ between
consecutive time steps; the long-time ergodic limit and the Zeno limit
are distinct operations.
\end{remark}
% ------------------------------------------------------------------
\section{Elementary contraction and adjoint formulas}\label{sec:basic}
% ------------------------------------------------------------------

The first analytic observation is that normalization of a positive matrix measure gives an $L^2$ contraction, even though it does not generally yield contractivity at other exponents.
\begin{lemma}\label{lem:contraction}
Let $\sigma\in M(G,M_n)$ be positive with
$\sigma(G)=I_n$. Then
\[
\|T_\sigma f\|_2\leq \|f\|_2
\qquad(f\in L^2(G,M_n)).
\]
\end{lemma}

\begin{proof}
Let $\lambda=\operatorname{Tr}\sigma$. Since $\sigma$ is positive, there is a measurable function $W:G\longrightarrow M_n$ such that
$$W(y)\geq 0,\qquad d\sigma(y)=W(y)\,d\lambda(y),\qquad \int_G W(y)\,d\lambda(y)=I_n.$$

We first establish the matrix-valued Jensen inequality
\begin{equation}\label{eq:matrix-jensen}
\left\|\int_G W(y)F(y)\,d\lambda(y)\right\|_2^2
\leq
\int_G\operatorname{Tr}
\bigl(F(y)^*W(y)F(y)\bigr)\,d\lambda(y).
\end{equation}

Here $F$ is any measurable matrix-valued function for which the right-hand side is finite. Indeed, define $V:\mathbb C^n\longrightarrow L^2(G,\lambda;\mathbb C^n)$ by $(V\xi)(y)=W(y)^{1/2}\xi$. The normalization of $\sigma$ implies that
$$\|V\xi\|^2=\left\langle\left(\int_G W(y)\,d\lambda(y)\right)\xi,\xi\right\rangle=\|\xi\|^2.$$

Thus $V$ is an isometry. Applying the contraction $V^*$ to the columns of $W^{1/2}F$, we obtain
\[
\begin{aligned}
\left\|\int_G W(y)F(y)\,d\lambda(y)\right\|_2^2
&=
\left\|V^*(W^{1/2}F)\right\|_2^2\\
&\leq
\int_G\|W(y)^{1/2}F(y)\|_2^2
\,d\lambda(y)\\
&=
\int_G\operatorname{Tr}
\bigl(F(y)^*W(y)F(y)\bigr)
\,d\lambda(y),
\end{aligned}
\]
which proves \eqref{eq:matrix-jensen}. Now let $f\in C_c(G,M_n)$ and apply \eqref{eq:matrix-jensen} with $F(y)=f(y^{-1}x)$. We obtain
$$\|T_\sigma f(x)\|_2^2\leq\int_G\operatorname{Tr}\left(f(y^{-1}x)^*W(y)f(y^{-1}x)\right)\,d\lambda(y).$$

Integrating over $G$, using Tonelli's theorem and the left invariance of Haar measure, implies
\[
\begin{aligned}
\|T_\sigma f\|_2^2
&\leq
\int_G\int_G
\operatorname{Tr}\bigl(
f(y^{-1}x)^*W(y)f(y^{-1}x)
\bigr)\,dx\,d\lambda(y)\\
&=
\int_G\int_G
\operatorname{Tr}\bigl(
f(u)^*W(y)f(u)
\bigr)\,du\,d\lambda(y)\\
&=
\int_G\operatorname{Tr}\left(
f(u)^*
\left(\int_G W(y)\,d\lambda(y)\right)
f(u)\right)\,du\\
&=\|f\|_2^2.
\end{aligned}
\]

Since $C_c(G,M_n)$ is dense in $L^2(G,M_n)$
and convolution by a finite matrix-valued measure
is bounded on $L^2$, the inequality extends to
every $f\in L^2(G,M_n)$.
\end{proof}

\begin{lemma}\label{lem:adjoint}
Let $\sigma\in M(G,M_n)$ be positive. Then $T_\sigma^*=T_{\check\sigma}$ on $L^2(G,M_n)$. In particular, if $\check\sigma=\sigma$, then $T_\sigma$ is self-adjoint.
\end{lemma}

\begin{proof}
Write $d\sigma=W\,d\lambda$, where $\lambda=\operatorname{Tr}\sigma$ and $W(y)=W(y)^*$ for $\lambda$-almost every $y$. For $f,g\in C_c(G,M_n)$, Fubini's theorem, the substitution $x=yu$, and left invariance of Haar measure yield
\[
\begin{aligned}
\langle T_\sigma f,g\rangle
&=
\int_G\int_G
\operatorname{Tr}\bigl(
W(y)f(y^{-1}x)g(x)^*
\bigr)\,d\lambda(y)\,dx\\
&=
\int_G\int_G
\operatorname{Tr}\bigl(
W(y)f(u)g(yu)^*
\bigr)\,d\lambda(y)\,du\\
&=
\int_G
\operatorname{Tr}\left(
f(u)
\left(\int_G W(y)g(yu)\,
d\lambda(y)\right)^*
\right)\,du.
\end{aligned}
\]

By the definition of the inverse-image measure,
$$\int_G W(y)g(yu)\,d\lambda(y)=\int_G d\check\sigma(z)\,g(z^{-1}u)=(T_{\check\sigma}g)(u).$$

Thus, for all $f,g\in C_c(G,M_n)$, we obtain
$$\langle T_\sigma f,g\rangle=\langle f,T_{\check\sigma}g\rangle.$$

The above identity extends to all of $L^2(G,M_n)$ by density and boundedness of both convolution operators. Hence $T_\sigma^*=T_{\check\sigma}$.
\end{proof}

\section{Positive idempotent matrix measures}\label{sec:idempotents}
% ------------------------------------------------------------------

The scalar Kawada--It\^o--Kelley theorem identifies a nonzero positive convolution-idempotent scalar measure with normalized Haar measure on a compact subgroup; see \cite{Heyer,KawadaIto,Kelley,Loynes,Pym,Salmi,Wendel54} and the references therein. The matrix-valued analogue does not force a single compact subgroup. The following results give precise alternatives and describe the admissible convolution measurements used in the bounded Zeno application below. Throughout this section positivity means positivity of the matrix $\mu(A)$ for every Borel set $A$; it is stronger than positivity of an associated convolution operator.

\begin{example}\label{ex:positive-convolution-obstruction}
Let $G$ be the free group on three generators $g_1,g_2,g_3$. In $\mathbb C^2$ set
\[
  v_1=\binom{1}{0},\qquad
  v_2=\binom{-1/2}{\sqrt3/2},\qquad
  v_3=\binom{-1/2}{-\sqrt3/2},\qquad
  A_j=\frac23 v_jv_j^*.
\]

The matrices $A_j$ are positive and $A_1+A_2+A_3=I_2$. Thus, $\nu=\sum_{j=1}^3 A_j\delta_{g_j}$ is a normalized positive matrix measure.  Since $g_1g_2$ is the product of exactly one ordered pair of the listed generators, however,
\[
 (\nu*\nu)(\{g_1g_2\})=A_1A_2
  =\begin{pmatrix}1/9&-\sqrt3/9\\0&0\end{pmatrix},
\]
which is not self-adjoint.  Thus positivity of $\nu$ and $\nu(G)=I_2$ do not ensure positivity of $\nu*\nu$. The classification theorem concerns precisely the rigid case in
which the convolution square is again the original positive measure.
\end{example}

\begin{proposition}\label{prop:idempotent-mass}
Let $\mu\in M(G,M_n)$ be a positive convolution-idempotent matrix measure. Then $P:=\mu(G)=P^*=P^2$, $\mu(A)=P\mu(A)P$, for every Borel subset $A\subseteq G$. Thus either $\mu=0$ (when $P=0$), or $\mu$ reduces canonically to a positive convolution-idempotent measure on $P\C^n$ whose total mass is the identity of that space. In particular, a positive convolution idempotent cannot have a nonorthogonal total-mass projection.
\end{proposition}
\begin{proof}
The total-mass map is multiplicative for matrix convolution, so
\[
 P^2=(\mu*\mu)(G)=\mu(G)=P.
\]

Since $\mu$ is positive, $P\ge0$, hence $P=P^*$ and $P$ is an orthogonal projection. For every Borel set $A$, we obtain
\[
 0\le\mu(A)\le P,
 \qquad P-\mu(A)=\mu(G\setminus A).
\]

If $v\in\ker P$, then $0\le\langle\mu(A)v,v\rangle\le\langle Pv,v\rangle=0$. Since $\mu(A)\ge0$, this implies $\mu(A)^{1/2}v=0$ and hence $\mu(A)v=0$. Taking adjoints gives
\[
 \mu(A)(I-P)=0=(I-P)\mu(A),
\]
which proves $\mu(A)=P\mu(A)P$. The compressed measure on
$P\C^n$ has total mass $I_{P\C^n}$ and retains positivity and
convolution idempotence. If $P=0$, the preceding order inequality
forces $\mu(A)=0$ for every $A$.
\end{proof}

\begin{corollary}\label{cor:matrix-idempotent-l2}
For every positive convolution-idempotent $\mu\in M(G,M_n)$, $T_\mu$ is an orthogonal projection on $L^2(G,M_n)$. Moreover, $\check\mu=\mu$, equivalently $\mu(A^{-1})=\mu(A)$, for every Borel subset $A\subseteq G$. Neither assertion requires an a priori commutativity assumption on the matrix values $\mu(A)$.
\end{corollary}
\begin{proof}
Put $P=\mu(G)$. By Proposition~\ref{prop:idempotent-mass}, we have $\mu(A)=P\mu(A)P$. For $P\ne0$, apply the matrix-valued Jensen inequality from Lemma~\ref{lem:contraction} in the corner $P M_nP$, where the total mass is the identity on $P\C^n$. Apply this inequality separately to each column of $Pf$ (the columns lie in $P\C^n$); it remains valid for the rectangular corner $PM_n$.
Since $T_\mu f=T_\mu(Pf)$, it follows that
\[
 \|T_\mu f\|_{L^2(G,M_n)}
 \le\|Pf\|_{L^2(G,M_n)}
 \le\|f\|_{L^2(G,M_n)}.
\]

The same estimate is immediate for $P=0$. The composition identity \eqref{eq:composition} implies that $T_\mu^2=T_{\mu*\mu}=T_\mu$. A contractive idempotent on a Hilbert space is an orthogonal
projection: for $x$ in its range and $y$ in its kernel, $\|x\|=\|T_\mu(x+zy)\|\le\|x+zy\|$, for all $z\in\C$, which implies $\langle x,y\rangle=0$. Therefore $T_\mu=T_\mu^*$. By Lemma~\ref{lem:adjoint}, positivity of $\mu$ implies $T_\mu^*=T_{\check\mu}$. Faithfulness of convolution on $L^2$, proved in Lemma~\ref{lem:faithful-measure-convolution}, now implies $\check\mu=\mu$.
\end{proof}

\begin{example}\label{ex:idempotent-two-subgroups}
Let $G=\{e,s\}\cong\mathbb Z/2\mathbb Z$ and
$m_G=\tfrac12(\delta_e+\delta_s)$. Set
\[
 \mu=\begin{pmatrix}\delta_e&0\\0&m_G\end{pmatrix},
 \qquad
 \mu(\{e\})=\begin{pmatrix}1&0\\0&1/2\end{pmatrix},
 \qquad
 \mu(\{s\})=\begin{pmatrix}0&0\\0&1/2\end{pmatrix}.
\]

Then $\mu\ge0$, $\mu(G)=I_2$, and $\mu*\mu=\mu$, since $\delta_e*\delta_e=\delta_e$, $m_G*m_G=m_G$, and the two channel projections are orthogonal. More explicitly, we have
\[
 (T_\mu f)(x)=e_{11}f(x)
       +e_{22}\,\frac{f(e)+f(s)}2,
 \qquad x\in G.
\]

Thus $T_\mu$ fixes the first row and projects the second row onto constant functions, providing a concrete orthogonal $L^2$-projection. It cannot equal $m_K\otimes E$ for any subgroup $K\le G$ and any projection $E$: total mass would force $E=I_2$, whereas $\mu(\{s\})$ has rank one and is nonzero, while every nonzero atom of $m_K\otimes I_2$ has rank two. The rank obstruction persists under
constant unitary conjugation and translation of the measure; ordinary translation need not preserve convolution idempotence. Hence a
single-subgroup classification is false even for finite abelian $G$.
\end{example}

\begin{remark}
The orthogonal-projection conclusion in Corollary~\ref{cor:matrix-idempotent-l2} does not assert that $T_\mu$ preserves the pointwise positive cone of matrix-valued
functions. Indeed, in Example~\ref{ex:idempotent-two-subgroups}, let
\[
 f(e)=\begin{pmatrix}1&1\\1&1\end{pmatrix},
 \qquad f(s)=0.
\]

Both values of $f$ are positive semidefinite, whereas
\[
 (T_\mu f)(e)=
 \begin{pmatrix}1&1\\1/2&1/2\end{pmatrix}
\]
is not self-adjoint. Thus ``measurement'' refers here to an orthogonal projection in the Hilbert--Schmidt Hilbert space, not to a positive or completely positive map on the algebra of matrix-valued functions.
\end{remark}

\begin{theorem}\label{thm:commuting-idempotents}
Let $G$ be a locally compact group and let $\mu\in M(G,M_n)$ be a positive convolution-idempotent measure. Suppose that $\mu(A)\mu(B)=\mu(B)\mu(A)$ for all Borel sets $A,B\subseteq G$.
Then there exist {pairwise distinct} compact subgroups $K_1,\ldots,K_r\le G$ and mutually orthogonal nonzero projections $E_1,\ldots,E_r\in M_n$ such that
\begin{equation}\label{eq:orthogonal-haar-decomp}
 \mu=\sum_{j=1}^{r}m_{K_j}\otimes E_j,
 \qquad \sum_{j=1}^{r}E_j=\mu(G),\qquad r\le n.
\end{equation}

After equal subgroups have been combined, the pairs $(K_j,E_j)$ are uniquely determined up to permutation. Moreover,
\[
 \operatorname{supp}(\operatorname{Tr}\mu)=\bigcup_{j=1}^r K_j,
\]
so every such measure has compact support, even when $G$ is noncompact.
Conversely, any finite sum of normalized Haar measures of compact
subgroups with mutually orthogonal coefficient projections is a positive
convolution idempotent. For $\mu=0$ the decomposition is the empty sum.
\end{theorem}
\begin{proof}
Every $\mu(A)$ is positive and hence Hermitian. The commuting family $\{\mu(A): A\subseteq G\text{ Borel}\}$ can be simultaneously unitarily diagonalized, even though the family may be infinite, since it acts on the finite-dimensional space $\C^n$. Choose a unitary $U$ with
\[
 U^*\mu(A)U=\operatorname{diag}(\nu_1(A),\ldots,\nu_n(A))
 \qquad(A\subseteq G\text{ Borel}).
\]

The $\nu_i$ are finite positive regular scalar measures. Constant unitary conjugation respects matrix convolution is as follows
\[
 U^*(\mu*\mu)U=(U^*\mu U)*(U^*\mu U).
\]

The right-hand side is $\operatorname{diag}(\nu_1*\nu_1,\ldots,\nu_n*\nu_n)$, so idempotence of $\mu$ implies $\nu_i*\nu_i=\nu_i$ for each $i$.
Taking total masses implies  $\nu_i(G)^2=\nu_i(G)$, hence $\nu_i(G)\in\{0,1\}$. A positive measure of zero total mass is zero. If $\nu_i(G)=1$, then $\nu_i$ is an idempotent probability measure on the locally compact group $G$. The scalar Kawada--It\^o--Kelley theorem implies that $\nu_i=m_{H_i}$ for a compact subgroup $H_i\le G$; see \cite{Heyer,KasprzakKhosraviSoltan,KawadaIto,Kelley,Loynes,Pym}, for more details. No countability assumption on $G$ is needed for it. List the distinct subgroups among the $H_i$ as $K_1,\ldots,K_r$, and let $E_j$ be the projection onto $\operatorname{span}\{Ue_i:\nu_i=m_{K_j}\}$. Then the $E_j$ are nonzero and mutually orthogonal, and
\[
 \mu(A)=\sum_{j=1}^r m_{K_j}(A)E_j,
 \qquad \sum_{j=1}^rE_j=\mu(G).
\]

In particular, $r\le\operatorname{rank}\mu(G)\le n$. Taking the matrix trace gives the positive scalar measure
\[
 \operatorname{Tr}\mu=\sum_{j=1}^r
       (\operatorname{rank}E_j)\,m_{K_j}.
\]

Every normalized Haar measure $m_{K_j}$ has support precisely $K_j$. Since this is a finite sum with strictly positive coefficients, $\operatorname{supp}(\operatorname{Tr}\mu)=\bigcup_jK_j$, a compact set. For uniqueness, assume that $\mu\ne0$ and regard $\mathcal A=C^*(\mu(A):A\subseteq G\text{ Borel})$
as a unital algebra on $P\C^n$, where $P=\mu(G)$ is its identity.
The preceding decomposition gives
$\mathcal A\subseteq\bigoplus_{j=1}^r\C E_j$.
For $i\ne j$, the measures $m_{K_i}$ and $m_{K_j}$ are distinct:
they have different supports. Consequently, some Borel set $A$
satisfies $m_{K_i}(A)\ne m_{K_j}(A)$.
Thus the coordinate functions supplied by the matrices $\mu(A)$
separate the $r$ summands. The elementary finite-dimensional
Stone--Weierstrass argument (or polynomial interpolation) implies
\[
 \mathcal A=\bigoplus_{j=1}^r\C E_j.
\]

It follows that the $E_j$ are precisely the minimal nonzero projections of the intrinsically defined abelian algebra $\mathcal A$. They are therefore determined by $\mu$ up to permutation. Finally, 
\[
 m_{K_j}(A)=\frac{\operatorname{Tr}(E_j\mu(A))}
                       {\operatorname{rank}E_j}
\]
recovers the corresponding probability measure, and $K_j=\operatorname{supp}(m_{K_j})$ recovers its subgroup. We now consider the converse case. For any compact subgroups $K_j$ and mutually orthogonal projections $E_j$, the measure $\mu=\sum_jm_{K_j}\otimes E_j$ is positive. By bilinearity of convolution,
\[
 \mu*\mu
 =\sum_{i,j}(m_{K_i}*m_{K_j})\otimes E_iE_j
 =\sum_j(m_{K_j}*m_{K_j})\otimes E_j
 =\sum_jm_{K_j}\otimes E_j=\mu,
\]
since $E_iE_j=0$ for $i\ne j$, $E_j^2=E_j$, and $m_K*m_K=m_K$ for every compact subgroup $K$. The empty sum gives the zero measure.
\end{proof}

\begin{remark}\label{rem:commuting-haar-intrinsic}
In Theorem~\ref{thm:commuting-idempotents}, the channel projections $E_j$ are intrinsic, rather than artifacts of the diagonalizing unitary. On $L^2(G,M_n)$ the corresponding convolution operator  is 
\[
 (T_\mu f)(x)=\sum_{j=1}^r E_j
       \int_{K_j}f(k^{-1}x)\,dm_{K_j}(k).
\]

Each summand is the orthogonal subgroup-averaging projection on its row corner, and distinct summands have orthogonal ranges. Thus the commuting-value case also gives a direct explanation of Corollary~\ref{cor:matrix-idempotent-l2}. 
\end{remark}

For abelian groups, the commuting-value hypothesis in Theorem~\ref{thm:commuting-idempotents} is automatic. The following positive-definite-function lemma is actually valid for every group; its  application to matrix-valued Fourier--Stieltjes transforms requires the locally compact abelian setting. Matrix-valued positive type and positive definite functions attached to the convolution  algebra $L^1(G,M_n)$, together with a GNS construction and a Bochner theorem on abelian groups, are studied in \cite{EbadianJabbariBochner,JabbariPositiveType}. There positivity is tested through  the trace pairing with $M_n$-valued test functions; here we use the block condition $\sum_{i,j}v_i^*F(h_i^{-1}h_j)v_j\ge0$ ($v_i\in\C^r$), which comes from an operator-valued GNS dilation.

\begin{lemma}\label{lem:pd-projections-abelian}
Let $H$ be an arbitrary group, and let $F:H\longrightarrow M_r$ be a normalized matrix-valued positive-definite function, $F(e)=I_r$. If every $F(h)$ is an orthogonal projection, then
\[
 F(g)F(h)=F(h)F(g)\qquad(g,h\in H).
\]

In particular, commutativity of $H$ is not required.
\end{lemma}
\begin{proof}
By the operator-valued GNS construction, there are a Hilbert space $\mathcal H$, a unitary representation $U:H\longrightarrow\mathcal U(\mathcal H)$, and an isometry $V:\C^r\longrightarrow\mathcal  H$ such that $F(h)=V^*U(h)V$. If $v\in\Ran F(h)$, then $F(h)v=v$, and hence
\[
 \langle U(h)Vv,Vv\rangle
   =\langle F(h)v,v\rangle=\|v\|^2
   =\|U(h)Vv\|\,\|Vv\|.
\]

Equality in Cauchy--Schwarz, with the displayed positive real inner product, implies $U(h)Vv=Vv$. Thus
$U(h)VF(h)=VF(h)$, and consequently
\begin{equation}\label{eq:pd-projection-identity}
 F(kh)F(h)=F(k)F(h)\qquad(k,h\in H).
\end{equation}

Positive definiteness also implies $F(h^{-1})=F(h)^*=F(h)$. Fix $g,h\in H$ and put 
\[
 P=F(g),\qquad S=F(h),\qquad
 R=F(hg^{-1})=F(gh^{-1}),
\]
where the last equality uses inversion symmetry, not commutativity of the group. Applying \eqref{eq:pd-projection-identity} first with $(k,h)=(hg^{-1},g)$ and then with $(k,h)=(gh^{-1},h)$ implies $SP=RP$, and $PS=RS$. Taking adjoints in the first identity yields $PS=PR$. Therefore
\[
 (PS)^2=(PR)(RS)=PR^2S=PRS=P(PS)=PS.
\]

The product $PS$ is a contraction, since $P$ and $S$ are orthogonal projections. A contractive idempotent on a Hilbert space is an orthogonal projection, hence $PS=(PS)^*=SP$.
\end{proof}

\begin{corollary}\label{cor:pd-subgroup-indicators}
Under the hypotheses of Lemma~\ref{lem:pd-projections-abelian}, there are pairwise distinct subgroups $H_1,\ldots,H_m\le H$ and mutually orthogonal nonzero projections $E_1,\ldots,E_m\in M_r$
with $\sum_jE_j=I_r$ such that
\begin{equation}\label{eq:pd-subgroup-indicators}
 F(h)=\sum_{j=1}^m\mathbf 1_{H_j}(h)E_j\qquad(h\in H).
\end{equation}

The decomposition with distinct subgroups is unique up to permutation. Conversely, every such sum is normalized and matrix-valued positive definite and has projection values.
\end{corollary}
\begin{proof}
By the Lemma \ref{lem:pd-projections-abelian} the values of $F$ are commuting Hermitian matrices, so a single unitary simultaneously diagonalizes them. Thus
$F(h)=\sum_{i=1}^r\varphi_i(h)P_i$, where the $P_i$ are pairwise orthogonal rank-one projections with sum $I_r$ and $\varphi_i(h)\in\{0,1\}$. Each $\varphi_i$ is normalized positive definite. Let  $\xi_i$ denote the corresponding unit vector in $\C^r$, so that $P_i\xi_i=\xi_i$. In the GNS dilation above, 
\[
 \varphi_i(h)=1\quad\Longleftrightarrow\quad
 U(h)V\xi_i=V\xi_i.
\]

Therefore $H_i=\{h:\varphi_i(h)=1\}$ is the stabilizer of $V\xi_i$, hence a subgroup, and $\varphi_i=\mathbf 1_{H_i}$. Grouping equal subgroups gives \eqref{eq:pd-subgroup-indicators}. For  uniqueness, the matrices $F(h)$ generate a finite-dimensional commutative unital $C^*$-algebra. Distinct indicator functions $\mathbf 1_{H_j}$ separate the $m$ coordinates, so the minimal projections of that algebra are exactly the grouped $E_j$. They are intrinsic, and their scalar coefficient functions recover $H_j$ uniquely.

Conversely, for any subgroup $K\le H$, the function $\mathbf 1_K$ is positive definite: it is the matrix coefficient of the quasi-regular representation on $\ell^2(H/K)$ at the unit vector  $\delta_K$. A finite orthogonal sum of these scalar positive-definite functions with coefficients $E_j$ is matrix-valued positive definite. Its value at the identity is $I_r$, and orthogonality of the $E_j$ makes all values projections.
\end{proof}

\begin{remark}
If $H$ is a topological group and $F$ is continuous, then each subgroup $H_j$ in \eqref{eq:pd-subgroup-indicators} is clopen, since its indicator function is continuous. Without a continuity
hypothesis the subgroups need not be closed.
\end{remark}

\begin{theorem}\label{thm:lca-matrix-idempotents}
Let $G$ be a locally compact abelian group and let $\mu\in M(G,M_n)$ be positive with $\mu*\mu=\mu$. Writing $P=\mu(G)$, there exist pairwise distinct compact subgroups $K_1,\ldots,K_r\le G$ and  mutually orthogonal nonzero projections $E_1,\ldots,E_r\in M_n$ such that
\[
 \mu=\sum_{j=1}^r m_{K_j}\otimes E_j,
 \qquad \sum_{j=1}^rE_j=P,
 \qquad r\le\operatorname{rank}P.
\]

The decomposition with distinct subgroups is unique up to permutation. Equivalently, for every $\gamma\in\widehat G$,
\begin{equation}\label{eq:lca-fourier-haar}
 \widehat\mu(\gamma)=\sum_{j=1}^r
        \mathbf1_{K_j^\perp}(\gamma)E_j,
 \quad K_j^\perp=\{\gamma\in\widehat G:\gamma(k)=1
                                  \text{ for all }k\in K_j\}.
\end{equation}

In particular $\operatorname{supp}(\operatorname{Tr}\mu) =\bigcup_{j=1}^rK_j$ is compact.  The zero measure corresponds  to the empty sum. Conversely, every such sum is a positive convolution
idempotent.
\end{theorem}
\begin{proof}
Proposition~\ref{prop:idempotent-mass} gives $P=P^2=P^*$ and $\mu(A)=P\mu(A)P$ for all Borel sets $A$.  If $P=0$, then $\mu=0$ and there is nothing to prove. Otherwise, identify $P\C^n$ with $\C^d$, $d=\operatorname{rank}P$, and work in that corner, so $\mu(G)=I_d$. For $\gamma\in\widehat G$, define
\[
 F(\gamma)=\widehat\mu(\gamma)
      =\int_G\overline{\gamma(x)}\,d\mu(x).
\]

For characters $\gamma_1,\dots,\gamma_N$ and vectors $v_1,\dots,v_N\in\C^d$, positivity of the matrix measure implies 
\[
 \begin{aligned}
 &\sum_{i,j=1}^N v_i^*F(\gamma_i^{-1}\gamma_j)v_j=&\int_G\Bigl(\sum_i\overline{\gamma_i(x)}v_i\Bigr)^*
       d\mu(x)
       \Bigl(\sum_j\overline{\gamma_j(x)}v_j\Bigr)\ge0.
 \end{aligned}
\]

Thus $F$ is a normalized matrix-valued positive-definite function
on $\widehat G$, with $F(1)=I_d$.  The positive-definite
$2\times2$ block matrix
\[
 \begin{pmatrix}I_d&F(\gamma)\\F(\gamma)^*&I_d\end{pmatrix}\ge0
\]
implies $\|F(\gamma)\|_{\mathrm{op}}\le1$. Since $\mu$ is a finite regular matrix measure, its Fourier-Stieltjes transform $F$ is norm-continuous on $\widehat G$: characters converge uniformly on  compact subsets of $G$, and the variation of $\mu$ has arbitrarily small mass outside a compact set. Because $G$ is abelian, the Fourier--Stieltjes transform takes matrix convolution to matrix multiplication, and hence
\[
 F(\gamma)^2=\widehat{\mu*\mu}(\gamma)=F(\gamma).
\]

Each $F(\gamma)$ is therefore a contractive idempotent on $\C^d$, hence an orthogonal projection. Applying Corollary~\ref{cor:pd-subgroup-indicators} to $F$ on $\widehat G$ gives pairwise distinct  subgroups $H_1,\dots,H_r\le\widehat G$ and mutually orthogonal nonzero projections $E_1,\dots,E_r$ with
$\sum_jE_j=I_d$ and
\[
 F(\gamma)=\sum_{j=1}^r\mathbf1_{H_j}(\gamma)E_j.
\]

For each $j$, the coefficient
\[
 \mathbf1_{H_j}(\gamma)
   =\frac{\operatorname{Tr}(E_jF(\gamma))}
                  {\operatorname{rank}E_j}
\]
is continuous, so $H_j$ is an open (and closed) subgroup of $\widehat G$. By Pontryagin duality for locally compact abelian groups \cite{Folland,HewittRoss}, the annihilator 
\[
 K_j=H_j^\perp
  =\{x\in G:\gamma(x)=1\text{ for all }\gamma\in H_j\}
\]
is a compact subgroup of $G$ and $K_j^\perp=H_j$. The normalized Haar measure on $K_j$ has Fourier--Stieltjes transform 
\[
 \widehat{m_{K_j}}(\gamma)
     =\int_{K_j}\overline{\gamma(x)}\,dm_{K_j}(x)
     =\mathbf1_{K_j^\perp}(\gamma).
\]

Thus $F$ is the transform of $\sum_jm_{K_j}\otimes E_j$. Uniqueness of the Fourier--Stieltjes transform of finite complex regular measures, applied to each matrix entry, yields the claimed identity for $\mu$.  Returning to $M_n$, the $E_j$ act on $P\C^n$ and are extended by zero on $(I-P)\C^n$. The grouped projections $E_j$ and distinct subgroups $H_j$ are unique up to permutation  by Corollary~\ref{cor:pd-subgroup-indicators}; $K_j=H_j^\perp$ then makes the Haar-channel decomposition unique. Since
\[
 \operatorname{Tr}\mu
  =\sum_j(\operatorname{rank}E_j)m_{K_j},
\]
its support is $\bigcup_jK_j$, a finite union of compact sets. Finally, $E_iE_j=0$ for $i\ne j$ and $m_{K_j}*m_{K_j}=m_{K_j}$, so any such sum is positive and convolution-idempotent, as already  verified in Theorem~\ref{thm:commuting-idempotents}.
\end{proof}

\begin{remark}[When the group has no nontrivial compact subgroup]
\label{rem:lca-no-compact-subgroups}
If $G$ has no nontrivial compact subgroup, the only possible Haar
channel in Theorem~\ref{thm:lca-matrix-idempotents} is $\delta_e$.
Consequently every positive convolution idempotent has the form
$\mu=P\delta_e$, with $P\in M_n$ an orthogonal projection.
This applies, in particular, when $\widehat G$ is connected
(e.g. $G=\R^d$): a continuous subgroup indicator on connected
$\widehat G$ must be constant, and its value at the identity is $1$.
\end{remark}

\subsection{The matrix Kawada--It\^o theorem on arbitrary locally compact groups}
\label{subsec:general-kawada}

The abelian proof rests on two ingredients: a Naimark dilation turns each Fourier coefficient into a compression of a unitary, and idempotence turns that compression into a projection. Both  ingredients survive for arbitrary unitary representations. The group structure of $\widehat G$ is then replaced by inner tensor products, and the final step invokes faithfulness of the integrated  left regular representation on finite measures over $G\times G$.

\begin{lemma}\label{lem:faithful-measure-convolution}
Let $G$ be a locally compact group and $n\ge1$. For each $1\le p<\infty$, the bounded convolution representation 
\[
 M(G,M_n)\longrightarrow\mathcal B(L^p(G,M_n)),\qquad
 \nu\longmapsto T_\nu,
\]
is injective. The same conclusion holds for the standard weak-* convolution representation on $L^\infty(G,M_n)$. In particular, for any locally compact group $H$ and any finite complex Radon  measure $\nu$ on $H$,
\[
 \lambda_H(\nu):=\int_H\lambda_H(y)\,d\nu(y)=0
 \quad\text{on }L^2(H)
 \quad\Longrightarrow\quad\nu=0.
\]
%
%No second-countability, separability, or $\sigma$-compactness assumption is required.
\end{lemma}
\begin{proof}
Left translation is an isometry on each $L^p(G,M_n)$, so $\|T_\nu\|_{p\longrightarrow p}\le |\nu|(G)$. The same bound holds on $L^\infty$ for the weak-$*$ integral of translations. Suppose that $T_\nu=0$ on one of these spaces. For $\varphi\in C_c(G)$, put $f_\varphi(x)=\varphi(x)I_n$; this is a valid test function for every $1\le p\le\infty$. Then
\[
 F_\varphi(x):=(T_\nu f_\varphi)(x)
       =\int_G\varphi(y^{-1}x)\,d\nu(y)
\]
is the canonical continuous representative of the convolution. Indeed, compactly supported continuous functions on a locally compact group are uniformly continuous for right translations, and, for $z\longrightarrow e$,
\[
 \sup_{x\in G}\|F_\varphi(xz)-F_\varphi(x)\|
 \le |\nu|(G)\sup_{u\in G}|\varphi(uz)-\varphi(u)|\longrightarrow0.
\]

As $T_\nu f_\varphi=0$ in $L^p$ (or $L^\infty$), the continuous function $F_\varphi$ vanishes Haar-almost everywhere. Every nonempty open set has positive Haar measure, so $F_\varphi$ vanishes  everywhere. In particular,
\[
 0=F_\varphi(e)
   =\int_G\varphi(y^{-1})\,d\nu(y).
\]

Since inversion is a homeomorphism and $\{y\mapsto\varphi(y^{-1}):\varphi\in C_c(G)\}=C_c(G)$, every matrix entry of $\nu$ annihilates $C_c(G)$. The uniqueness part of the Riesz-Markov representation theorem yields $\nu=0$. For $n=1$ and $p=2$, the operator $T_\nu$ is the integrated left regular representation $\lambda_H(\nu)$, giving the last assertion.
\end{proof}

\begin{remark}\label{rem:faithfulness-measures-vs-full}
The faithfulness established above concerns the measure algebra $M(G,M_n)$ (and, in particular, finite complex measures on $G\times G$). It should not be confused with faithfulness of the  canonical map from the full to the reduced group $C^*$-algebra, which can fail for nonamenable groups. In applying the lemma to the tensor-leg argument, one must first identify the two expressions  as integrated operators of finite matrix-valued Radon measures on $G\times G$; only then does faithfulness imply equality of those measures.
\end{remark}

\begin{theorem}[Matrix Kawada--It\^o theorem]
\label{thm:general-matrix-kawada}
Let $G$ be any locally compact group and let $\mu\in M(G,M_n)$ be
positive with $\mu*\mu=\mu$. Then the values of $\mu$ commute, and there are
compact subgroups $K_1,\ldots,K_r$ and pairwise orthogonal nonzero
projections $E_1,\ldots,E_r$ with
\[
 \mu=\sum_{j=1}^rm_{K_j}\otimes E_j,\qquad\sum_jE_j=\mu(G).
\]

If the subgroups $K_1,\ldots,K_r$ are taken pairwise distinct, the decomposition is unique up to the order of the summands, and
\begin{gather*}
 \Ran E_j=\{u\in\C^n:\ \mu(D)u=m_{K_j}(D)u\ \text{for all Borel }D\}
 \subseteq\Ran\mu(G),\\
 K_j=\{h\in G:\ \delta_h*(\mu u)=\mu u\}\qquad(0\ne u\in\Ran E_j).
\end{gather*}

On the complement, $\ker\mu(G)=\bigcap_j\ker E_j$, every value of $\mu$ vanishes (Proposition~\ref{prop:idempotent-mass}); this summand carries no subgroup, and it is the whole space exactly when $\mu=0$.
\end{theorem}

The decisive noncommutative step is isolated in the following lemma. It concerns an arbitrary projection-valued dilation of $\mu$, and its proof uses only the idempotence of $\mu$.  For a Hilbert space $\mathcal H$, a \emph{regular projection-valued measure} $P$ on $G$ is a map from Borel sets to orthogonal projections on $\mathcal H$ with $P(G)=I$, $P(D\cap D')=P(D)P(D')$, such that every  $\nu_{\eta,\eta'}=\langle P(\cdot)\eta,\eta'\rangle$ is a regular complex Borel measure. Matrices act on $\C^n$ with standard basis $(e_a)$; we write $e_{ab}=e_ae_b^*$ and  $\mu_{ab}=\langle\mu(\cdot)e_b,e_a\rangle$, so that $(\mu*\nu)_{ab}=\sum_c\mu_{ac}*\nu_{cb}$. For a strongly continuous unitary representation $\pi$ of $G$ and a finite complex measure $\nu$, $\pi(\nu)=\int\pi(x)\,d\nu(x)$ is the integrated representation, defined weakly. Recall that $\pi(\nu*\nu')=\pi(\nu)\pi(\nu')$ \cite[Ch.~3]{Folland}. Tensor factors are always listed in the order  written. For three Hilbert spaces and an operator $X$ on $\mathcal K_i\otimes\mathcal K_j$, $X^{(ij)}$ denotes the operator on $\mathcal K_1\otimes\mathcal K_2\otimes\mathcal K_3$ acting as $X$ on factors $i,j$ and as the identity on the remaining factor.

\begin{lemma}\label{lem:radon-prerequisites}
Let $G$ be a locally compact group, with no countability assumption.
\begin{enumerate}
\item[(a)] Suppose that $\nu$ is a finite positive Radon measure on $G$ and that $\pi:G\longrightarrow\mathcal U(\mathcal H_\pi)$ is a strongly continuous
unitary representation, where $\mathcal H_\pi$ need not be separable. For every $F\in L^2(G,\nu;\mathcal H_\pi)$, the field $x\mapsto\pi(x)F(x)$ is strongly $\nu$-measurable and belongs to the
same Bochner space. Moreover,
\[
 (M_\pi F)(x)=\pi(x)F(x)
\]
defines a unitary operator on $L^2(G,\nu;\mathcal H_\pi)$ with inverse
$(M_\pi^{-1}F)(x)=\pi(x)^*F(x)$.
\item[(b)] For any two finite complex Radon measures $\nu,\nu'$ on
$G$, there exists a unique finite complex Radon measure
$\nu\widehat\otimes\nu'$ on $G\times G$ satisfying
\begin{equation}\label{eq:radon-product-tensors}
 \int_{G\times G} f(x)g(y)\,d(\nu\widehat\otimes\nu')(x,y)
 =\left(\int_G f\,d\nu\right)\left(\int_G g\,d\nu'\right)
 \qquad(f,g\in C_c(G)).
\end{equation}
It agrees with the ordinary product on Borel rectangles:
$(\nu\widehat\otimes\nu')(A\times B)=\nu(A)\nu'(B)$.
In addition, the linear span of
$\{f\otimes g:f,g\in C_c(G)\}$ is uniformly dense in
$C_0(G\times G)$, so a finite complex Radon measure on $G\times G$
is determined by its integrals against these functions.
\end{enumerate}

The Radon product is used in (b) because, for an arbitrary locally compact $G$, the Borel $\sigma$-algebra of $G\times G$ need not equal $\mathcal B(G)\otimes\mathcal B(G)$.
\end{lemma}

\begin{proof}
(a) By inner regularity and finiteness of $\nu$, choose compact sets $K_m\subseteq G$ with $\nu(G\setminus\bigcup_m K_m)=0$. If $F=\sum_{j=1}^N\mathbf 1_{D_j}\xi_j$ is simple, then
\[
 x\longmapsto\pi(x)F(x)
   =\sum_{j=1}^N\mathbf 1_{D_j}(x)\pi(x)\xi_j
\]
is weakly measurable. Its range, outside a $\nu$-null set, lies in the finite union of the sets $\bigcup_m\pi(K_m)\xi_j$. Each $\pi(K_m)\xi_j$ is compact in the norm topology of $\mathcal H_\pi$,  by strong continuity, and therefore separable. Thus the field is essentially separably valued, and Pettis' measurability theorem makes it strongly measurable. For arbitrary $F$ choose simple $F_j$  converging to $F$ pointwise $\nu$-almost everywhere. Since $\pi(x)$ is unitary, $\pi(x)F_j(x)\longrightarrow\pi(x)F(x)$ almost everywhere; consequently the latter is strongly measurable. Pointwise unitarity implies
\[
 \int_G\|\pi(x)F(x)\|^2\,d\nu(x)
 =\int_G\|F(x)\|^2\,d\nu(x).
\]

The same argument applies to the strongly continuous map $x\mapsto\pi(x)^*=\pi(x^{-1})$, and the corresponding multiplier is the two-sided inverse of $M_\pi$.

(b) Existence of the finite Radon product, including agreement with the product on Borel rectangles, follows from the Radon product construction (see \cite[\S7.4]{FollandRA}); for complex measures  one may use linear combinations of finite positive Radon measures. For completeness, existence can also be seen directly from the Riesz representation theorem: for $\varphi\in C_c(G\times G)$ define
\[
 \Lambda(\varphi)
  =\int_G\left(\int_G\varphi(x,y)\,d\nu(x)\right)d\nu'(y).
\]

The inner integral is a continuous compactly supported function of $y$, because $\operatorname{supp}\varphi$ is compact. Moreover, $|\Lambda(\varphi)|\le|\nu|(G)|\nu'|(G)\|\varphi\|_\infty$. The  Riesz theorem therefore supplies a unique finite complex Radon measure representing $\Lambda$, and its action on tensors is \eqref{eq:radon-product-tensors}.

We spell out why tensor tests determine arbitrary finite Radon measures on the product. Let $\kappa$ annihilate every $f\otimes g$ with $f,g\in C_c(G)$, and let $\varphi\in C_c(G\times G)$. Choose  a compact $K\subseteq G$ such that $\operatorname{supp}\varphi\subseteq K\times K$, and choose $\chi\in C_c(G)$ with $0\le\chi\le1$ and $\chi=1$ on $K$. Set $L=\operatorname{supp}\chi$. The  restrictions to $L\times L$ of finite sums of $C_c(G)$-tensors form a self-adjoint algebra that contains constants and separates points. By Stone--Weierstrass, for each $\varepsilon>0$ there is  $\psi=\sum_{j=1}^N f_j\otimes g_j$, with $f_j,g_j\in C_c(G)$, such that $\sup_{L\times L}|\varphi-\psi|<\varepsilon$. Because $\varphi=(\chi\otimes\chi)\varphi$ everywhere, the function $\widetilde\psi=(\chi\otimes\chi)\psi$ is still a finite sum of $C_c(G)$-tensors and satisfies $\|\varphi-\widetilde\psi\|_\infty\le\varepsilon$. It follows that
\[
 \left|\int_{G\times G}\varphi\,d\kappa\right|
 \le\varepsilon|\kappa|(G\times G).
\]

Letting $\varepsilon\downarrow0$ and using uniqueness in the Riesz theorem demonstrates that $\kappa=0$. The same approximation establishes the asserted uniform density in $C_0(G\times G)$.
\end{proof}

\begin{remark}
For a finite $M_n$-valued Radon measure $\mu$ one may define finite $M_n$-valued Radon measures $\kappa_1,\kappa_2$ on $G\times G$ by
\[
 (\kappa_1)_{ab}=\sum_{c=1}^n\mu_{ac}\widehat\otimes\mu_{cb},
 \qquad
 (\kappa_2)_{ab}=\sum_{c=1}^n\mu_{cb}\widehat\otimes\mu_{ac}.
\]

Then, for all Borel subsets $A,B\subseteq G$, we have
\[
 \kappa_1(A\times B)=\mu(A)\mu(B),\qquad
 \kappa_2(A\times B)=\mu(B)\mu(A).
\]

In particular, no identification of the Borel $\sigma$-algebra of $G\times G$ with a product $\sigma$-algebra is required when these measures are used in Lemma~\ref{lem:tensor-leg}.
\end{remark}

\textbf{Proof roadmap for the noncommutative step.}
Theorem~\ref{thm:general-matrix-kawada} reduces, after restricting to the range of the projection $\mu(G)$, to showing that all the values of a normalized positive convolution idempotent commute.   The following lemma performs that reduction in four logically distinct steps.  First, a Naimark dilation $\mu(D)=V^*P(D)V$ provides a unitary integrated operator $\mathcal U_\pi$ for each unitary representation $\pi$; the compressed coefficient $A_\pi$ is a contractive idempotent, hence an orthogonal projection.  Second, the two coefficient operators $P_1,P_2$ for representations on separate group legs share one matrix-coefficient leg.  Their joint compression $R=A_{\pi\otimes\rho}$ satisfies the {range identities} $RP_2=P_1P_2$, and $RP_1=P_2P_1$. These identities---rather than any general rule for multiplying compressions-make $P_1P_2$ an orthogonal projection and force $P_1P_2=P_2P_1$.  Third, the products of matrix coefficients are realized as finite Radon  measures $\kappa_1,\kappa_2$ on $G\times G$, so no identification of two potentially different Borel $\sigma$-algebras is required.  Finally, taking the outer tensor product of two regular  representations and applying Lemma~\ref{lem:faithful-measure-convolution} on $G\times G$ gives $\kappa_1=\kappa_2$ and thus $\mu(D)\mu(D')=\mu(D')\mu(D)$ for all Borel sets $D,D'$. The preceding  Radon-product lemma provides the measure-theoretic justification, including without second-countability; no property of the harmonic semigroup is invoked.

\begin{lemma}
\label{lem:tensor-leg}
Let $G$ be a locally compact group, let $P$ be a regular projection-valued
measure on $G$ acting on a Hilbert space $\mathcal H$, and let
$V:\C^n\longrightarrow\mathcal H$ be an isometry. Assume that
$\mu(D)=V^*P(D)V$ satisfies $\mu*\mu=\mu$. For every strongly continuous
unitary representation $\pi:G\longrightarrow\mathcal U(\mathcal H_\pi)$ the following
assertions hold.
\begin{enumerate}
\item[(a)] There is a unique unitary $\mathcal U_\pi$ on
$\mathcal H_\pi\otimes\mathcal H$ whose matrix coefficients are
\begin{equation}\label{eq:tensor-leg-U}
 \langle\mathcal U_\pi(\xi\otimes\eta),\xi'\otimes\eta'\rangle
 =\int_G\langle\pi(x)\xi,\xi'\rangle\,
            d\langle P(x)\eta,\eta'\rangle.
\end{equation}
For a second strongly continuous unitary representation $\rho$,
write $\pi\otimes\rho$ for the \emph{inner} tensor product representation
$x\mapsto\pi(x)\otimes\rho(x)$. On
$\mathcal H_\pi\otimes\mathcal H_\rho\otimes\mathcal H$,
\begin{equation}\label{eq:tensor-leg-factor}
 \mathcal U_{\pi\otimes\rho}
 =\mathcal U_\pi^{(13)}\mathcal U_\rho^{(23)}
 =\mathcal U_\rho^{(23)}\mathcal U_\pi^{(13)}.
\end{equation}
\item[(b)] With $W_\pi=I_{\mathcal H_\pi}\otimes V$, the compression
\[
 A_\pi=W_\pi^*\mathcal U_\pi W_\pi
       =\sum_{a,b=1}^n\pi(\mu_{ab})\otimes e_{ab}
\]
is an orthogonal projection on $\mathcal H_\pi\otimes\C^n$.
\item[(c)] If $\zeta\in\Ran A_\pi$, then
$\mathcal U_\pi W_\pi\zeta=W_\pi\zeta$.
\item[(d)] On $\mathcal K=\mathcal H_\pi\otimes\mathcal H_\rho\otimes\C^n$ put $P_1=A_\pi^{(13)}$, $P_2=A_\rho^{(23)}$, and $R=A_{\pi\otimes\rho}$. Then $RP_2=P_1P_2$, $RP_1=P_2P_1$, and hence $P_1P_2=P_2P_1$. 
\item[(e)] Let $\pi\boxtimes\rho$ be the outer tensor product representation of $G\times G$ on $\mathcal H_\pi\otimes\mathcal H_\rho$. Define finite matrix-valued Radon measures $\kappa_1,\kappa_2$ by
\[
 (\kappa_1)_{ad}=\sum_{b=1}^n\mu_{ab}\widehat\otimes\mu_{bd},
 \qquad
 (\kappa_2)_{ad}=\sum_{b=1}^n\mu_{bd}\widehat\otimes\mu_{ab},
\]
where $\widehat\otimes$ denotes the Radon product of scalar measures from Lemma~\ref{lem:radon-prerequisites}(b). Thus, on Borel rectangles,
\[
 \kappa_1(D\times D')=\mu(D)\mu(D'),\qquad
 \kappa_2(D\times D')=\mu(D')\mu(D),
\]
and
\[
 P_1P_2=\bigl((\pi\boxtimes\rho)\otimes\mathrm{id}\bigr)(\kappa_1),
 \qquad
 P_2P_1=\bigl((\pi\boxtimes\rho)\otimes\mathrm{id}\bigr)(\kappa_2).
\]

For a representation $\Sigma$ of $G\times G$, the notation on the right means $\sum_{a,d}\Sigma((\kappa_i)_{ad})\otimes e_{ad}$.
\end{enumerate}
\end{lemma}

\begin{proof} 
(a)  Since the projections $P(D)$ commute, a maximal orthogonal family of cyclic subspaces gives $\mathcal H=\bigoplus_\iota\mathcal H_\iota$, with each $\mathcal H_\iota$ reducing every $P(D)$. For a cyclic vector $\eta_\iota$, let $\nu_\iota(D)=\langle P(D)\eta_\iota,\eta_\iota\rangle$. The assignment $P(D)\eta_\iota\mapsto\mathbf1_D$ extends to a unitary $\mathcal H_\iota\simeq L^2(G,\nu_\iota)$ under which $P(D)$ acts as multiplication by $\mathbf1_D$. This is valid without any separability assumption. Hence
\[
 \mathcal H_\pi\otimes\mathcal H_\iota
 \simeq L^2(G,\nu_\iota;\mathcal H_\pi).
\]

On this summand define $(\mathcal U_\pi F)(x)=\pi(x)F(x)$. Lemma~\ref{lem:radon-prerequisites}(a) shows that this is a unitary, including when $\mathcal H_\pi$ is nonseparable. Take the orthogonal  sum of these unitaries over $\iota$.

For vectors in one cyclic summand, the coefficient measure of $P$ has Radon--Nikodym density $f\overline g$ with respect to $\nu_\iota$; for vectors in distinct summands the coefficient measure  vanishes. This proves \eqref{eq:tensor-leg-U}. For completeness, the total variation of each coefficient measure satisfies 
\[
 \bigl|\langle P(\cdot)\eta,\eta'\rangle\bigr|(G)
 \le\|\eta\|\,\|\eta'\|,
\]
by Cauchy--Schwarz over finite Borel partitions. Thus the right-hand side of \eqref{eq:tensor-leg-U} is bounded by $\|\xi\|\|\xi'\|\|\eta\|\|\eta'\|$. Since elementary tensors span a dense  subspace, the coefficients uniquely determine $\mathcal U_\pi$.

On $\mathcal H_\pi\otimes\mathcal H_\rho\otimes\mathcal H_\iota$, viewed as $L^2(G,\nu_\iota;\mathcal H_\pi\otimes\mathcal H_\rho)$, the operators $\mathcal U_\pi^{(13)}$ and $\mathcal  U_\rho^{(23)}$ are multiplication respectively by $\pi(x)\otimes I$ and $I\otimes\rho(x)$. These multipliers commute pointwise and have product $\pi(x)\otimes\rho(x)$, proving \eqref{eq:tensor-leg-factor}.

(b) Applying \eqref{eq:tensor-leg-U} to $Vu,Vu'$ gives, with the stated matrix-entry convention, 
\[
 A_\pi=\sum_{a,b}\pi(\mu_{ab})\otimes e_{ab}.
\]

Since $\pi(\nu*\nu')=\pi(\nu)\pi(\nu')$, we have
\[
 A_\pi^2
 =\sum_{a,b}\pi\Bigl(\sum_c\mu_{ac}*\mu_{cb}\Bigr)\otimes e_{ab}
 =\sum_{a,b}\pi((\mu*\mu)_{ab})\otimes e_{ab}=A_\pi.
\]

Moreover $\|A_\pi\|\le1$, since it is a compression of a unitary by an isometry. Every contractive Hilbert-space idempotent is an orthogonal projection; hence $A_\pi=A_\pi^*=A_\pi^2$.

(c) If $A_\pi\zeta=\zeta$ and $w=W_\pi\zeta$, then
\[
 \langle\mathcal U_\pi w,w\rangle
 =\langle A_\pi\zeta,\zeta\rangle
 =\|\zeta\|^2=\|w\|^2=\|\mathcal U_\pi w\|^2.
\]

It follows that $\|\mathcal U_\pi w-w\|^2=0$.

(d) Write $W_3=I_{\mathcal H_\pi}\otimes I_{\mathcal H_\rho}\otimes V$, $U_1=\mathcal U_\pi^{(13)}$ and $U_2=\mathcal U_\rho^{(23)}$. Then $U_1U_2=U_2U_1=\mathcal U_{\pi\otimes\rho}$ and
\[
 P_i=W_3^*U_iW_3\quad(i=1,2),\qquad
 R=W_3^*U_1U_2W_3.
\]

By (c), applied to the appropriate tensor factors, we obtain $U_1W_3P_1=W_3P_1$, and $U_2W_3P_2=W_3P_2$. Thus the two {range identities} are
\[
 RP_2=W_3^*U_1U_2W_3P_2
     =W_3^*U_1W_3P_2=P_1P_2
\]
and, using the reversed factorization,
\[
 RP_1=W_3^*U_2U_1W_3P_1=P_2P_1.
\]

All three operators $P_1,P_2,R$ are orthogonal projections by (b). Taking adjoints gives $P_1R=P_1P_2$, so
\[
 (P_1P_2)^2=(P_1R)(RP_2)=P_1R^2P_2
   =P_1RP_2=P_1(P_1P_2)=P_1P_2.
\]

Since $\|P_1P_2\|\le1$, the product is an orthogonal projection. Therefore $P_1P_2=(P_1P_2)^*=P_2P_1$. Notice that no unrestricted identity
$W_3^*U_1U_2W_3=(W_3^*U_1W_3)(W_3^*U_2W_3)$ is used.

(e) Matrix-unit multiplication gives
\[
 P_1P_2
 =\sum_{a,d}\left(\sum_b
        \pi(\mu_{ab})\otimes\rho(\mu_{bd})\right)\otimes e_{ad},
\]
and
\[
 P_2P_1
 =\sum_{a,d}\left(\sum_b
        \pi(\mu_{bd})\otimes\rho(\mu_{ab})\right)\otimes e_{ad}.
\]

For finite scalar Radon measures $\nu,\nu'$,
\[
 (\pi\boxtimes\rho)(\nu\widehat\otimes\nu')
   =\pi(\nu)\otimes\rho(\nu').
\]

Both sides are bounded operators, so it suffices to compare their coefficients at elementary tensors $\xi\otimes\xi'$ and $\eta\otimes\eta'$. These are the integral of the bounded continuous  function $(x,y)\mapsto\langle\pi(x)\xi,\eta\rangle\langle\rho(y)\xi',\eta'\rangle$ against $\nu\widehat\otimes\nu'$ and the product of the two one-variable integrals. They agree by the Fubini  theorem for Radon products of finite Radon measures \cite[\S7.4]{FollandRA}, applied to the positive and negative parts of the real and imaginary parts of $\nu$ and $\nu'$. Lemma~\ref{lem:radon-prerequisites}(b) makes every product measure on the right a finite Radon measure on $G\times G$, so the stated formulas for $\kappa_1,\kappa_2$ and their integrated  representations follow. On rectangles their matrix entries have exactly the stated ordered products. No identification of the Borel $\sigma$-algebra of $G\times G$ with the product Borel $\sigma$-algebra is required.
\end{proof}

\begin{remark}\label{rem:tensor-leg-range-identity}
In general $R\ne P_1P_2$, even though $P_1$ and $P_2$ commute. For example, let $G=C_2$ and let $\mu=m_G$ be the scalar Haar measure. Take both $\pi$ and $\rho$ to be the nontrivial  one-dimensional character. Then $A_\pi=A_\rho=0$, while $A_{\pi\otimes\rho}=1$. Hence $P_1P_2=0$ but $R=1$. The range identities in part~(d), rather than a product identity for
compressions, are exactly what the commutativity proof requires.
\end{remark}

\begin{proof}[Proof of Theorem~\ref{thm:general-matrix-kawada}]
Put $P_0=\mu(G)$. By Proposition~\ref{prop:idempotent-mass}, $P_0$ is an orthogonal projection and $\mu(D)=P_0\mu(D)P_0$ for every Borel set $D$. If $P_0=0$, positivity gives $\mu=0$, and the  asserted decomposition is empty. Otherwise, restrict the matrix coefficients to the $r$-dimensional space $H_0=P_0\C^n$; on this space the total mass is $I_{H_0}$. We prove the result for this  normalized measure and extend the resulting coefficient projections by zero on $H_0^\perp$. For notational simplicity, write again $\mu$ and $\C^r$ for the compressed measure and its coefficient space, so that $\mu(G)=I_r$.

If $G$ is compact, put $X=G$. Otherwise let $X=G_\infty=G\cup\{\infty\}$ be its one-point compactification and extend $\mu$ by zero at $\infty$, namely $\widetilde\mu(B)=\mu(B\cap G)$ for Borel  $B\subseteq X$. In either case, $\widetilde\mu$ is a finite regular positive $M_r$-valued measure on the compact Hausdorff space $X$ with $\widetilde\mu(X)=I_r$. The map
\[
 \Phi:C(X)\longrightarrow M_r,\qquad
 \Phi(f)=\int_X f(x)\,d\widetilde\mu(x)
\]
is unital and positive, hence completely positive, because $C(X)$ is commutative. Naimark's dilation theorem \cite[Theorem~4.6]{Paulsen} supplies a Hilbert space $\mathcal H$, a regular projection-valued measure $\widetilde P$ on $X$, and an isometry $V:\C^r\longrightarrow\mathcal H$ with
\[
 \widetilde\mu(B)=V^*\widetilde P(B)V
 \qquad(B\subseteq X\ \text{Borel}).
\]

If $X=G_\infty$, then
$V^*\widetilde P(\{\infty\})V=0$. Since
$\widetilde P(\{\infty\})$ is an orthogonal projection,
\[
 \|\widetilde P(\{\infty\})Vu\|^2
 =\langle V^*\widetilde P(\{\infty\})Vu,u\rangle=0
 \qquad(u\in\C^r).
\]

Thus $\widetilde P(\{\infty\})V=0$. The closed subspace $\mathcal H_0=\Ran\widetilde P(G)$ therefore contains $\Ran V$ and reduces every $\widetilde P(D)$, because the projections of a spectral  measure commute. Restricting $\widetilde P(D)$ to $\mathcal H_0$, for Borel $D\subseteq G$, produces a regular projection-valued measure $P$ on $G$ with $P(G)=I_{\mathcal H_0}$ and  $\mu(D)=V^*P(D)V$. For compact $G$ use $P=\widetilde P$ directly. The cyclic spectral construction in Lemma~\ref{lem:tensor-leg} requires no separability assumption.

We now investigate the automatic commutativity. Let $\lambda_G$ be the left regular representation of $G$ on $L^2(G)$. Apply Lemma~\ref{lem:tensor-leg}(d)--(e) with $\pi=\rho=\lambda_G$. The outer  tensor product $\lambda_G\boxtimes\lambda_G$ identifies with the left regular representation $\lambda_{G\times G}$, also without second countability. Indeed, the Radon product $m_G\widehat\otimes  m_G$ of left Haar measures \cite[\S7.4]{FollandRA} is a nonzero left-invariant Radon measure on $G\times G$, hence a left Haar measure there. The map $f\otimes g\mapsto f(x)g(y)$ is isometric on  $C_c(G)\odot C_c(G)$, by the defining property of the Radon product on $C_c(G\times G)$, and its range is dense in $L^2(G\times G)$: by the Stone--Weierstrass argument in the proof of  Lemma~\ref{lem:radon-prerequisites}(b), every $\varphi\in C_c(G\times G)$ is a uniform limit of such tensors supported in a fixed compact set. This gives a unitary $L^2(G)\otimes L^2(G)\cong  L^2(G\times G)$ intertwining $\lambda_G\boxtimes\lambda_G$ with $\lambda_{G\times G}$. With $\kappa=\kappa_1-\kappa_2$ from that lemma, 
\[
 \bigl(\lambda_{G\times G}\otimes\mathrm{id}_{M_r}\bigr)(\kappa)
 =P_1P_2-P_2P_1=0.
\]

Taking each matrix entry gives $\lambda_{G\times G}(\kappa_{ab})=0$. The $\kappa_{ab}$ are finite complex Radon measures on $G\times G$ by Lemma~\ref{lem:radon-prerequisites}(b). By the  faithfulness of left convolution (Lemma~\ref{lem:faithful-measure-convolution}, applied to $G\times G$), each $\kappa_{ab}$ vanishes. Hence $\kappa_1=\kappa_2$. Evaluating these measures on Borel rectangles gives
\[
 \mu(D)\mu(D')=\kappa_1(D\times D')
 =\kappa_2(D\times D')=\mu(D')\mu(D)
 \qquad(D,D'\subseteq G\text{ Borel}).
\]

No identification of the product Borel $\sigma$-algebra with the Borel $\sigma$-algebra of $G\times G$ is needed here. The family of matrix values is now commuting, so Theorem~\ref{thm:commuting-idempotents} gives the orthogonal Haar decomposition
\[
 \mu=\sum_{j=1}^s m_{K_j}\otimes E_j,
 \qquad \sum_{j=1}^s E_j=I_r,
\]
with pairwise distinct compact subgroups $K_j$ and mutually
orthogonal nonzero projections $E_j$. Returning to the original
coefficient space means extending each $E_j$ by zero on
$\ker P_0$, and yields $\sum_jE_j=P_0$.

For completeness, we verify the intrinsic description of the channels.
For any compact subgroup $K\le G$, define
\[
 H_K=\{u\in\C^n:\mu(D)u=m_K(D)u
                 \text{ for every Borel }D\subseteq G\}.
\]

Taking $D=G$ shows that $H_K\subseteq\Ran P_0$. For $u\in H_K$, decompose $u=\sum_jE_ju$. Since the summands are orthogonal, for each $j$ with $E_ju\ne0$ one has $m_{K_j}=m_K$ as scalar measures.  These Haar measures have supports $K_j$ and $K$, respectively, so $K_j=K$. It follows that $H_{K_j}=\Ran E_j$ and that $H_K=\{0\}$ if $K$ is not among the listed subgroups. Thus the channels and their projections are uniquely determined, up to permutation. Finally, for $0\ne u\in\Ran E_j$, the vector measure $\mu_u(D):=\mu(D)u$ is $m_{K_j}(D)u$. Thus
\[
 \delta_h*\mu_u=\mu_u
 \quad\Longleftrightarrow\quad
 \delta_h*m_{K_j}=m_{K_j}
 \quad\Longleftrightarrow\quad hK_j=K_j
 \quad\Longleftrightarrow\quad h\in K_j.
\]

This proves the stated stabilizer formula and the uniqueness.
\end{proof}

\begin{remark}\label{rem:kawada-countability}
The proof of Theorem~\ref{thm:general-matrix-kawada} uses no countability assumption on $G$. Such an assumption could enter at four points, and each is handled directly. First, convolution and the measures $\kappa_1,\kappa_2$ on $G\times G$ are defined through Radon products, so the possible difference between the Borel $\sigma$-algebra of $G\times G$ and $\mathcal B(G)\otimes\mathcal B(G)$ plays no role (Lemma~\ref{lem:radon-prerequisites}(b)). Second, the dilation space and the representation spaces may be nonseparable; the multipliers $x\mapsto\pi(x)$ are handled by Pettis'  theorem (Lemma~\ref{lem:radon-prerequisites}(a)) and the cyclic decomposition in Lemma~\ref{lem:tensor-leg}(a). Third, $L^2(G\times G)=L^2(G)\otimes L^2(G)$ holds for the Radon product of Haar  measures, as shown in the proof. Fourth, faithfulness of $\lambda_{G\times G}$ on finite Radon measures is proved by evaluating a continuous function at the identity (Lemma~\ref{lem:faithful-measure-convolution}), not by a disintegration or a countable separating family. The scalar theorem used in Theorem~\ref{thm:commuting-idempotents} holds on every locally compact group \cite{Heyer,Kelley,Pym}.
\end{remark}

Figure~\ref{fig:kawada-channels} summarizes the structure given by Theorem~\ref{thm:general-matrix-kawada}.

\begin{figure}[htbp]
\centering
\begin{tikzpicture}[>=Stealth, line cap=round, line join=round,
                    every node/.style={font=\small}]
  \draw[thick] (0,0) rectangle (3.0,5.2);
  \fill[gray!18] (0,0) rectangle (3.0,1.05);
  \draw (0,1.05) -- (3.0,1.05);
  \draw (0,2.15) -- (3.0,2.15);
  \draw (0,3.05) -- (3.0,3.05);
  \draw (0,4.15) -- (3.0,4.15);
  \node at (1.5,4.67) {$E_1\C^n$};
  \node at (1.5,3.60) {$E_2\C^n$};
  \node at (1.5,2.62) {$\vdots$};
  \node at (1.5,1.60) {$E_r\C^n$};
  \node[align=center, font=\footnotesize] at (1.5,0.52) {$\ker\mu(G)$};
  \node[anchor=south] at (1.5,5.25) {coefficient space $\C^n$};
  \node[draw, rounded corners=2pt, minimum width=3.4cm, minimum height=0.72cm]
       (h1) at (7.3,4.67) {$m_{K_1}$ on $K_1\le G$};
  \node[draw, rounded corners=2pt, minimum width=3.4cm, minimum height=0.72cm]
       (h2) at (7.3,3.60) {$m_{K_2}$ on $K_2\le G$};
  \node at (7.3,2.62) {$\vdots$};
  \node[draw, rounded corners=2pt, minimum width=3.4cm, minimum height=0.72cm]
       (hr) at (7.3,1.60) {$m_{K_r}$ on $K_r\le G$};
  \node (h0) at (7.3,0.52) {$0$};
  \draw[->, thick] (3.05,4.67) -- (h1.west);
  \draw[->, thick] (3.05,3.60) -- (h2.west);
  \draw[->, thick] (3.05,1.60) -- (hr.west);
  \draw[->, thick, gray] (3.05,0.52) -- (h0.west);
  \node[anchor=south] at (7.3,5.25) {compact subgroups of $G$};
  \node[anchor=north, font=\footnotesize] at (4.4,-0.25)
       {$\mu(D)u=\sum_{j=1}^r m_{K_j}(D)\,E_ju$ \quad for all Borel $D\subseteq G$ and $u\in\C^n$};
\end{tikzpicture}
\caption{The conclusion of Theorem~\ref{thm:general-matrix-kawada}. The
coefficient space splits orthogonally into the channels $E_j\C^n$ and
$\ker\mu(G)$. On the $j$th channel, $\mu$ acts as the normalized Haar
measure of a compact subgroup $K_j$, and it vanishes on $\ker\mu(G)$. Only
the coefficient channels are orthogonal: the subgroups $K_j$ can be taken
distinct, but they need not be disjoint or nested
(Example~\ref{ex:idempotent-two-subgroups}).}
\label{fig:kawada-channels}
\end{figure}

\begin{corollary}\label{cor:canonical-kawada}
Let $G$ be a locally compact group and let $n\ge1$.  For each finite set $\mathcal F$ of \emph{distinct} compact subgroups $K$ of $G$, let $(E_K)_{K\in\mathcal F}$ be pairwise orthogonal nonzero  orthogonal projections in $M_n$.  The assignment
\[
 (\mathcal F,(E_K)_{K\in\mathcal F})\longmapsto
 \mu=\sum_{K\in\mathcal F}m_K\otimes E_K
\]
is a bijection onto the positive convolution idempotents in $M(G,M_n)$.  Necessarily $|\mathcal F|\le n$. Moreover, for a compact subgroup $K$ define the intrinsic subspace
\[
 H_K(\mu)=\{u\in\C^n:\ \mu(D)u=m_K(D)u
                \ \text{for every Borel }D\subseteq G\}.
\]

Then the inverse correspondence is
\[
 \mathcal F_\mu=\{K\le G:K\text{ compact and }H_K(\mu)\ne\{0\}\},
 \qquad E_K=\operatorname{proj}_{H_K(\mu)}.
\]

Furthermore,
\[
 \mu(G)=\sum_{K\in\mathcal F_\mu}E_K,\qquad
 \ker\mu(G)=\bigcap_{K\in\mathcal F_\mu}\ker E_K.
\]

The empty family corresponds to $\mu=0$, the only positive idempotent with $\mu(G)=0$.  A \emph{nonzero} idempotent has a single nonzero Haar channel $m_K\otimes E$ precisely when $|\mathcal F_\mu|=1$.
\end{corollary}

\begin{proof}
For any such family, the sum is positive and idempotent by the converse of Theorem~\ref{thm:commuting-idempotents}.  Nonzero mutually orthogonal projections have orthogonal nonzero ranges, so their number is at most $n$.  Theorem~\ref{thm:general-matrix-kawada} gives existence and uniqueness of the decomposition after equal subgroups have been grouped.  Its intrinsic channel formula  gives $H_K(\mu)=\Ran E_K$ for an occurring subgroup and $H_K(\mu)=\{0\}$ otherwise, proving injectivity and the inverse formula.  The mass and kernel identities follow by summing the mutually  orthogonal projections.  Finally, $0\le\mu(D)\le\mu(G)$ for each Borel $D$, so $\mu(G)=0$ implies $\mu=0$; conversely the empty sum is zero.
\end{proof}

\begin{corollary}\label{cor:classification-trace-support}
Let $\mu\ge0$ satisfy $\mu*\mu=\mu$, and write its unique Haar-channel decomposition as in Theorem~\ref{thm:general-matrix-kawada}.  Then
\[
 \operatorname{Tr}\mu=\sum_{j=1}^r (\operatorname{rank}E_j)m_{K_j},
 \qquad
 \operatorname{supp}\mu=\bigcup_{j=1}^rK_j.
\]

In particular $\mu$ is compactly supported and $\operatorname{Tr}\mu(G)=\operatorname{rank}\mu(G)$.  Furthermore,
\[
 \mu(\{e\})=\sum_{\substack{1\le j\le r\\K_j\text{ finite}}}
       \frac{1}{|K_j|}E_j,
\]
so $\mu$ has an atom if and only if $\mu(\{e\})\ne0$.
\end{corollary}
\begin{proof}
Taking traces gives the displayed scalar positive measure. The support of a finite sum of positive measures is the union of their supports, since each $K_j$ is closed and the union is finite. Every compact subgroup contains $e$; its Haar probability measure has an atom at $e$ exactly when it is finite, in which case $m_{K_j}(\{e\})=|K_j|^{-1}$.  If $\mu$ has an atom at a point, at  least one of the scalar Haar summands has an atom, so one $K_j$ is finite and $\mu(\{e\})\ne0$.  The total-mass trace identity is immediate from the orthogonality of the projections. 
\end{proof}

\begin{corollary}\label{cor:finite-dimensional-cstar-kawada}
Let $\mathcal A$ be a finite-dimensional $C^*$-algebra, and let $\mu$ be a finite regular $\mathcal A$-valued Borel measure on a locally compact group $G$.  Suppose that $\mu(D)$ is positive in $\mathcal A$ for every Borel set $D$ and that $\mu*\mu=\mu$, with convolution using multiplication in $\mathcal A$. Then there is a unique finite family of distinct compact subgroups $K_j\le G$  and pairwise orthogonal nonzero projections $p_j\in\mathcal A$ such that
\[
   \mu=\sum_j m_{K_j}\otimes p_j,
   \qquad \sum_jp_j=\mu(G).
\]

In particular the values of $\mu$ generate a finite-dimensional commutative $C^*$-subalgebra of $\mathcal A$.
\end{corollary}
\begin{proof}
Write $\mathcal A\cong\bigoplus_{\ell=1}^sM_{n_\ell}$. The coordinate measures are positive convolution idempotents, since convolution is computed componentwise in this decomposition. Theorem~\ref{thm:general-matrix-kawada} classifies each coordinate measure.  Form the finite union of their compact subgroup families, and for each subgroup combine its coefficient projections  into the corresponding direct-sum projection $p_j\in\mathcal A$ (using zero projections on coordinates where the subgroup does not occur). The $p_j$ are pairwise orthogonal, and their sum is  $\mu(G)$. Coordinatewise uniqueness gives uniqueness of the grouped family. The concluding commutativity statement follows from the decomposition.
\end{proof}

\begin{corollary}\label{cor:general-kawada-consequences}
Let $G$ be locally compact.  Then 
\begin{enumerate}
\item[(a)] If $\mu\ge0$ and $T_\mu^2=T_\mu$ on a space on which the
convolution representation is faithful, then, with the finite Haar
family of Corollary~\ref{cor:canonical-kawada},
\[
 (T_\mu f)(x)=\sum_{K\in\mathcal F_\mu}
 E_K\int_K f(k^{-1}x)\,dm_K(k).
\]

In tensor notation this is $\sum_K T_{m_K}\otimes L_{E_K}$, where
$L_{E_K}(A)=E_KA$ on $M_n$.  In particular, the conclusion applies to
harmonic projections represented by positive convolution-idempotent
measures in Theorem~\ref{thm:harmonic-measure-representation}.

\item[(b)] The existence and commutativity conclusions in Theorems~\ref{thm:lca-matrix-idempotents}, \ref{thm:discrete-matrix-idempotents}, \ref{thm:ac-matrix-idempotents}, \ref{thm:transverse-matrix-idempotents}, \ref{thm:atomic-splitting}, \ref{thm:countably-transverse-idempotents}, and Proposition~\ref{prop:homogeneous-idempotents} are subsumed by  Theorem~\ref{thm:general-matrix-kawada} where the respective hypotheses apply.  Their more detailed conclusions--finite support, subgroup openness under absolute continuity, and transverse
or atomic splittings---remain useful and follow from the canonical form together with the additional support hypotheses.

\item[(c)] Every \emph{nonzero} positive convolution idempotent has a nonzero common eigenvector for its matrix values; every unit vector in an occurring channel $\Ran E_j$ is one.  \item[(d)]  Let $(\sigma_t)_{t>0}$ be a normalized positive matrix convolution semigroup whose operators are strongly continuous on $L^2(G,M_n)$, as required by Theorem~\ref{thm:harmonic-measure-representation}.  If $\mu$ denotes its harmonic idempotent, then
\[
 \mu=\sum_jm_{K_j}\otimes E_j.
\]

For every occurring $j$, every $t>0$, and every unit vector
$v\in\Ran E_j$,
\[
 \sigma_t(G\setminus K_j)E_j=E_j\sigma_t(G\setminus K_j)=0,
 \qquad \sigma_t(K_j)E_j=E_j=E_j\sigma_t(K_j),
\]
and
\[
 K_j=\overline{\bigcup_{t>0}
       \operatorname{supp}\langle\sigma_t(\cdot)v,v\rangle}.
\]
No commutativity of the values of $\sigma_t$ is assumed.  If the
harmonic idempotent is zero, the sum is empty and there are no nonzero
harmonic channels.
\end{enumerate}
\end{corollary}

\begin{proof}
Let $\mu\ge0$ be convolution-idempotent.  By
Theorem~\ref{thm:general-matrix-kawada}, write
$\mu=\sum_jm_{K_j}\otimes E_j$ with pairwise orthogonal $E_j$.  If
$u_j=E_ju$, then
\begin{equation}\label{eq:kawada-on-vectors}
 \mu(D)u=\sum_jm_{K_j}(D)u_j,\qquad
 \langle\mu(D)u,u\rangle=\sum_j\|u_j\|^2m_{K_j}(D).
\end{equation}

{(a)} If $T_\mu^2=T_\mu$, the composition formula gives $T_{\mu*\mu}=T_\mu$; faithfulness (Lemma~\ref{lem:faithful-measure-convolution}) implies $\mu*\mu=\mu$.  The decomposition above gives the  displayed action. The harmonic projection in the cited theorem has such a representation. The adjective ``positive'' here refers to the underlying matrix measure; left multiplication by a positive  matrix need not be a positive map on the $C^*$-algebra $M_n$.

{(b)} The universal theorem proves the general Haar-channel existence statement.  Under the extra hypotheses, familiar elementary facts about Haar measures recover the corresponding refinements.   In particular, $m_K$ has atoms exactly when the compact group $K$ is finite.  Thus the atomic part is
\[
 \mu_a=\sum_{j:K_j\text{ finite}}m_{K_j}\otimes E_j,
\]
whereas the remaining Haar channels are atomless.  When $G$ is discrete every compact $K_j$ is finite.  When a Haar channel is absolutely continuous with respect to Haar measure on a carrying closed subgroup $H$, its subgroup $K_j$ must be open in $H$. The more elaborate transverse and countable splittings require the particular support hypotheses of their original statements.

{(c)} If $\mu\ne0$, let $v$ be a unit vector in an occurring $\Ran E_j$. Equation~\eqref{eq:kawada-on-vectors} gives $\mu(D)v=m_{K_j}(D)v$ for all Borel $D$, so $v$ is a common eigenvector.

{(d)} Under the stated strong-continuity hypothesis, Theorem~\ref{thm:harmonic-measure-representation} supplies a positive convolution idempotent $\mu$ as the harmonic-measure limit.  The unrestricted classification gives its Haar-channel decomposition. The assertions on support and confinement are then exactly Theorem~\ref{thm:harmonic-channel-structure}(b).  The claim is about the channels of the harmonic idempotent; it does not assert that the individual matrix measures $\sigma_t$ have commuting values.
\end{proof}

\begin{remark}
\label{rem:general-kawada-states}
For a positive matrix measure $\mu$ set $\phi_\mu(f)=\int_Gf\,d\mu$, a completely positive map $C_0(G)\longrightarrow M_n$.  Matrix convolution is represented formally, with the usual multiplier  extensions in the noncompact case, by 
\[
 \phi_{\mu*\nu}
   =m\circ(\phi_\mu\otimes\phi_\nu)\circ\Delta,
 \qquad m(A\otimes B)=AB.
\]

For $n\ge2$, the linear multiplication map $m:M_n\otimes M_n\longrightarrow M_n$ is not positive: two positive noncommuting matrices $A,B$ can have $AB\ne(AB)^*$.  Therefore matrix convolution of  {arbitrary} positive matrix measures need not be positive, and the displayed factorization is not a general completely positive construction. {For a positive idempotent}, however,
$\phi_{\mu*\mu}=\phi_\mu$ is completely positive.  Thus one must not infer non-complete-positivity of that particular composite merely from the failure of positivity of $m$.

There is a natural scalar positive functional on $C_0(G)\otimes M_n$ defined, for positive $\mu$, by
\[
 \omega_\mu(f\otimes a)
  =\frac1n\operatorname{Tr}\bigl(\phi_\mu(f)a\bigr),
 \qquad
 \|\omega_\mu\|=\frac{\operatorname{Tr}\mu(G)}n.
\]

It is a {state} when $\mu(G)=I_n$.  For a nonzero positive idempotent with $P=\mu(G)$ of rank $d<n$, the corresponding normalized state is obtained by using $d^{-1}$ instead of $n^{-1}$; for  $\mu=0$ the functional is zero, not a state.  The matrix multiplication in the preceding formula is not a $*$-homomorphism and does not by itself define a coproduct on $M_n$.  In particular,  writing ``$\Delta\otimes\mathrm{id}$'' does not specify a comultiplication $C_0(G)\otimes M_n\longrightarrow M((C_0(G)\otimes M_n)^{\otimes2})$: one would additionally have to specify a suitable  map on the matrix factor.  Consequently the ordinary scalar convolution of states on a quantum semigroup should not be identified, without further structure, with the present matrix-valued convolution.

Kasprzak, Khosravi and So\l tan \cite{KasprzakKhosraviSoltan} study scalar idempotent states on locally compact quantum semigroups with weak cancellation.  Under their compact-support hypothesis  they obtain a Haar-state realization on a  {proto-compact quantum hypergroup} described by an operator system; they also characterize when the relevant object is a compact quantum subgroup.  Our  problem concerns instead positive $M_n$-valued measures whose convolution multiplies matrix coefficients.  The two settings are related in spirit but have different convolution structures, so the  cited state theorem does not directly yield Theorem~\ref{thm:general-matrix-kawada}.

The classification should likewise be distinguished from Greenleaf's results on scalar contractive idempotent measures \cite{Greenleaf} and their locally compact quantum-group extension by  Neufang, Salmi, Skalski and Spronk \cite{NeufangSalmiSkalskiSpronk}.  For example, even after the present classification, a scalar compression of a matrix idempotent need not itself be idempotent.   If $\mu=m_{K_1}\otimes E_1+m_{K_2}\otimes E_2$ and $\|E_1v\|^2,\|E_2v\|^2>0$, then
\[
 \langle\mu(\cdot)v,v\rangle
 =\|E_1v\|^2m_{K_1}+\|E_2v\|^2m_{K_2},
\]
which is generally not convolution-idempotent when $K_1\ne K_2$. Thus one cannot classify a general positive matrix idempotent by applying the scalar theorem separately to all quadratic-form compressions before the matrix structure has been analyzed.  Cohen's idempotent theorem for abelian groups \cite{CohenIdempotent} and Host's theorem on the Fourier--Stieltjes algebra  \cite{HostIdempotents} address different scalar idempotent problems.  Spronk \cite{SpronkCommuting} studies when {two} contractive idempotent Haar measures commute: for compact subgroups $K_1,K_2$ the pertinent condition is that $K_1K_2$ be a subgroup. In the present orthogonal Haar-channel sum, mixed products vanish  because $E_iE_j=0$ for $i\ne j$, whether or not $K_iK_j$ is a subgroup.   Thus no such compatibility is required between distinct channels.

The operator-valued positive-definite-function lemma (Lemma~\ref{lem:pd-projections-abelian}) supplies an auxiliary Fourier argument in the abelian case, but does not on its own prove the nonabelian classification.  There is also a scalar Fourier-operator approach to idempotent probability measures due to Loynes \cite{Loynes}.  The compressed operators $A_\pi$ of Lemma~\ref{lem:tensor-leg} are analogous integrated coefficients; the additional matrix problem is to show that their different tensor legs commute.  The approaches of Kelley \cite{Kelley} and Pym
\cite{Pym} likewise work with scalar measures.  Chu's monograph \cite{Chu} develops spectral, semigroup and harmonic-function theory for matrix convolution operators, providing related background. In the literature reviewed for the present paper, we have not found Theorem~\ref{thm:general-matrix-kawada} for $n\ge2$, even for finite groups.  This novelty comparison is made with the specific ordered coefficient convolution used here.
\end{remark}

\begin{theorem}\label{thm:discrete-matrix-idempotents}
Let $G$ be an arbitrary discrete group and let $\mu\in M(G,M_n)$ be positive and convolution-idempotent. Then $\mu$ has finite support, its point-mass matrices commute pairwise, and
\begin{equation}\label{eq:finite-haar-decomp}
 \mu=\sum_{j=1}^{r}m_{H_j}\otimes E_j,
 \qquad \sum_{j=1}^{r}E_j=\mu(G),
\end{equation}
where the $H_j$ are pairwise distinct finite subgroups of $G$ and the $E_j$ are mutually orthogonal nonzero projections. The empty sum represents $\mu=0$. The decomposition is unique up to  permutation. Moreover, if $P=\mu(G)\ne0$ and $A_e=\mu(\{e\})$, then $A_e$ is positive and invertible on $P\C^n$, and 
\begin{equation}\label{eq:discrete-identity-atom}
 A_e=\sum_{j=1}^{r}|H_j|^{-1}E_j.
\end{equation}
\end{theorem}
\begin{proof}
By Proposition~\ref{prop:idempotent-mass}, $P=\mu(G)$ is an orthogonal projection and all values of $\mu$ vanish on $\ker P$. If $P=0$, positivity gives $\mu=0$. Otherwise, restrict to $P\C^n\cong\C^d$ and assume $\mu(G)=I_d$. Set $A_g=\mu(\{g\})\ge0$. The positive matrices $A_g$ satisfy $\sum_{g\in G}A_g=I_d$ in norm; only countably many are nonzero. On $\ell^2(G,\C^d)$ define
\[
 (Qf)(x)=\sum_{g\in G}A_g f(g^{-1}x).
\]

For every $f\in\ell^2(G,\C^d)$, the matrix Jensen/variance identity is
\begin{equation}\label{eq:discrete-variance-equality}
 \sum_{x,g\in G}\left\|A_g^{1/2}
 \bigl(f(g^{-1}x)-Qf(x)\bigr)\right\|^2
 =\|f\|_2^2-\|Qf\|_2^2.
\end{equation}

Indeed, expand the squares, use $\sum_gA_g=I_d$ and $\sum_g A_g f(g^{-1}x)=Qf(x)$, then change variables $x\mapsto g^{-1}x$ in the first term. Tonelli's theorem justifies the sums. In particular,  $Q$ is a contraction. Since $Q^2=Q$ by convolution idempotence, $Q$ is an orthogonal projection. As the $A_g$ are self-adjoint, $Q^*=T_{\check\mu}$. Applying $Q=Q^*$ to the delta vector
$\delta_e v$ shows that
\begin{equation}\label{eq:discrete-symmetry}
 A_{g^{-1}}=A_g \qquad (g\in G).
\end{equation}

For $f=Qf$, equation~\eqref{eq:discrete-variance-equality} reduces to a sum of nonnegative terms equal to zero. Hence
\[
 A_g\bigl(f(g^{-1}x)-f(x)\bigr)=0
 \qquad(g,x\in G).
\]

Take $f=Qh$, with $h$ finitely supported. Comparing the
coefficient of the freely chosen vector $h(z^{-1}x)$ on both
sides gives
\begin{equation}\label{eq:discrete-variance-identity}
 A_gA_{g^{-1}z}=A_gA_z\qquad(g,z\in G).
\end{equation}

With $z=e$, equation~\eqref{eq:discrete-symmetry}  implies $A_g^2=A_gA_e$. Taking adjoints also gives $A_g^2=A_eA_g$, so $A_e$ commutes with every $A_g$.
The coefficient at $e$ of $\mu*\mu=\mu$ gives
\[
 A_e=\sum_gA_gA_{g^{-1}}=\sum_gA_g^2.
\]

If $A_ev=0$, then $A_gv=0$ for every $g$; but $v=\sum_gA_gv=0$. Thus $A_e$ is strictly positive and invertible on $\C^d$. Set $P_g=A_e^{-1}A_g$. The commutation and
$A_g^2=A_eA_g$ imply that $P_g$ is a self-adjoint idempotent; in particular $P_g$ is the support projection of $A_g$. Moreover,
\[
 A_g=A_eP_g,\qquad P_e=I_d,\qquad
 P_{g^{-1}}=P_g,\qquad \sum_{g\in G}P_g=A_e^{-1}.
\]

Every nonzero $P_g$ has trace at least one, whereas $\operatorname{Tr}(A_e^{-1})<\infty$. It follows that only finitely many $P_g$, and hence only finitely many $A_g$, are
nonzero. Multiplying~\eqref{eq:discrete-variance-identity} by
$A_e^{-2}$ gives
\begin{equation}\label{eq:discrete-projection-identity}
 P_gP_{g^{-1}z}=P_gP_z.
\end{equation}

For $g,h\in G$, put $P=P_g$, $S=P_h$ and $R=P_{g^{-1}h}=P_{h^{-1}g}$. Then $PS=PR$ and $SP=SR$, so taking adjoints gives $PS=RS$.
Thus
\[
 (PS)^2=(PR)(RS)=PRS=P(PS)=PS.
\]

Since $\|PS\|\le1$, it is a contractive idempotent on $\C^d$, hence an orthogonal projection. Thus $PS=SP$. All $P_g$ commute; because each commutes with $A_e$, all $A_g$
commute as well. Simultaneously diagonalize the commuting family $\{A_e\}\cup\{P_g:g\in G\}$, and let $u_1,\ldots,u_d$ be an orthonormal common eigenbasis, with 
\[
 A_eu_i=\alpha_i u_i\quad(\alpha_i>0),\qquad
 P_gu_i=p_i(g)u_i\quad(p_i(g)\in\{0,1\}).
\]

Since $P_e=I_d$, we have $p_i(e)=1$. Define $H_i=\{g\in G:p_i(g)=1\}$. The scalar version of \eqref{eq:discrete-projection-identity} says that if $g\in H_i$, then $p_i(g^{-1}z)=p_i(z)$ for every  $z$. It follows that $H_i$ is a subgroup. It is finite because $\mu$ has finite support. Also, $A_gu_i=\alpha_i\mathbf1_{H_i}(g)u_i$; normalization gives $1=\sum_g\alpha_i\mathbf1_{H_i}(g)=\alpha_i|H_i|$. Thus $\alpha_i=|H_i|^{-1}$. Grouping equal subgroups and writing $E_j$ for the projection onto the corresponding common eigenspaces gives the  asserted orthogonal Haar decomposition and~\eqref{eq:discrete-identity-atom}. Uniqueness follows from the intrinsic recovery of the distinct Haar channels in Theorem~\ref{thm:general-matrix-kawada}.
\end{proof}

\begin{remark}\label{rem:idempotent-states}
A positive normalized matrix measure $\mu(G)=I_n$ determines a
completely positive map
$\phi_\mu:C_0(G)\longrightarrow M_n$, $\phi_\mu(f)=\int f\,d\mu$.
For noncompact $G$, its integration extension to $C_b(G)$ is
unital, as is its extension to the one-point compactification.
For $f\in C_0(G)$, matrix convolution is encoded by
\[
 \phi_{\mu*\nu}(f)
   =m\bigl((\phi_\mu\otimes\phi_\nu)(\Delta f)\bigr),
 \qquad (\Delta f)(x,y)=f(xy),
\]
where the right-hand side is interpreted via the integration
extensions when $\Delta f\notin C_0(G\times G)$ and
$m:M_n\otimes M_n\longrightarrow M_n$ denotes matrix multiplication.
For $n=1$, normalized positive idempotents are scalar idempotent
probability measures, which also define idempotent states in the
classical group setting. For $n\ge2$, $m$ is not a positive map:
products of positive matrices need not be self-adjoint. Thus matrix
convolution is not in general the convolution of states associated
with a coproduct on $C_0(G)\otimes M_n$.
Notice, however, that for a positive \emph{idempotent} $\mu$ the
particular composite above equals $\phi_\mu$ and is therefore
completely positive. The issue is the general multiplication map,
not positivity of this particular composite. See also
Remark~\ref{rem:general-kawada-states}.
\end{remark}

\begin{corollary}\label{cor:multisubgroup-zeno}
For a positive idempotent $\mu=\sum_jm_{K_j}\otimes E_j$, written as in
Theorem~\ref{thm:general-matrix-kawada}, its convolution operator on
$L^2(G,S_2^n)$ is the orthogonal projection
\[
 (Qf)(x)=\sum_{j=1}^{r}E_j
 \int_{K_j}f(k^{-1}x)\,dm_{K_j}(k).
\]
It is therefore a contractive Hilbert-space Zeno projection.
It need not be contractive on $L^p(G,S_2^n)$ for $p\ne2$;
indeed, the example below fails contractivity for $p=8$ and $p=\infty$.
\end{corollary}
\begin{proof}
By Lemma~\ref{lem:contraction}, $Q$ is an $L^2(S_2)$ contraction,
and $Q^2=Q$. Hence it is an orthogonal projection.
For a direct norm-expanding test, take $G=C_2=\{e,s\}$ with
\emph{normalized} Haar measure and
$\mu=\operatorname{diag}(\delta_e,m_G)$, as in
Example~\ref{ex:idempotent-two-subgroups}. Let
\[
 f(e)=\begin{pmatrix}1&0\\0&0\end{pmatrix},\qquad
 f(s)=\begin{pmatrix}0&0\\0&1\end{pmatrix}.
\]
Then
$Qf(e)=\operatorname{diag}(1,1/2)$ and
$Qf(s)=\operatorname{diag}(0,1/2)$. In the
$L^8(G,S_2^2)$ norm,
\[
 \|f\|_8^8=\tfrac12(1+1)=1,\qquad
 \|Qf\|_8^8
 =\tfrac12\left((\sqrt5/2)^8+(1/2)^8\right)
 =\frac{313}{256}>1.
\]
Thus $\|Q\|_{L^8(S_2)\longrightarrow L^8(S_2)}>1$.
The $L^\infty(S_2)$ norm also increases because
$\|Qf(e)\|_{S_2}=\sqrt5/2>\|f\|_\infty=1$.
The word \emph{measurement} here means an orthogonal projection
on the indicated Hilbert space, not necessarily a positive or
completely positive map on the algebra of matrix-valued functions.
\end{proof}

\begin{remark}[Non-projection mass and translation]
\label{rem:nonprojection-idempotent-mass}
Convolution idempotence implies $\mu(G)^2=\mu(G)$ even without
positivity, but the idempotent need not be self-adjoint. For example,
\[
 J=\begin{pmatrix}1&1\\0&0\end{pmatrix},\qquad
 J^2=J\ne J^*,
\]
and $\delta_e\otimes J$ is a convolution idempotent which is not
positive. Therefore $\mu(G)^2\ne\mu(G)$ is impossible for any
convolution idempotent. Translation need not preserve idempotence:
$\delta_g*\delta_g=\delta_{g^2}\ne\delta_g$ when $g\ne e$.
The orthogonal Haar decomposition is an intrinsic equality of
matrix measures. A single fixed unitary change of basis may
simultaneously diagonalize its channel projections, but neither
unitary equivalence nor translation is required in the classification.
\end{remark}

% ------------------------------------------------------------------

\begin{lemma}[Uniform $L^p$ boundedness]\label{lem:uniform-p-bound}
If $\sigma\ge0$ and $\sigma(G)=I_n$, then for all $1\le p\le\infty$
\[
\|T_\sigma\|_{L^p(G,S_2^n)\to L^p(G,S_2^n)}\le n.
\]
Thus a normalized positive matrix convolution semigroup is uniformly
bounded on every spatial $L^p$ where it is strongly continuous.
\end{lemma}
\begin{proof}
Since $\sigma$ is positive, its total variation is dominated by
$\operatorname{Tr}\sigma$; hence $|\sigma|(G)\le n$.
Apply the convolution variation inequality.
\end{proof}
\begin{proposition}\label{prop:ergodic}
Let $1<p<\infty$ and suppose $(T_t)_{t\ge0}$ is a strongly continuous semigroup on $L^p(G,M_n)$ with $\sup_{t\ge0}\|T_t\|_p<\infty$. Then the Ces\`aro means
\[
 C_r=\frac1r\int_0^rT_t\,dt
\]
converge strongly as $r\longrightarrow\infty$ to a bounded projection $P_p$
onto $\hsp=\ker A$, with
$$ \ker P_p=\overline{\Ran A},\qquad P_pT_t=T_tP_p=P_p.$$

Here the closure is taken in $L^p(G,M_n)$, and
$$ \Lp=\hsp\oplus\overline{\Ran A}$$
is a topological direct sum. In general, the closure of $\Ran A$ cannot be omitted. The projection $P_p$ commutes with right translations and right multiplication by constant matrices; it is contractive whenever each $T_t$ is contractive. If $p,q\in(1,\infty)$ and the same convolution semigroup induces compatible strongly continuous semigroups on both $L^p(G,M_n)$ and $L^q(G,M_n)$, then
$$ P_pf=P_qf\quad\text{almost everywhere} \qquad(f\in L^p(G,M_n)\cap L^q(G,M_n)).$$
\end{proposition}

\begin{proof}
Since $L^p(G,M_n)$ is reflexive for $1<p<\infty$, the mean ergodic theorem for uniformly bounded $C_0$-semigroups applies \cite[Ch.~V.4]{EngelNagel}, \cite[\S1.2]{Krengel}. The uniform bound is automatic for the normalized positive matrix convolution semigroups considered here, by the elementary total-variation bound in Lemma~\ref{lem:uniform-p-bound}. Thus
$$ L^p(G,M_n)=\ker A\oplus\overline{\Ran A},$$
and $C_r$ converges strongly to the projection onto $\ker A$ along $\overline{\Ran A}$. Since $\ker A=\hsp$, the assertions about the range, kernel, and invariance follow. Commutation with right  translations and constant right matrix multiplications passes from $T_t$ to $C_r$ and then to the strong limit.

For consistency across exponents, fix $f\in L^p(G,M_n)\cap L^q(G,M_n)$. The Ces\`aro means $C_rf$ are the same measurable functions in the two spaces. Along any sequence $r_k\longrightarrow\infty$, these  functions converge in $L^p$ to $P_pf$ and in $L^q$ to $P_qf$. A subsequence converges almost everywhere to $P_pf$, and a further subsequence converges almost everywhere to $P_qf$. Uniqueness of
almost-everywhere limits gives $P_pf=P_qf$ almost everywhere.

 Take $n=1$, $G=\mathbb R$, and the Gaussian heat semigroup, whose generator is $A=d^2/dx^2$ on $D(A)=W^{2,p}(\mathbb R)$. An $L^p$ function with distributional second derivative zero is affine,  hence zero. Therefore $\ker A=\{0\}$ and $P_p=0$, so $\overline{\Ran A}=L^p(\mathbb R)$. Nevertheless $\Ran A\ne L^p(\mathbb R)$. If $A$ were surjective, then, since it is closed and injective,  the bounded inverse theorem applied to $A:D(A)\longrightarrow L^p(\mathbb R)$ (with $D(A)$ equipped with its graph norm) would imply
$$ \|f\|_p\le c\|f''\|_p\qquad(f\in W^{2,p}(\mathbb R)).$$

For nonzero $f\in C_c^\infty(\mathbb R)$ and $f_\lambda(x)=f(\lambda x)$,
$$ \|f_\lambda\|_p=\lambda^{-1/p}\|f\|_p, \qquad \|f_\lambda''\|_p=\lambda^{2-1/p}\|f''\|_p.$$

The proposed inequality fails as $\lambda\downarrow0$. Hence $\Ran A$ is a dense proper, and therefore nonclosed, subspace of $\ker P_p=L^p(\mathbb R)$.
\end{proof}

\section{Harmonic projections}\label{sec:harmonic}
% ------------------------------------------------------------------

\subsection{Harmonic projections represented by matrix measures}
\label{subsec:harmonic-measures}

Although a bounded operator on $L^p(G)$ commuting with right translations
need not be convolution by a finite measure when $1<p<\infty$, the
normalized positive measures defining our semigroup have uniformly
bounded variation.  This additional information forces the ergodic
projection itself to be a finite-measure convolution, even if mass escapes
to infinity.  The following result strengthens Proposition~\ref{prop:ergodic}.

The faithfulness of left convolution on matrix measures, used
repeatedly below, is Lemma~\ref{lem:faithful-measure-convolution} in
Section~\ref{subsec:general-kawada}.

\begin{theorem}\label{thm:harmonic-measure-representation}
Let $\{\sigma_t\}_{t>0}$ be a normalized positive matrix convolution
semigroup, and assume that $T_t=T_{\sigma_t}$ is strongly continuous on
$L^2(G,M_n)$. For $R>0$ form the Ces\`aro-averaged matrix measures
\begin{equation}\label{eq:ergodic-measure-average}
 \mu_R=\frac1R\int_0^R\sigma_t\,dt,
\end{equation}
where the integral is interpreted weak-$*$ in
$M(G,M_n)=C_0(G,M_n)^*$ using the trace pairing.
Then there is a unique $\mu\in M(G,M_n)$ satisfying
\begin{equation}\label{eq:ergodic-measure-weakstar}
 \mu_R\xrightarrow[R\longrightarrow\infty]{\ \sigma(M(G,M_n),C_0(G,M_n))\ }\mu.
\end{equation}

This measure is positive, $0\le\mu(G)\le I_n$, and
\begin{equation}\label{eq:ergodic-idempotent-absorption}
 \mu*\mu=\mu,\qquad
 \sigma_t*\mu=\mu*\sigma_t=\mu\quad(t>0),\qquad
 \check\mu=\mu.
\end{equation}

In particular $E:=\mu(G)$ is an orthogonal projection, and
$\mu(D)=E\mu(D)E$ for every Borel set $D$.
For every $1<p<\infty$ the semigroup is strongly continuous on
$L^p(G,M_n)$, and the mean-ergodic projection of
Proposition~\ref{prop:ergodic} is the \emph{same} convolution operator:
\begin{equation}\label{eq:harmonic-Pp-is-measure}
 P_p=T_\mu,\qquad
 H_S^p(G,M_n)=\Ran T_\mu,\qquad
 \ker P_p=\overline{\Ran A_p}.
\end{equation}

Here $A_p$ is the generator on $L^p(G,M_n)$.
\end{theorem}

\begin{proof}
We first consider the existence of the measure averages. For each $t>0$, positivity and $\sigma_t(G)=I_n$ imply
\[
 \|\sigma_t\|_{\mathrm{var}}
 \le\operatorname{Tr}\sigma_t(G)=n.
\]

Indeed, for any finite Borel partition $(D_j)$, $\sum_j\|\sigma_t(D_j)\|\le\sum_j\operatorname{Tr}\sigma_t(D_j)=n$. The same bound holds for the averaged measures once defined. To justify the  weak-$*$ integral, note that $t\mapsto\sigma_t$ is weak-$*$ continuous. In fact, if $t_i\longrightarrow t$, any weak-$*$ cluster point $\nu$ of the uniformly bounded measures $\sigma_{t_i}$  satisfies, for $f,h\in C_c(G,M_n)$,
\[
 \langle T_\nu f,h\rangle_2
 =\lim_i\langle T_{t_i}f,h\rangle_2
 =\langle T_tf,h\rangle_2,
\]
where the first equality holds along the convergent subnet. The lemma implies $\nu=\sigma_t$.  Hence all cluster points agree, which proves weak-$*$ continuity and makes \eqref{eq:ergodic-measure-average} well-defined by the Riesz representation theorem. Each $\mu_R$ is positive and $\mu_R(G)=I_n$; in particular $\|\mu_R\|_{\mathrm{var}}\le n$. By Banach--Alaoglu the family $(\mu_R)$ has weak-$*$ cluster points along subnets with $R\longrightarrow\infty$.  Fix such a point $\mu$. Positivity passes to the vague limit, so $\mu$ is positive; moreover, $0\le\mu(G)\le I_n$, as is seen by testing against positive compactly supported scalar functions bounded by one. For $f,h\in C_c(G,M_n)$ define 
\[
 K_{f,h}(y)=\int_G f(y^{-1}x)h(x)^*\,dx\in C_c(G,M_n).
\]

The continuity and compact support follow from the compact supports of $f,h$. With the Hilbert--Schmidt pairing,
\begin{equation}\label{eq:matrix-correlation-pairing}
 \langle T_{\mu_R}f,h\rangle_2
 =\int_G\operatorname{Tr}\bigl(d\mu_R(y)K_{f,h}(y)\bigr).
\end{equation}

Vague convergence therefore gives weak convergence of the associated convolution operators when tested on $C_c$ functions. But $T_{\mu_R}=R^{-1}\int_0^R T_t\,dt$; by Proposition~\ref{prop:ergodic}  on $L^2$ these operators converge strongly to $P_2$. Thus, $\langle T_\mu f,h\rangle_2=\langle P_2f,h\rangle_2$, for all $f,h\in C_c(G,M_n)$, hence $T_\mu=P_2$. Any other weak-$*$ cluster point  has the same convolution operator, and is equal to $\mu$ by the lemma. A bounded family in a weak-$*$ compact set with a unique cluster point converges to it, proving \eqref{eq:ergodic-measure-weakstar} for the full family. 

The mean-ergodic projection satisfies $P_2^2=P_2$ and $T_tP_2=P_2T_t=P_2$. Because convolution representations are faithful, we obtain
\[
 \mu*\mu=\mu,\qquad\sigma_t*\mu=\mu*\sigma_t=\mu.
\]

As every $T_t$ is a contraction on $L^2$, $P_2$ is a contractive projection on a Hilbert space and therefore orthogonal. By Lemma~\ref{lem:adjoint}, $T_{\check\mu}=T_\mu^*=T_\mu$, so  $\check\mu=\mu$. Proposition~\ref{prop:idempotent-mass} proves that $E=\mu(G)$ is an orthogonal projection and $\mu(D)=E\mu(D)E$. Fix $1<p<\infty$, $p\ne2$. Let $q=1$ if
$p<2$ and $q=\infty$ if $p>2$, and define $\theta\in(0,1)$ by $1/p=\theta/2+(1-\theta)/q$. For $g\in L^2\cap L^q(G,S_2^n)$, H\"older's inequality applied to the scalar function $\|g(\cdot)\|_{S_2}$ gives
\[
 \|g\|_p\le\|g\|_2^{\theta}\,\|g\|_q^{1-\theta}.
\]

Let $f\in C_c(G,M_n)$, a dense subspace of $L^p(G,M_n)$ contained in $L^2\cap L^q$, and put $g_t=T_tf-f$. By Lemma~\ref{lem:uniform-p-bound}, $\|T_t\|_{L^q(S_2)\longrightarrow L^q(S_2)}\le n$, so  $\|g_t\|_q\le(n+1)\|f\|_q$, for all $t$. By strong continuity on $L^2$, $\|g_t\|_2\longrightarrow0$ as $t\downarrow0$. Hence $\|T_tf-f\|_p\longrightarrow0$. Since $\sup_t\|T_t\|_{p\longrightarrow  p}\le n$ (Lemma~\ref{lem:uniform-p-bound}), an $\varepsilon/3$-argument extends the convergence to all $f\in L^p(G,M_n)$. So $(T_t)$ is a bounded $C_0$-semigroup on $L^p$, and  Proposition~\ref{prop:ergodic} applies for each $p$.

For $f,h\in C_c(G,M_n)$, the same identity \eqref{eq:matrix-correlation-pairing}, together with \eqref{eq:ergodic-measure-weakstar}, shows that the $L^p$-strong limit of $T_{\mu_R}f$ has the same  pairing with $h$ as $T_\mu f$. Hence $P_pf=T_\mu f$ on $C_c$, and boundedness implies that the identity on all of $L^p$. The remaining assertions follow from Proposition~\ref{prop:ergodic}.
\end{proof}

\begin{remark}\label{rem:convolution}
For $p=1$, bounded operators commuting with right translations are convolutions by finite measures (Wendel's theorem \cite{Wendel}). For $1<p<\infty$ such a representation is false for general  bounded translation-commuting operators: $L^p$-multipliers need not be finite measures. For abelian $G$ the bounded operators on $L^2(G,M_n)$ that commute with translations and right matrix  multiplication are exactly the matrix Fourier multipliers with symbols in $L^\infty(\widehat G,M_n)$ \cite[Theorem~3.1]{EbadianJabbariOaM}; see also \cite{EbadianJabbariMIA} for  translation-invariant operators on $L^p(G,M_n)$. Theorem~\ref{thm:harmonic-measure-representation} uses the extra positivity and normalization of the {underlying semigroup measures}, not  translation commutation alone. The measure $\mu$ is independent of $p$ even though its convolution operator acts on
different $L^p$ spaces.
\end{remark}

\begin{corollary}\label{cor:harmonic-tightness}
Under Theorem~\ref{thm:harmonic-measure-representation}, the following are equivalent
\begin{enumerate}
\item[(i)] $\mu(G)=I_n$;
\item[(ii)] the scalar positive measures $\{\operatorname{Tr}\mu_R:R>0\}$ are \emph{asymptotically tight}: for every $\varepsilon>0$ there exist a compact set $K\subseteq G$ and $R_0>0$ such that
\[
 \operatorname{Tr}\mu_R(G\setminus K)<\varepsilon
 \quad(R\ge R_0).
\]
\end{enumerate}

In particular, (i) holds for compact $G$. Without tightness the projection $\mu(G)$ can be proper or zero.
\end{corollary}
\begin{proof}
Put $\lambda_R=\operatorname{Tr}\mu_R$ and $\lambda=\operatorname{Tr}\mu$, so $\lambda_R(G)=n$ and $\lambda_R\longrightarrow\lambda$ vaguely. Suppose that $\lambda(G)=n$. Choose $\phi\in C_c(G)$ with  $0\le\phi\le1$ and $\int\phi\,d\lambda>n-\varepsilon$. Vague convergence implies $\int\phi\,d\lambda_R>n-\varepsilon$ for all large $R$, so $\lambda_R(G\setminus\operatorname{supp}\phi) <\varepsilon$. Conversely, if the family is asymptotically tight, choose a compact $K$ with $\lambda_R(G\setminus K)<\varepsilon$ for all large $R$ and a function $\phi\in C_c(G)$ equal to one on  $K$, with $0\le\phi\le1$. Then $\int\phi\,d\lambda=\lim_R\int\phi\,d\lambda_R \ge n-\varepsilon$, so $\lambda(G)=n$. Since $0\le I_n-\mu(G)$ and $\operatorname{Tr}(I_n-\mu(G))=n-\lambda(G)$,
this is equivalent to $\mu(G)=I_n$.
\end{proof}

\subsection{Noncompact harmonic spaces and compact matrix confinement}
\label{subsec:noncompact-matrix-harmonic}

\begin{example}\label{ex:compact}
Let $G$ be compact with normalized Haar measure $m_G$, and $\sigma_t=\sigma_t^{0}I_n$ with $\{\sigma_t^0\}$ an adapted (support generates a dense subgroup) scalar convolution semigroup of  probability measures. Then $\hsp$ consists of the constant functions $G\longrightarrow M_n$, and $P_pf=(m_GI_n)\star f$. For genuinely matrix-valued semigroups the harmonic space can be larger: if  $n=2$ and $\sigma_t=\operatorname{diag}(\sigma^0_t,\delta_e)$, then the second row of every $f$ is harmonic. This does not contradict the Liouville theorem of Chu on compact groups  \cite[Proposition~21]{ChuMatrixHarmonic}, which requires the support of every pure-state compression $\langle\sigma(\cdot)\xi,\xi\rangle$ to generate a dense subgroup; here the compression along $e_2$ is $\delta_e$.
\end{example}

\begin{remark}\label{rem:noncompact}
Let $n=1$ and $G$ non-compact, and let $\{\sigma_t\}$ be adapted. If $f\in L^2(G)$ is harmonic, then $\|\sigma_t\star f\|_2=\|f\|_2$ and strict convexity of $L^2$ force $f$ to be invariant under  left translation by every element of the closed subgroup generated by $\bigcup_t\operatorname{supp}\sigma_t$, i.e.\ by $G$; so $f$ is constant, hence $0$ because $G$ has infinite Haar measure.  Thus $P_2=0$. For abelian $G$ and $1\le p<\infty$ the same follows from the Choquet--Deny theorem \cite{ChoquetDeny} applied to $\varphi*f$, $\varphi\in C_c(G)$. Every non-unimodular group is  non-compact, so this applies there. The matrix-valued counterpart, with its exact compact-confinement criterion, is given in Theorem~\ref{thm:matrix-harmonic-compact-confinement}.
\end{remark}

For scalar adapted convolution semigroups on noncompact groups, Remark~\ref{rem:noncompact} gives $H_S^2=\{0\}$.  In the matrix-valued case it is not enough that the union of the supports of $\sigma_t$ generate $G$: one matrix direction may remain inside a compact region even while other directions explore the whole group. The following characterization applies on arbitrary locally  compact groups, without symmetry or unimodularity assumptions.  It concerns finite-$p$ harmonic functions rather than the bounded Choquet--Deny problem; compare the matrix-valued abelian theorem  of \cite{ChuHilberdinkHowroyd} and the discussion in \cite{Chu}.

For $C\Subset G$ (a compact subset), put 
\begin{equation}\label{eq:trapped-matrix-subspace}
 V_C=\bigcap_{t>0}\ker\bigl(\sigma_t(G\setminus C)\bigr),
 \qquad V_{\mathrm{trap}}=\bigcup_{C\Subset G}V_C
 \subseteq\mathbb C^n.
\end{equation}

Here $V_{\mathrm{trap}}$ is a vector subspace: if $v\in V_C$ and
$w\in V_D$, positivity and monotonicity imply
$v,w\in V_{C\cup D}$ and hence $v+w\in V_{C\cup D}$.

\begin{theorem}%[Compact-confinement characterization of matrix harmonic spaces]
\label{thm:matrix-harmonic-compact-confinement}
Assume that the hypotheses of Theorem~\ref{thm:harmonic-measure-representation}, and write $E=\mu(G)$ for the total-mass projection of its representing positive convolution idempotent. Then
\begin{equation}\label{eq:harmonic-trapped-equality}
 {\ V_{\mathrm{trap}}=\Ran E.\ }
\end{equation}

The following are equivalent 
\begin{enumerate}
\item[(i)] $H_S^2(G,M_n)\ne\{0\}$;
\item[(ii)] $H_S^p(G,M_n)\ne\{0\}$ for some (equivalently every)
$1<p<\infty$;
\item[(iii)] $E\ne0$;
\item[(iv)] there are $v\in\mathbb C^n\setminus\{0\}$ and a compact
$C\subseteq G$ such that
\begin{equation}\label{eq:trapped-matrix-condition}
 \sigma_t(G\setminus C)v=0\qquad(t>0).
\end{equation}
\end{enumerate}

Equivalently, for the scalar probability measures $\nu_{t,v}(D)=\langle\sigma_t(D)v,v\rangle/\|v\|^2$, condition \eqref{eq:trapped-matrix-condition} is $\nu_{t,v}(C)=1$ for every $t>0$. 
\end{theorem}

\begin{proof}
Let $0\ne f\in H_S^2(G,M_n)$.  Set
\[
 M_f=\int_Gf(x)f(x)^*\,dx\ge0,
 \qquad F_f=\text{the orthogonal projection onto }\Ran M_f.
\]

For $y\in G$ define the positive matrix
\[
 D_f(y)=\int_G\bigl(f(y^{-1}x)-f(x)\bigr)
             \bigl(f(y^{-1}x)-f(x)\bigr)^*\,dx.
\]

Expanding the square, using left invariance of Haar measure and $\sigma_t(G)=I_n$, implies
\begin{equation}\label{eq:harmonic-matrix-energy}
\begin{split}
 \int_G\operatorname{Tr}\bigl(d\sigma_t(y)D_f(y)\bigr)
 &=2\|f\|_2^2-2\operatorname{Re}\langle T_tf,f\rangle_2\\
 &=0\qquad(t>0),
\end{split}
\end{equation}
since $T_tf=f$. The matrix autocorrelation
\[
 R_f(y)=\int_G f(y^{-1}x)f(x)^*\,dx
\]
belongs to $C_0(G,M_n)$: its entries are coefficients of the left-regular representation, which vanish at infinity by approximation with compactly supported $L^2$ functions.  Moreover,
\[
 D_f(y)=2M_f-R_f(y)-R_f(y)^*\longrightarrow2M_f
 \quad(y\longrightarrow\infty).
\]

The matrices $D_f(y)$ vanish on $\ker M_f$.  Since $M_f$ is strictly positive on $\Ran F_f$, there are a compact $C_f\subseteq G$ and $c_f>0$ for which
\begin{equation}\label{eq:harmonic-energy-coercivity}
 D_f(y)\ge c_f F_f\qquad(y\notin C_f).
\end{equation}

The integrands in \eqref{eq:harmonic-matrix-energy} are nonnegative. Equations \eqref{eq:harmonic-matrix-energy}-\eqref{eq:harmonic-energy-coercivity} therefore imply
\[
 0\ge c_f\operatorname{Tr}\bigl(F_f\sigma_t(G\setminus C_f)\bigr).
\]

Positivity of $\sigma_t(G\setminus C_f)$ implies that $\sigma_t(G\setminus C_f)F_f=0$ for every $t>0$, and hence
\begin{equation}\label{eq:harmonic-row-trap}
 \Ran F_f\subseteq V_{\mathrm{trap}}.
\end{equation}

Take $v\in V_C$.  Choose $\varphi\in C_c(G)$ with $0\le\varphi\le1$ and $\varphi=1$ on $C$. Positivity and normalization show that
\[
 \left\langle\left(\int_G\varphi\,d\sigma_t\right)v,v\right\rangle
 =\|v\|^2\qquad(t>0).
\]

The same holds for the Ces\`aro measures $\mu_R$; their weak-* convergence to $\mu$ implies 
\[
 \left\langle\left(\int_G\varphi\,d\mu\right)v,v\right\rangle
 =\|v\|^2.
\]

As $0\le\int\varphi\,d\mu\le E\le I_n$ and $E$ is an orthogonal projection, $Ev=v$. Thus $V_{\mathrm{trap}}\subseteq\Ran E$. Every $f\in H_S^2$ satisfies $f=T_\mu f$.
Since $\mu(D)=E\mu(D)E$, the range of $f(x)$ is contained in $\Ran E$ for almost all $x$.  Let $W$ be the span of the column spaces $\Ran M_f$ over all $f\in H_S^2$. Then $W\subseteq\Ran E$, and  \eqref{eq:harmonic-row-trap} shows $W\subseteq V_{\mathrm{trap}}$. Suppose that some nonzero $v\in\Ran E$ is orthogonal to $W$. For every scalar $\psi\in C_c(G)$ the harmonic function $T_\mu(\psi I_n)$ consequently satisfies
\[
 v^*T_\mu(\psi I_n)(x)=0\quad\text{for almost all }x.
\]

The left side is continuous, so it vanishes everywhere.  An approximate identity $\psi$ then shows that the row-valued measure $v^*\mu$ is zero.  In particular $v^*E=0$, contradicting $Ev=v$. Therefore $W=\Ran E$, and $\Ran E\subseteq V_{\mathrm{trap}}$. This proves \eqref{eq:harmonic-trapped-equality}.

Finally, Theorem~\ref{thm:harmonic-measure-representation} implies that $H_S^p=\Ran T_\mu$ for every $1<p<\infty$; faithfulness of left convolution (Lemma~\ref{lem:faithful-measure-convolution}) implies that $T_\mu\ne0$ if and only if $\mu\ne0$ if and only if $E\ne0$. The equivalences follow.
\end{proof}

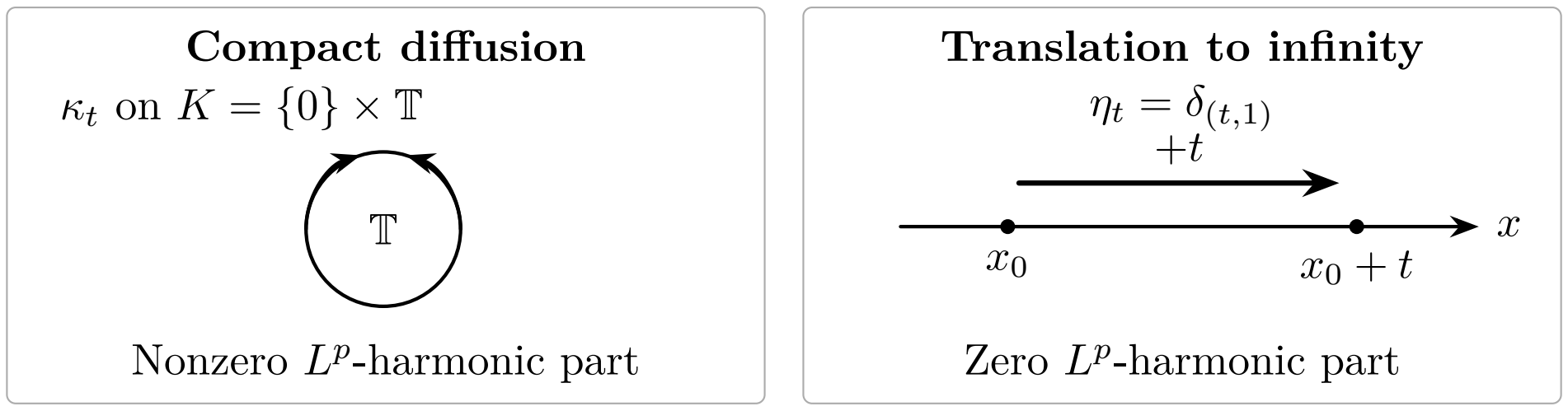
\begin{figure}[htbp]
\centering
\begin{tikzpicture}[
    >=Stealth,
    line cap=round,
    line join=round,
    every node/.style={font=\small}
]
  % Compact heat channel: motion is periodic in the angular coordinate.
  \draw[gray!65, thin, rounded corners=2pt] (0,0) rectangle (6.08,3.18);
  \node[font=\small\bfseries] at (3.04,2.84)
       {Compact diffusion};
  \node[anchor=west] at (0.29,2.37)
       {$\kappa_t$ on $K=\{0\}\times\mathbb T$};
  \draw[thick] (3.02,1.40) circle[radius=0.62];
  \draw[->, thick] (3.02+0.62,1.40)
       arc[start angle=0,end angle=73,radius=0.62];
  \draw[->, thick] (3.02-0.62,1.40)
       arc[start angle=180,end angle=107,radius=0.62];
  \node at (3.02,1.40) {$\mathbb T$};
  \node[anchor=center] at (3.04,0.32)
       {Nonzero $L^p$-harmonic part};
  % Noncompact translation channel: the orbit escapes in R.
  \draw[gray!65, thin, rounded corners=2pt] (6.39,0) rectangle (12.47,3.18);
  \node[font=\small\bfseries] at (9.43,2.84)
       {Translation to infinity};
  \node at (9.43,2.37) {$\eta_t=\delta_{(t,1)}$};
  \draw[->, thick] (7.17,1.42) -- (11.81,1.42)
       node[right] {$x$};
  \fill (8.03,1.42) circle[radius=1.7pt];
  \fill (10.83,1.42) circle[radius=1.7pt];
  \draw[->, very thick] (8.12,1.77) -- (10.69,1.77)
       node[midway,above] {$+t$};
  \node[below] at (8.03,1.36) {$x_0$};
  \node[below] at (10.83,1.36) {$x_0+t$};
  \node at (9.43,0.32)
       {Zero $L^p$-harmonic part};
\end{tikzpicture}
\caption{The two channels in Example~\ref{ex:noncompact-evolving-harmonic}.
The first diffuses in the compact angular coordinate (for each fixed
$x\in\mathbb R$), so an $L^p$-harmonic component survives. The second
translates along the noncompact $\mathbb R$-coordinate; its finite-$p$
harmonic component vanishes.}
\label{fig:two-matrix-channels}
\end{figure}

\medskip\noindent\emph{Why no tightness is needed.} The proof never uses
tightness of $\{\sigma_t\}$ or of the Ces\`aro means $\{\mu_R\}$. By
Corollary~\ref{cor:harmonic-tightness}, tightness is equivalent to
$\mu(G)=I_n$, so it holds exactly when every direction is confined, and
it fails in Example~\ref{ex:noncompact-evolving-harmonic}
(Figure~\ref{fig:two-matrix-channels}). Two ingredients replace it. The
first is the exact energy identity \eqref{eq:harmonic-matrix-energy},
valid for each fixed $t$. The second is the coercivity
\eqref{eq:harmonic-energy-coercivity} at infinity, which comes from the
fact that the autocorrelation $R_f$ of a square-integrable function lies
in $C_0(G,M_n)$. This is where finite-$p$ harmonicity enters. Together
they confine the directions $\Ran F_f$ that carry harmonic functions,
while mass in all other directions may escape to infinity. Tightness is
then a \emph{conclusion} in the confined directions: for $v\in V_C$ the
scalar laws $\langle\sigma_t(\cdot)v,v\rangle$ are carried by the fixed compact set
$C$.

\begin{corollary}%[Compact support of the harmonic idempotent]
\label{cor:harmonic-idempotent-compact-support}
Under Theorem~\ref{thm:matrix-harmonic-compact-confinement}, the
representing measure $\mu$ has compact support. More precisely, there
is a compact $C_*\subseteq G$ such that
\begin{equation}\label{eq:entire-harmonic-channel-trap}
 \sigma_t(G\setminus C_*)E=0\quad(t>0),
 \qquad\operatorname{supp}\mu\subseteq C_*.
\end{equation}
\end{corollary}
\begin{proof}
By \eqref{eq:harmonic-trapped-equality}, each vector in the
finite-dimensional space $\Ran E$ belongs to some $V_C$.
Choose a basis $v_1,\ldots,v_r$ of $\Ran E$ and corresponding
compact sets $C_1,\ldots,C_r$. The finite union
$C_*=\bigcup_{j=1}^rC_j$ satisfies
$\sigma_t(G\setminus C_*)E=0$ for every $t>0$.
(The case $E=0$ is immediate.)
This property passes to $\mu_R$ and, after testing against functions
supported in $G\setminus C_*$, to $\mu$ by weak-$*$ convergence.
Since $\mu(D)=E\mu(D)E$, it follows that $\mu$ vanishes on
$G\setminus C_*$. Hence $\operatorname{supp}\mu\subseteq C_*$.
\end{proof}

\begin{corollary}%[An exact vanishing criterion and an escape test]
\label{cor:noncompact-vanishing-test}
Under the preceding hypotheses, for $1<p<\infty$, $H_S^p(G,M_n)=\{0\}$  if and only if $\bigcap_{t>0}\ker\sigma_t(G\setminus C)=\{0\}$, for every compact $C\subseteq G$. In particular, it suffices that for every compact $C\subseteq G$ there exists $t_C>0$ for which the positive matrix $\sigma_{t_C}(G\setminus C)$ is positive definite.
\end{corollary}
\begin{proof}
By Theorem~\ref{thm:matrix-harmonic-compact-confinement}, for every $1<p<\infty$, $H_S^p(G,M_n)\ne\{0\}$ if and only if $E\ne0$. By \eqref{eq:harmonic-trapped-equality}, this holds if and only if $V_{\mathrm{trap}}=\Ran E\ne\{0\}$. Since $V_{\mathrm{trap}}=\bigcup_{C}V_C$ and $V_C=\bigcap_{t>0}\ker\sigma_t(G\setminus C)$ by \eqref{eq:trapped-matrix-subspace}, we get $H_S^p(G,M_n)=\{0\}$ if  and only if $V_C=\{0\}$ for every compact $C$, which is the displayed condition. For the sufficient condition: if $\sigma_{t_C}(G\setminus C)$ is positive definite, its kernel is $\{0\}$, and
$V_C\subseteq\ker\sigma_{t_C}(G\setminus C)=\{0\}$.
\end{proof}

\begin{example}%[Compact diffusion versus noncompact translation]
\label{ex:noncompact-evolving-harmonic}
Let $G=\mathbb R\times\mathbb T$, with componentwise multiplication,
and let $K=\{0\}\times\mathbb T$. Let $(\kappa_t)_{t\ge0}$ be
the heat probability convolution semigroup on $K$, extended to $G$,
and put $\eta_t=\delta_{(t,1)}$. Define
\[
 \sigma_t=\begin{pmatrix}\kappa_t&0\\0&\eta_t\end{pmatrix}.
\]

This is a normalized positive matrix convolution semigroup on the noncompact group $G$. Its two channels behave differently: the first undergoes genuine diffusion confined to $K$, whereas the  second translates along the noncompact factor $\mathbb R$. For $1<p<\infty$, the only $L^p$ functions fixed by the first channel are constant in the $\mathbb T$-variable; the second channel has no  nonzero $L^p$ fixed functions. Thus, writing the coordinates as $(x,\theta)$,
\[
 H_S^p(G,M_2)=
 \left\{\begin{pmatrix}a(x)&b(x)\\0&0\end{pmatrix}:
 a,b\in L^p(\mathbb R)\right\}.
\]

The Ces\`aro measures converge weak-$*$ to
\[
 \mu=\begin{pmatrix}m_K&0\\0&0\end{pmatrix},
 \qquad E=\mu(G)=\begin{pmatrix}1&0\\0&0\end{pmatrix},
 \qquad V_{\mathrm{trap}}=\mathbb Ce_1.
\]

The union of the supports of $\sigma_t$ nevertheless generates $G$. Thus adaptation formulated only in terms of the union of all matrix supports does not imply the vanishing of $H_S^p$. The  confined matrix direction need not be stationary; the compact heat channel is dynamically nontrivial. Figure~\ref{fig:two-matrix-channels} illustrates the distinction between the two channels.
\end{example}

\begin{remark}\label{rem:noncompact-harmonic-scope}
Theorem~\ref{thm:matrix-harmonic-compact-confinement} applies to non-unimodular groups as well, because left translation acts unitarily on $L^2$ for left Haar measure. It is a classification of \emph{finite-$p$} fixed functions, not of the bounded harmonic space: constants and more general bounded harmonic functions may survive without a compactly confined direction. The decomposition of  $\mu$ into subgroup Haar measures holds on every locally compact group by Theorem~\ref{thm:general-matrix-kawada}.
\end{remark}

\subsection{Structure and confinement of harmonic channels}\label{subsec:channel-structure}
\begin{theorem}\label{thm:harmonic-channel-structure}
\begin{enumerate}
\item[(a)] Let $v\in\Ran E$ be a unit common eigenvector of the values of
$\mu$, and put $G_v=\{h\in G:\delta_h*(\mu v)=\mu v\}$. Then $v$ lies in
exactly one channel $\Ran E_j$ of the decomposition in (b), $G_v=K_j$, and
$\mu(\cdot)v=m_{G_v}(\cdot)v$. For all $t>0$,
\[
 \sigma_t(G\setminus G_v)v=0,\qquad \sigma_t(G_v)v=v,
\]
the probability measures $\rho_t=\langle\sigma_t(\cdot)v,v\rangle$ live
on $G_v$, their Ces\`aro means converge to $m_{G_v}$, and
\[
 G_v=\overline{\textstyle\bigcup_{t>0}\operatorname{supp}\rho_t} .
\]

Under TV mixing, $\|\rho_t-m_{G_v}\|\le a(t)$.
\item[(b)] By Theorem~\ref{thm:general-matrix-kawada},
$\mu=\sum_jm_{K_j}\otimes E_j$ with pairwise distinct compact subgroups
$K_j$ and mutually orthogonal nonzero projections $E_j$ summing to $E$;
no further hypothesis is needed. Every channel $\Ran E_j$ is confined to
its subgroup:
$\sigma_t(G\setminus K_j)E_j=E_j\sigma_t(G\setminus K_j)=0$, $\sigma_t(K_j)$
commutes with $E_j$, and
$K_j=\overline{\bigcup_t\operatorname{supp}\langle\sigma_t(\cdot)v,v\rangle}$
for every unit $v\in\Ran E_j$.
\item[(c)] If $\mu_R\longrightarrow\mu$ in total variation and, for a.e.\ $t$,
$\sigma_t\ll m_H$ for a closed subgroup $H$ carrying all $\sigma_t$, then
$\mu\ll m_H$, and the subgroups $K_j$ of (b) are open in $H$.
\end{enumerate}
\end{theorem}

\begin{proof}
(a) Write $v=\sum_jE_jv$. If $\mu(D)v=c(D)v$ for all Borel $D$, then
(b) gives $m_{K_j}(D)=c(D)$ whenever $E_jv\ne0$. Distinct subgroups have
distinct Haar measures, since $\operatorname{supp}m_{K_j}=K_j$, so exactly
one $E_jv$ is nonzero and $v\in\Ran E_j$. Hence $\mu(\cdot)v=m_{K_j}(\cdot)v$,
and Theorem~\ref{thm:general-matrix-kawada} gives $G_v=K_j$.
Absorption $\sigma_t*\mu=\mu$ applied to $v$ gives
$\int d\sigma_t(y)\,m_{G_v}(y^{-1}A)\,v=m_{G_v}(A)v$. For $A=G_v$ we have
$m_{G_v}(y^{-1}G_v)=\mathbf 1_{G_v}(y)$, so $\sigma_t(G_v)v=v$, and
positivity gives $\sigma_t(G\setminus G_v)v=0$. The Ces\`aro means of
$\rho_t$ are $\langle\mu_R(\cdot)v,v\rangle\longrightarrow\langle\mu(\cdot)v,v\rangle
=m_{G_v}$. If $S$ is the closed union of the supports, then $S\subseteq G_v$,
and the limit $m_{G_v}$ of measures carried by $S$ is carried by $S$, so
$G_v\subseteq S$. The TV bound follows from
$|\langle\nu(A)v,v\rangle|\le\|\nu\|$.
(b) Apply (a) to every unit $v\in\Ran E_j$, for which $G_v=K_j$. Then
$\sigma_t(K_j)E_j=E_j$, and taking adjoints gives
$E_j\sigma_t(K_j)=E_j$.
(c) A total-variation limit of measures absolutely continuous with respect
to $m_H$ is absolutely continuous, so $\mu\ll m_H$. The positive measure
$\tr\mu=\sum_j(\operatorname{rank}E_j)\,m_{K_j}$ is carried by the closed
set $H$ and is absolutely continuous with respect to $m_H$. Hence
$K_j\subseteq H$ and $m_H(K_j)>0$ for every $j$. A measurable subgroup of
positive Haar measure is open, because $K_jK_j^{-1}=K_j$ contains a
neighbourhood of the identity by the Steinhaus--Weil theorem
\cite{HewittRoss}; thus each $K_j$ is open in $H$. For second-countable
$G$ this is also Theorem~\ref{thm:ac-matrix-idempotents}.
\end{proof}

\begin{corollary}\label{cor:connected-ac}
Let $G$ be connected, and suppose $\sigma_t\ll m_G$ for $t>0$ and $\mu_R\longrightarrow\mu$ in total variation.% (for instance, under TV mixing). 
Then $G$ is compact and $\mu=m_G\otimes I_n$, so the harmonic projection is averaging over $G$ and there are no proper subgroup channels.
\end{corollary}

\begin{proof}
By Theorem~\ref{thm:harmonic-channel-structure}(c) every $K_j$ is a compact open subgroup of the connected group $G$, hence $K_j=G$. Full mass
gives $\sum_jE_j=I_n$.
\end{proof}

\begin{remark}
\label{rem:channel-structure-open}
Parts (a) and (b) need no hypothesis beyond those of
Theorem~\ref{thm:harmonic-measure-representation}, on any locally compact
group. By Theorem~\ref{thm:general-matrix-kawada} the harmonic idempotent
always has an orthogonal Haar-channel decomposition, the unit common
eigenvectors of its values in $\Ran E$ are exactly the unit vectors of the
channels $\Ran E_j$, and each channel is confined to its subgroup (Corollary~\ref{cor:general-kawada-consequences}(d)).
No commutation, transversality, density or total-variation hypothesis is
needed. The special-case arguments of
Appendix~\ref{sec:special-idempotents} give additional information when
their hypotheses hold: part (c) and Theorem~\ref{thm:ac-matrix-idempotents}
show that the channel subgroups are open in a carrying subgroup when the
dynamics has densities;
Theorems~\ref{thm:transverse-matrix-idempotents}
and~\ref{thm:countably-transverse-idempotents} place each channel in a
unique carrying subgroup when the carriers are transverse (for instance
finitely many maximal tori of $SU(2)$,
Corollary~\ref{cor:su2-tori-idempotents}); and
Theorem~\ref{thm:atomic-splitting} splits off the atomic channels, which
are Haar measures of finite subgroups and have finite support.
\end{remark}

\begin{theorem}\label{thm:geometric-harmonic-space}
Assume the hypotheses of Theorem~\ref{thm:harmonic-measure-representation}, and write $\mu=\sum_{j=1}^r m_{K_j}\otimes E_j$ for its harmonic idempotent. For $1<p<\infty$ define
\[
 \mathcal H_j^p=
 \bigl\{f\in L^p(G,M_n):\ E_jf=f,\quad
       f(k^{-1}\cdot)=f\text{ in }L^p\text{ for every }k\in K_j\bigr\}.
\]

Then the common harmonic space has the canonical topological direct sum decomposition
\begin{equation}\label{eq:geometric-harmonic-decomp}
 H_S^p(G,M_n)=\bigoplus_{j=1}^r\mathcal H_j^p,
 \qquad
 P_pf=\sum_{j=1}^r E_j\int_{K_j}f(k^{-1}\cdot)\,dm_{K_j}(k).
\end{equation}

For $p=2$ the sum is orthogonal.  If $r=0$, the harmonic spaces are zero for every $1<p<\infty$.
\end{theorem}
\begin{proof}
Put $Q_jf=E_j\int_{K_j}f(k^{-1}\cdot)\,dm_{K_j}(k)$. Left Haar translation is isometric on $L^p(G,M_n)$ and $K_j$ is compact, so the Bochner integral defines a bounded operator. Its range lies in  $E_jL^p$ and is invariant under all left translations from $K_j$; hence $\operatorname{Ran}Q_j\subseteq \mathcal H_j^p$.

  Conversely, a function in $\mathcal H_j^p$ is fixed by the averaging and by $E_j$, giving $\operatorname{Ran}Q_j=\mathcal H_j^p$ and $Q_j^2=Q_j$. For $i\ne j$, the scalar averaging commutes with constant left matrix multiplication, and $E_iE_j=0$; thus $Q_iQ_j=0$. The classification and Theorem~\ref{thm:harmonic-measure-representation} imply $P_p=\sum_jQ_j$.  It follows that its range is the algebraic sum of the $\mathcal H_j^p$, direct because left multiplication by $E_i$ selects its $i$th component.  The sum is topological, since it is finite and each $Q_j$ is bounded.  For  $p=2$, the averaging over a compact subgroup is self-adjoint, as is left multiplication by an orthogonal projection, and the two operations commute.  Thus each $Q_j$ is orthogonal and their ranges are mutually orthogonal.
\end{proof}

\section{A direct bounded Zeno application}\label{sec:bounded-zeno-application}
\begin{corollary}\label{cor:classified-bounded-zeno}
Let $\mu\in M(G,M_n)$ be a positive convolution idempotent and $Q=T_\mu$ on $L^2(G,M_n)$.  Then $Q$ is an orthogonal projection of the form $Q=\sum_j(T_{m_{K_j}}\otimes E_j)$.  If $B$ is a bounded  operator generating a contraction semigroup, then, for every $t\ge0$,
\[
  \lim_{N\to\infty}(Qe^{tB/N})^NQ
   =e^{t(QB)|_{\operatorname{Ran}Q}}Q
\]
in operator norm, locally uniformly in $t$.  Moreover, 
\begin{equation}\label{eq:classified-zeno-quantitative}
\bigl\|(Qe^{tB/N})^NQ-e^{t(QB)|_{\operatorname{Ran}Q}}Q\bigr\|
\le \frac{t^2\|B\|^2}{N}e^{t\|B\|}.
\end{equation}
for every $N\ge1$ and $t\ge0$. 
\end{corollary}
\begin{proof}
The projection statement follows from Theorem~\ref{thm:general-matrix-kawada} and the symmetry of Haar measure on compact subgroups.  On $H_Q=\operatorname{Ran}Q$, let $C=QB|_{H_Q}$ and $F(h)=Qe^{hB}|_{H_Q}$.  Then $\|F(h)\|\le1$, since $e^{hB}$ is contractive.  Moreover, $C$ is dissipative because $\operatorname{Re}\langle Cu,u\rangle =\operatorname{Re}\langle Bu,u\rangle\le0$ for $u\in H_Q$; hence $\|e^{hC}\|\le1$. The exponential Taylor remainder and $\|C\|\le\|B\|$ imply
\[
\|F(h)-e^{hC}\|
 \le\|Q(e^{hB}-I-hB)|_{H_Q}\|
       +\|e^{hC}-I-hC\|
 \le h^2\|B\|^2 e^{h\|B\|}.
\]

The telescoping identity for the $N$th powers of the two contractions therefore yields, with $h=t/N$,
\[
 \|F(t/N)^N-e^{tC}\|
 \le N\|F(t/N)-e^{(t/N)C}\|
 \le \frac{t^2\|B\|^2}{N}e^{t\|B\|/N}.
\]

Since $(Qe^{hB})^NQ=F(h)^NQ$, this implies \eqref{eq:classified-zeno-quantitative} and the convergence.
\end{proof}

\appendix
% ------------------------------------------------------------------
\section{Independent proofs and finer structure in special cases}
\label{sec:special-idempotents}
% ------------------------------------------------------------------

Theorem~\ref{thm:general-matrix-kawada} classifies all positive
convolution idempotents. This appendix presents a complementary approach
through stabilizers, common eigenvectors and densities. Without the
tensor-leg lemma, it proves the classification under absolute
continuity, for transverse families of subgroups and for the atomic part
of every idempotent. Together with Theorem~\ref{thm:general-matrix-kawada}
it also also gives refinements of the general theorem. The density and support-splitting arguments are independent
of the tensor-leg lemma; when a final channel classification is invoked
below, the proof indicates whether it uses the separate commuting-value
theorem or the unrestricted matrix classification. In particular, the
channel subgroups are open in a carrying subgroup under absolute
continuity, with bounded densities (Theorem~\ref{thm:ac-matrix-idempotents});
each channel lies in a unique carrying subgroup when the carriers are
transverse (Theorems~\ref{thm:transverse-matrix-idempotents}
and~\ref{thm:countably-transverse-idempotents}); and the atomic part
splits off orthogonally and has finite support
(Theorem~\ref{thm:atomic-splitting}). The main text does not depend on
this appendix.

\subsection{A complementary approach: stabilizers, eigenvectors and densities}
\label{subsec:kawada-general}

Throughout this subsection $G$ is second countable and
$\mu\in M(G,M_n)$ is a nonzero positive convolution idempotent.
By Proposition~\ref{prop:idempotent-mass}, we compress to
$\Ran\mu(G)$ and relabel its dimension $n$, so $\mu(G)=I_n$. Put
$\lambda=\tr\mu$, $W=d\mu/d\lambda$ (so $W\ge0$ and $\tr W=1$
$\lambda$-a.e.), let $P_y$ be the support projection of $W(y)$,
$T=T_\mu$ on $L^2(G,\C^n)$, and $\mathcal K=\Ran T$.
The columnwise version of the contraction proof in
Corollary~\ref{cor:matrix-idempotent-l2} shows that $T$ is an
orthogonal projection. Hence $T=T^*=T_{\check\mu}$ and
$\check\mu=\mu$ by
Lemma~\ref{lem:faithful-measure-convolution}. For $u\in\C^n$ write
$\mu u$ for the vector measure $E\mapsto\mu(E)u$ and
$\nu_u=\langle\mu(\cdot)u,u\rangle$.

\begin{lemma}\label{lem:variance-identity}
In the preceding normalized setting, for every
$f\in L^2(G,\C^n)$ one has
\begin{equation}\label{eq:matrix-jensen-variance}
 \|f\|_2^2-\|Tf\|_2^2
 =\int_G\int_G
 \bigl\|W(y)^{1/2}\bigl(f(y^{-1}x)-Tf(x)\bigr)\bigr\|^2
 \,dx\,d\lambda(y).
\end{equation}

There exists a single Borel set $Y\subseteq G$ with $\lambda(G\setminus Y)=0$ such that, for every $y\in Y$ and {every} $f\in\mathcal K=\Ran T$,
\begin{equation}\label{eq:matrix-variance-fixed-space}
 W(y)\bigl(f(y^{-1}\cdot)-f\bigr)=0
 \quad\text{in }L^2(G,\C^n).
\end{equation}

In particular, for every $y\in Y$,
\begin{equation}\label{eq:translation-absorption}
 P_y(\delta_y*\mu)=P_y\mu
 \quad\text{as }M_n\text{-valued Radon measures}.
\end{equation}

Equivalently, for every Borel $D\subseteq G$,
\begin{equation}\label{eq:translation-absorption-both-sides}
 P_y\mu(y^{-1}D)=P_y\mu(D),\qquad
 \mu(y^{-1}D)P_y=\mu(D)P_y.
\end{equation}
\end{lemma}

\begin{proof}
Choose a measurable representative of $f$ and write $L_yf(x)=f(y^{-1}x)$. Since $\lambda(G)=n$, $0\le W(y)\le I_n$ for $\lambda$-almost every $y$, and
\[
 \int_G\int_G\langle W(y)L_yf(x),L_yf(x)\rangle
 \,dx\,d\lambda(y)=\|f\|_2^2,
\]

Tonelli's theorem and left invariance of Haar measure justify all subsequent integrations. In particular, for almost every $x$ the vector integral $w(x)=\int_GW(y)L_yf(x)\,d\lambda(y)$ exists and  equals $Tf(x)$. For such $x$, expand the weighted square  
\begin{align*}
 \int_G\|W(y)^{1/2}(L_yf(x)-w(x))\|^2\,d\lambda(y) =&\int_G\langle W(y)L_yf(x),L_yf(x)\rangle\,d\lambda(y)-2\operatorname{Re}\langle w(x),w(x)\rangle\\
       &  +\langle I_nw(x),w(x)\rangle\\
 =&\int_G\langle W(y)L_yf(x),L_yf(x)\rangle\,d\lambda(y)
 -\|Tf(x)\|^2.
\end{align*}

Integrating over $x$ proves \eqref{eq:matrix-jensen-variance}. Since $T$ is an orthogonal projection, $Tf=f$ for $f\in\mathcal K$; thus the nonnegative right side of \eqref{eq:matrix-jensen-variance} is zero. Thus, for {each} $f\in\mathcal K$, we have
\[
 W(y)^{1/2}(L_yf-f)=0\quad\text{in }L^2(G,\C^n)
\]
for $\lambda$-almost every $y$. The second-countability assumption implies that $L^2(G,\C^n)$, and hence its closed subspace $\mathcal K$, is separable. Choose a countable dense family $(f_k)$ in  $\mathcal K$ and intersect the associated conull sets. We may also arrange that $W(y)\ge0$ and $\operatorname{Tr}W(y)=1$ throughout $Y$. For each $y\in Y$, the bounded operator $W(y)^{1/2}(L_y-I)$  vanishes on every $f_k$, hence on all of $\mathcal K$. Here $L_y$ is unitary, so the operator norm is at most $2\|W(y)^{1/2}\|$. Multiplying by $W(y)^{1/2}$ implies  \eqref{eq:matrix-variance-fixed-space}. The operator is bounded because $L_y$ is unitary and $\|W(y)\|\le1$ on $Y$.

Fix $y\in Y$ and let $h\in L^2(G,\C^n)$ be arbitrary. As $Th\in\mathcal K$, \eqref{eq:matrix-variance-fixed-space} implies
\begin{align*}
 0
 &=W(y)(L_y-I)Th\\
 &=W(y)\bigl((\delta_y*\mu-\mu)*h\bigr)
 =\bigl(W(y)(\delta_y*\mu-\mu)\bigr)*h.
\end{align*}

Here $L_y(\mu*h)=(\delta_y*\mu)*h$, with the order of the convolution factors as written, even for nonabelian $G$. Thus convolution by the finite matrix measure $\eta_y=W(y)(\delta_y*\mu-\mu)$  vanishes on $L^2(G,\C^n)$. It consequently vanishes on $L^2(G,M_n)$, column by column. Lemma~\ref{lem:faithful-measure-convolution} implies $\eta_y=0$.

Let $W(y)^\dagger$ denote the Moore--Penrose inverse of the positive matrix $W(y)$. Since $W(y)^\dagger W(y)=P_y$, left multiplication of the equality $\eta_y=0$ by $W(y)^\dagger$ gives  \eqref{eq:translation-absorption}. Finally $\mu(y^{-1}D)$ and $\mu(D)$ are self-adjoint, so taking adjoints gives the second equality in \eqref{eq:translation-absorption-both-sides}.
\end{proof}

\begin{remark}
\label{rem:variance-absorption-scope}
The exceptional set in Lemma~\ref{lem:variance-identity} is
independent of $f\in\mathcal K$ and of the Borel set $D$.
The absorption identity is not a claim that
$\delta_y*\mu=\mu$ for $\lambda$-almost every $y$;
it only asserts equality after left multiplication by the
support projection $P_y$. If $P_y=I_n$, it does imply
$\delta_y*\mu=\mu$ for that $y$.
\end{remark}

\begin{proposition}\label{prop:common-eigenvector}
Assume the normalized setting of this subsection: $G$ is second countable, $\mu\ge0$, $\mu*\mu=\mu$, $\mu(G)=I_n$, and $\mathcal K=\Ran T_\mu\subseteq L^2(G,\C^n)$. For $0\ne u\in\C^n$ put $\eta_u=\mu u$, $\nu_u=\langle\mu(\cdot)u,u\rangle$, and $G_u=\{h\in G:\delta_h*\eta_u=\eta_u\}$. 
\begin{enumerate}
\item[(a)] $G_u$ is a compact subgroup, and $\eta_u$ is both left and right $G_u$-invariant. Moreover,
\[
G_u=\Bigl\{h\in G:
 \langle f(h^{-1}\cdot),u\rangle=\langle f(\cdot),u\rangle
 \text{ in }L^2(G)\text{ for every }f\in\mathcal K\Bigr\}.
\]
\item[(b)] If $\|u\|=1$, then, in the order of self-adjoint operators on $L^2(G)$, $0\le T_{\nu_u}\le T_{m_{G_u}}$. 
\item[(c)] For $\|u\|=1$, the following are equivalent:
\begin{enumerate}
\item[(i)] $u$ is an eigenvector of every $\mu(D)$, $D\subseteq G$ Borel;
\item[(ii)] $\nu_u*\nu_u=\nu_u$;
\item[(iii)] $\nu_u=m_{G_u}$;
\item[(iv)] $\mu(G_u)u=u$.
\end{enumerate}
\item[(d)] Fix $n\ge1$. The following two assertions are equivalent:
\begin{enumerate}
\item[(i)] Every positive convolution idempotent in $M(G,M_r)$,
for each $1\le r\le n$, has an orthogonal Haar decomposition;
\item[(ii)] Every positive convolution idempotent in $M(G,M_r)$
with identity total mass, for each $1\le r\le n$, has a common
eigenvector.
\end{enumerate}
Both assertions hold by Theorem~\ref{thm:general-matrix-kawada}.
\end{enumerate}

In the canonical decomposition
$\mu=\sum_jm_{K_j}\otimes E_j$ with $\sum_jE_j=I_n$, one has
\begin{equation}\label{eq:common-eigenvector-stabilizer}
 G_u=\bigcap_{j:E_ju\ne0}K_j,
 \qquad
 \nu_u=\sum_j\|E_ju\|^2m_{K_j}.
\end{equation}
\end{proposition}

\begin{proof}
(a) The vector measure $\eta_u=\mu u$ is finite and nonzero, since $\eta_u(G)=u$. The stabilizer $G_u$ is a subgroup; it is closed because $h\mapsto\delta_h*\eta_u$ is continuous in the weak-*  topology of finite vector measures. To see compactness directly, observe that the total variation $|\eta_u|$ is a nonzero finite positive measure invariant under left translation by $G_u$. Choose  a compact set $C\subseteq G$ with $|\eta_u|(C)>0$. If $G_u$ were noncompact, an induction choosing $h_k\in G_u$ outside the finite unions of sets $h_iCC^{-1}$ would produce pairwise disjoint translates $h_kC$. Each would have the same positive $|\eta_u|$-mass, contradicting finiteness. Since $\check\mu=\mu$, inversion leaves $\eta_u$ invariant $\check\eta_u=\eta_u$. For $h\in G_u$, $\eta_u=\delta_h*\eta_u$ implies 
\[ \eta_u=\check\eta_u*\delta_{h^{-1}}
 =\eta_u*\delta_{h^{-1}},
\]
so right invariance follows as well. For the fixed-space characterization, put $L_hf(x)=f(h^{-1}x)$. For any $g\in L^2(G,\C^n)$ and $f=T_\mu g$, we have
\[
 \langle L_hf-f,u\rangle
 =\bigl\langle ((\delta_h*\mu-\mu)*g)(\cdot),u\bigr\rangle
 \quad\text{in }L^2(G).
\]

Thus the scalar function on the left vanishes for all $f\in\mathcal K$ if and only if every row measure $u^*(\delta_h*\mu-\mu)$ has zero left convolution representation. By scalar faithfulness  (Lemma~\ref{lem:faithful-measure-convolution}), that is equivalent to $u^*(\delta_h*\mu-\mu)=0$. Each value of $\delta_h*\mu-\mu$ is self-adjoint, so this is equivalent to $(\delta_h*\mu-\mu)u=0$, or $h\in G_u$.

(b)  Let $J_u:L^2(G)\to L^2(G,\C^n)$ be $J_u\varphi=\varphi u$ and let $A=T_{m_{G_u}}$, the orthogonal projection onto the left $G_u$-invariant subspace of $L^2(G)$. If $\|u\|=1$, then $T_{\nu_u}=J_u^*T_\mu J_u\ge0$. For $\varphi\in L^2(G)$ put $f=T_\mu J_u\varphi$. By (a), $J_u^*f=\langle f,u\rangle$ is left $G_u$-invariant, so $J_u^*f=A(J_u^*f)$. Since $T_\mu$ is an orthogonal  projection, 
\begin{align*}
 \|f\|_2^2
 &=\langle T_\mu J_u\varphi,J_u\varphi\rangle
 =\langle J_u^*f,\varphi\rangle\\
 &=\langle J_u^*f,A\varphi\rangle
 \le\|J_u^*f\|_2\,\|A\varphi\|_2
 \le\|f\|_2\,\|A\varphi\|_2.
\end{align*}

Therefore $\|f\|_2^2\le\|A\varphi\|_2^2$, including when $f=0$. As $\langle T_{\nu_u}\varphi,\varphi\rangle=\|f\|_2^2$ and $\langle A\varphi,\varphi\rangle=\|A\varphi\|_2^2$, the asserted operator inequality follows.

(c) we now characterize the common eigenvectors. Let $J=J_u$ and $B=J^*T_\mu J=T_{\nu_u}$. Since $\|u\|=1$, $J$ is an isometry, and since $T_\mu=T_\mu^2=T_\mu^*$, \begin{equation}\label{eq:common-eigenvector-defect}
 B-B^2
 =J^*T_\mu(I-JJ^*)T_\mu J
 =\bigl[(I-JJ^*)T_\mu J\bigr]^*
   \bigl[(I-JJ^*)T_\mu J\bigr].
\end{equation}

Thus,  $B^2=B$ if and only if $T_\mu$ leaves $L^2(G)\otimes\C u$ invariant. For each $v\perp u$, the latter condition means that the scalar convolution operator associated with $D\mapsto\langle\mu(D)u,v\rangle$ is zero. Faithfulness of scalar convolution then gives $\mu(D)u\in\C u$ for all Borel $D$; the converse is immediate. Also, scalar faithfulness implies that $B^2=B$ if and only if $\nu_u*\nu_u=\nu_u$. This proves (i)$\Leftrightarrow$(ii).

Under (i), normalization implies $\mu(D)u=\nu_u(D)u$ and $\nu_u(G)=1$. Thus $\nu_u$ is an idempotent probability measure, so the scalar Kawada--It\^o--Kelley theorem gives $\nu_u=m_K$ for a compact subgroup $K$. Its left stabilizer is $K$; since $\mu u=\nu_uu$, this is precisely $G_u$. Hence (i)$\Rightarrow$(iii). If (iii) holds, then $\langle\mu(G_u)u,u\rangle=\nu_u(G_u)=1$. Since $0\le\mu(G_u)\le I_n$, this forces $\mu(G_u)u=u$, proving (iv).

Now suppose that (iv) holds. Positivity gives $\mu(G\setminus G_u)u=0$. Thus the finite vector measure $\eta_u=\mu u$ is supported on the compact group $G_u$ and is left $G_u$-invariant by  definition. Uniqueness of Haar measure, applied to each coordinate of this vector measure, implies
\[
 \mu(D)u=m_{G_u}(D)\,\eta_u(G_u)
 =m_{G_u}(D)u
\]
for every Borel $D$. This proves (i).

(d) The implication (i)$\Rightarrow$(ii) is immediate: a normalized idempotent has at least one nonzero channel, and any vector in its range is a common eigenvector. Conversely, assume (ii) in  every dimension at most $n$. We prove (i) by induction on the coefficient dimension $r$. The zero idempotent gives the empty decomposition. For a nonzero positive idempotent,  Proposition~\ref{prop:idempotent-mass} reduces to the range of its total-mass projection, of dimension at most $r$, where it has identity total mass. Pick a unit common eigenvector $u$. Each  $\mu(D)$ is self-adjoint, so $\C u$ and $u^\perp$ reduce all $\mu(D)$; both diagonal corners are positive convolution idempotents. By (c) the first is $m_{G_u}\otimes uu^*$, while the second,  after reduction, has identity total mass in dimension one lower. The induction hypothesis provides its orthogonal Haar decomposition, which together with the
first corner proves (i).

By Theorem~\ref{thm:general-matrix-kawada}, write $\mu=\sum_jm_{K_j}\otimes E_j$ with distinct compact $K_j$. Then $\mu u=\sum_jm_{K_j}(E_ju)$ and the nonzero vectors $E_ju$ are mutually  orthogonal. Translation invariance of $\mu u$ is therefore equivalent to invariance of each active Haar measure, which occurs precisely for $h\in\bigcap_{j:E_ju\ne0}K_j$. Taking scalar matrix  coefficients gives the second formula in \eqref{eq:common-eigenvector-stabilizer}. Since the $K_j$ are distinct, a unit vector $u$ satisfies the equivalent conditions in (c) if and only if $u\in\Ran E_j$ for exactly one $j$.
\end{proof}

\begin{theorem}\label{thm:ac-matrix-idempotents}
Let $G$ be a second-countable locally compact group, and let $\mu\in M(G,M_n)$ be positive and convolution-idempotent. Suppose that $H\le G$ is closed, $\mu(G\setminus H)=0$, and $\mu\ll m_H$,  where $m_H$ is left Haar measure on $H$. Then the values of $\mu$ commute, and
\begin{equation}\label{eq:ac-orthogonal-haar}
 \mu=\sum_{j=1}^{r}m_{K_j}\otimes E_j,
 \qquad \sum_{j=1}^{r}E_j=\mu(G),
\end{equation}
where $K_j\le H$ are compact \emph{open} subgroups, the $E_j$ are mutually orthogonal nonzero projections, and the empty sum represents $\mu=0$. In particular, the $m_H$-density of $\mu$ is  essentially bounded and is given by
\begin{equation}\label{eq:ac-idempotent-density}
 \frac{d\mu}{dm_H}(h)
 =\sum_{j=1}^{r}\frac{\mathbf1_{K_j}(h)}{m_H(K_j)}E_j
 \quad (m_H\text{-a.e. }h\in H).
\end{equation}
\end{theorem}

\begin{proof}
The case $\mu=0$ is immediate. By Proposition~\ref{prop:idempotent-mass}, compress to $\operatorname{Ran}\mu(G)$ and assume that $\mu(G)=I_n$. All convolution operations take place in the closed  subgroup $H$; replacing $G$ by $H$, we may write $\mu=W\,dm$, where $m$ is left Haar measure, $W\in L^1(G,M_n)$ is positive almost everywhere, and $\int_GW\,dm=I_n$. The inversion symmetry $\check\mu=\mu$ follows from Corollary~\ref{cor:matrix-idempotent-l2} and Lemma~\ref{lem:faithful-measure-convolution}.

Put $\lambda=\operatorname{Tr}\mu=(\operatorname{Tr}W)m$ and let $P_y$ be the support projection of $W(y)$, with $P_y=0$ when $W(y)=0$. On $\{\operatorname{Tr}W>0\}$, the measures $\lambda$ and  $m$ have the same null sets, and the support projections of the $\lambda$- and $m$-densities coincide. Hence \eqref{eq:translation-absorption}, multiplied on the left by $W(y)$,
implies
\begin{equation}\label{eq:density-absorption}
 W(y)W(y^{-1}z)=W(y)W(z)
 \quad\text{for }(m\otimes m)\text{-almost every }(y,z)\in G^2.
\end{equation}

This also holds when $W(y)=0$, so the passage from $\lambda$-almost every $y$ to $m$-almost every $y$ is legitimate. Choose nonnegative $\varphi_k\in C_c(G)$ with
$\int_G\varphi_k\,dm=1$ and supports shrinking to the identity, and set
\[
 M_k=\int_GW(z)\varphi_k(z)\,dm(z)\ge0,
 \qquad
 B_k(y)=\int_GW(yz)\varphi_k(z)\,dm(z).
\]

Multiplying \eqref{eq:density-absorption} by $\varphi_k(y^{-1}z)$ and integrating in $z$ implies 
\begin{equation}\label{eq:density-mk}
 W(y)M_k=W(y)B_k(y)
 \qquad\text{for }m\text{-almost every }y.
\end{equation}

To justify the approximation on a possibly non-unimodular group, write $(R_zW)(y)=W(yz)$. Right translation is strongly continuous on $L^1(G,M_n)$ (although it need not be isometric), and its
operator norms are uniformly bounded for $z$ in a fixed compact neighborhood of $e$. Thus Minkowski's inequality gives
\[
 \|B_k-W\|_{L^1(m)}
 \le\int_G\varphi_k(z)\|R_zW-W\|_{L^1(m)}\,dm(z)
 \longrightarrow0.
\]

Passing to a subsequence, which we continue to index by $k$, we have $B_k(y)\to W(y)$ for almost every $y$. Thus,
\begin{equation}\label{eq:density-mk-pointwise}
 W(y)M_k\longrightarrow W(y)^2
 \qquad\text{for almost every }y.
\end{equation}

We do not assert that the products $W(\cdot)B_k(\cdot)$ converge in $L^1$; at this stage $W^2$ need not be integrable. Because $\int_GW\,dm=I_n$, the ranges $\operatorname{Ran}W(y)$,
for $y$ in any conull set, span $\mathbb C^n$. We may therefore choose finitely many points $y_1,\ldots,y_s$ where \eqref{eq:density-mk-pointwise} holds and $\sum_{i=1}^s\operatorname{Ran}W(y_i)=\mathbb C^n$. These points can be selected successively from that conull set: if their ranges did not span $\mathbb C^n$, a nonzero vector orthogonal to all  of them would satisfy $W(y)v=0$ almost everywhere, contrary to $\int_GW\,dm=I_n$. The linear map
\[
 \mathcal L:M_n\longrightarrow (M_n)^s,
 \qquad X\longmapsto(W(y_1)X,\ldots,W(y_s)X)
\]
is injective: its kernel consists of matrices whose ranges lie in $\bigcap_i\ker W(y_i)=\{0\}$. As its domain is finite-dimensional, $\mathcal L$ has a bounded inverse on its range. By \eqref{eq:density-mk-pointwise}, $\mathcal L(M_k)$ converges, hence $M_k$ converges to some $M\in M_n$; moreover, $M\ge0$. Returning to \eqref{eq:density-mk-pointwise}, we obtain
\begin{equation}\label{eq:density-constant-m}
 W(y)M=W(y)^2=MW(y)
 \qquad(m\text{-a.e. }y).
\end{equation}

Functional calculus for the commuting positive matrices $W(y)$ and $M$ implies that
\begin{equation}\label{eq:density-support-factor}
 W(y)=MP_y=P_yM\qquad(m\text{-a.e. }y).
\end{equation}

If $Mv=0$, then $W(y)v=0$ almost everywhere, so $v=(\int W\,dm)v=0$. Thus $M$ is invertible. Only {after} this observation do we infer that $P_y=M^{-1}W(y)$ is integrable and
\begin{equation}\label{eq:density-proj-integral}
 \int_GP_y\,dm(y)=M^{-1}.
\end{equation}

We now consider the commutation. Since $M$ is invertible and commutes with all $P_y$ off a null set, \eqref{eq:density-absorption} and \eqref{eq:density-support-factor} imply
\begin{equation}\label{eq:density-projection-absorption}
 P_yP_{y^{-1}z}=P_yP_z
 \qquad\text{for }(m\otimes m)\text{-almost every }(y,z).
\end{equation}

Inversion preserves the Haar measure class. Since $\check\mu=\mu$, its density at $y$ is a strictly positive scalar multiple of $W(y^{-1})$, and hence $P_{y^{-1}}=P_y$ for almost every $y$. The  substitutions $(g,h)\mapsto(h,g)$ and $(g,h)\mapsto(g,g^{-1}h)$ preserve product Haar-null sets: for the second map, fix $g$ and use left invariance in the $h$-variable, then apply Fubini.  Inversion is likewise nonsingular (with the usual modular Radon--Nikodym factor). It follows that, for almost every $(g,h)$, we may simultaneously apply \eqref{eq:density-projection-absorption} at $(g,h)$ and $(h,g)$ and put $P=P_g$, $S=P_h$, $R=P_{g^{-1}h}=P_{h^{-1}g}$. Then $PR=PS$ and $SR=SP$; taking adjoints gives $RS=PS$. Therefore
\[
 (PS)^2=(PR)(RS)=PR^2S=PRS=P(PS)=PS.
\]

As a product of orthogonal projections, $PS$ is contractive; a contractive idempotent is an orthogonal projection. Thus $PS=SP$. From $W(y)=MP_y=P_yM$ it follows that $W(g)W(h)=W(h)W(g)$ for  almost every $(g,h)$. Fubini and the $L^1$ integrability of $W$ now imply $\mu(D)\mu(D')=\mu(D')\mu(D)$ for all Borel sets $D,D'\subseteq G$.

Theorem~\ref{thm:commuting-idempotents} implies $\mu=\sum_jm_{K_j}\otimes E_j$ with compact subgroups $K_j$. Because $\mu$ is carried by $H$, positivity and $\operatorname{Tr}\mu=\sum_j(\operatorname{rank}E_j)m_{K_j}$ show that every $K_j\subseteq H$. Similarly, $\mu\ll m_H$ implies $m_{K_j}\ll m_H$. Therefore $m_H(K_j)>0$, and a measurable subgroup  of positive left Haar measure is open. Thus each $K_j$ is compact open in $H$. Normalized Haar measure on such a subgroup is $\mathbf1_{K_j}m_H/m_H(K_j)$, yielding \eqref{eq:ac-idempotent-density}. The formula shows in particular that this density is essentially bounded.
\end{proof}

\begin{remark}\label{rem:ac-idempotents-scope}
By the unrestricted matrix Kawada--It\^o theorem (Theorem~\ref{thm:general-matrix-kawada}), the conclusion of Theorem~\ref{thm:ac-matrix-idempotents} holds even \emph{without} second countability:  its Haar channels must be contained in $H$, and absolute continuity forces them to be open in $H$. The proof above supplies an independent density-and-variance argument. In particular it recovers the discrete case with counting measure, and applies to $L^1$ densities on compact subgroups, including compact Lie subgroups. If $H$ is connected, a compact open subgroup must equal
$H$; hence an absolutely continuous positive idempotent carried by a connected $H$ is either zero or (when $H$ is compact) of the form $m_H\otimes E$ for an orthogonal projection $E$.
\end{remark}

\begin{theorem}\label{thm:transverse-matrix-idempotents}
Let $G$ be a second-countable locally compact group and let $K_1,\ldots,K_r$ be compact subgroups, $r\geq1$.  Put $S=\bigcup_{j=1}^rK_j$.  Suppose that the probability measures $m_{K_j}$ are  pairwise mutually singular and that
\begin{equation}\label{eq:transverse-subgroup-condition}
(m_{K_i}*m_{K_j})(S)=0\qquad(i\neq j).
\end{equation}

If $\mu\in M(G,M_n)$ is positive and convolution-idempotent with
\begin{equation}\label{eq:transverse-subgroup-domination}
\mu\ll\sum_{j=1}^rm_{K_j},
\end{equation}
then all the values of $\mu$ commute. More precisely, there are pairwise distinct compact subgroups $H_1,\ldots,H_N$, mutually orthogonal nonzero projections $E_1,\ldots,E_N\in M_n$, and, for each  $\ell\in\{1,\ldots,N\}$, a unique index $j(\ell)$ such that $H_\ell$ is an open subgroup of $K_{j(\ell)}$ and 
\begin{equation}\label{eq:transverse-haar-decomposition}
\mu=\sum_{\ell=1}^Nm_{H_\ell}\otimes E_\ell,
\qquad \sum_{\ell=1}^NE_\ell=\mu(G).
\end{equation}

The empty sum is permitted if $\mu=0$.  If each $K_j$ is connected, then there are pairwise orthogonal (possibly zero) projections $E_j$ such that
\begin{equation}\label{eq:connected-transverse-haar}
\mu=\sum_{j=1}^rm_{K_j}\otimes E_j,
\qquad E_iE_j=0\ (i\neq j),\qquad
\sum_j E_j=\mu(G).
\end{equation}
\end{theorem}

\begin{proof}
Since the finite family $(m_{K_j})_j$ is pairwise mutually singular, there are pairwise disjoint Borel sets $B_j\subseteq K_j$ with $m_{K_j}(B_j)=1$.  Let $\mu_j(D)=\mu(D\cap B_j)$. Then $\mu_j$ are positive matrix measures, $\mu_j\ll m_{K_j}$, $\mu=\sum_j\mu_j$, and $A_j:=\mu_j(G)\geq0$.

For $i\ne j$, each entry of $\mu_i*\mu_j$ is absolutely continuous with respect to $m_{K_i}*m_{K_j}$. Indeed, writing $d\mu_j=W_j\,dm_{K_j}$, its entries are pushforwards of integrable matrix  products $W_i(x)W_j(y)$ under $(x,y)\mapsto xy$. Thus the transverse condition implies
\[
(\mu_i*\mu_j)(S)=0,
\qquad (\mu_i*\mu_j)(G\setminus S)=A_iA_j.
\]

The diagonal convolutions $\mu_j*\mu_j$ are carried by $K_j\subseteq S$.  Since $\mu*\mu=\mu$ is also carried by $S$, 
\[
0=\operatorname{Tr}\bigl((\mu*\mu)(G\setminus S)\bigr)
  =\sum_{i\ne j}\operatorname{Tr}(A_iA_j).
\]

Here the individual mixed matrix measures need not be positive; what matters is that
\[
\operatorname{Tr}(A_iA_j)
   =\|A_i^{1/2}A_j^{1/2}\|_{\mathrm{HS}}^2\geq0.
\]

Every summand must vanish. The vanishing is a consequence of {trace positivity}, not positivity of the generally non-self-adjoint mixed convolution measures. Therefore $A_i^{1/2}A_j^{1/2}=0$, and hence $A_iA_j=0$, for $i\ne j$. Let $Q_j$ denote the support projection of $A_j$ (and $Q_j=0$ when $A_j=0$). Because $0\leq\mu_j(D)\leq A_j$, we obtain
\[
\mu_j(D)=Q_j\mu_j(D)Q_j,
\qquad Q_iQ_j=0\quad(i\ne j).
\]

Thus $\mu_i(D)\mu_j(D')=0$ for $i\ne j$ and all Borel $D,D'$, so all the mixed convolution measures vanish. In particular, $Q_j\mu(D)=\mu(D)Q_j=\mu_j(D)$; hence, using idempotence of $\mu$, 
\[
\mu_j*\mu_j
 = (Q_j\mu Q_j)*(Q_j\mu Q_j)
 =Q_j(\mu*\mu)Q_j
 =\mu_j.
\]

Apply Theorem~\ref{thm:ac-matrix-idempotents} to each nonzero $\mu_j$ with carrying subgroup $K_j$. It gives a finite sum of orthogonal Haar channels on compact open subgroups of $K_j$. The  projections belonging to different $j$ are orthogonal, so these sums combine to give \eqref{eq:transverse-haar-decomposition}. There are no repeated subgroups across different carriers: if a  nonzero Haar probability $m_H$ were dominated by both $m_{K_i}$ and $m_{K_j}$, those two measures would not be mutually singular. This also makes the assignment $j(\ell)$ unique.

When each $K_j$ is connected, it has no proper open subgroup. Thus each nonzero $\mu_j$ is $m_{K_j}\otimes E_j$, where $E_j$ is its (orthogonal) total-mass projection. Assign $E_j=0$ for any vanishing component to obtain \eqref{eq:connected-transverse-haar}.
\end{proof}

\begin{corollary}\label{cor:su2-tori-idempotents}
Let $K_1,\ldots,K_r$ be pairwise distinct maximal tori of $SU(2)$. If a positive convolution-idempotent $\mu\in M(SU(2),M_n)$ satisfies $\mu\ll\sum_{j=1}^rm_{K_j}$, then
\begin{equation}\label{eq:su2-tori-idempotents}
\mu=\sum_{j=1}^rm_{K_j}\otimes E_j,
\qquad E_iE_j=0\quad(i\ne j),\qquad
\sum_jE_j=\mu(SU(2)),
\end{equation}
where the $E_j$ are orthogonal projections, possibly zero. In particular, all the matrix values of $\mu$ commute.
\end{corollary}

\begin{proof}
Each maximal torus in $SU(2)$ is a circle, and two distinct maximal
tori meet in the center $\{\pm I_2\}$. Hence their normalized Haar
measures are mutually singular. Fix $i\ne j$ and any $k$.
For each fixed $x\in K_i$ the fiber
\[
\{y\in K_j:xy\in K_k\}=K_j\cap x^{-1}K_k
\]
is either empty or a coset in $K_j$ of $K_j\cap K_k$.
If $k\ne j$, the intersection $K_j\cap K_k$ is finite and this
fiber has $m_{K_j}$-measure zero for every $x$. If $k=j$,
the fiber is nonempty precisely when $x\in K_j$, an
$m_{K_i}$-null event since $i\ne j$.
By Fubini,
$(m_{K_i}*m_{K_j})(K_k)=0$ for all $k$, and therefore
$(m_{K_i}*m_{K_j})(\bigcup_kK_k)=0$.
Theorem~\ref{thm:transverse-matrix-idempotents} applies.
Since the $K_j$ are connected, it gives the asserted formula.
\end{proof}

\begin{remark}\label{rem:transverse-subgroups}
The transverse condition is sufficient, not necessary, for the orthogonal Haar decomposition; all positive convolution idempotents are classified by Theorem~\ref{thm:general-matrix-kawada}. The  proof here supplies additional localization: every nonzero Haar channel lies in a {unique} dominating subgroup $K_j$ and is open there. Moreover, when $r\ge2$, none of the $K_j$ can be finite: if  $K_i$ were finite, the $x=e$ term in $m_{K_i}*m_{K_j}$ would give $(m_{K_i}*m_{K_j})(S)\ge 1/|K_i|$ for any $j\ne i$. The assumption of second countability can also be dropped if Theorem~\ref{thm:general-matrix-kawada} replaces the cited density classification at the final step; the support-separation proof itself requires no countability assumption.

If at least two channel projections are nonzero for distinct maximal tori of $SU(2)$, there is no single compact subgroup $L$ whose Haar measure dominates $\mu$. Indeed, domination would force both tori to be contained in $L$ and to be open in $L$; their intersection $\{\pm I_2\}$ would then be an open finite subgroup of $L$, making $L$ finite, an impossibility.
\end{remark}

\begin{proposition}\label{prop:homogeneous-idempotents}
Let $G$ be a locally compact group and let $\mu\in M(G,M_n)$ be positive and convolution-idempotent, with $\mu(G)=I_n$. For $0\ne u\in\C^n$ put $G_u=\{h\in G:\delta_h*(\mu u)=\mu u\}$. If $G_u=H$  is independent of $0\ne u\in\C^n$, then $H$ is a compact subgroup and $\mu=m_H\otimes I_n$. 
\end{proposition}

\begin{proof}
By Theorem~\ref{thm:general-matrix-kawada}, one can write
\[
 \mu=\sum_{j=1}^r m_{K_j}\otimes E_j,
 \qquad \sum_{j=1}^rE_j=I_n,
\]
with pairwise distinct compact subgroups $K_j$ and mutually orthogonal nonzero projections $E_j$. For each $j$ choose $0\ne u_j\in\Ran E_j$. By the stabilizer formula \eqref{eq:common-eigenvector-stabilizer}, $G_{u_j}=K_j$. By hypothesis $K_j=H$ for every $j$. As the subgroups in the canonical decomposition are distinct, $r=1$, and the sole projection is $I_n$.  Hence $\mu=m_H\otimes I_n$.

There is also an independent variance-based argument under the second-countability hypothesis of Lemma~\ref{lem:variance-identity}. For a common conull set of $y$ and for all $f\in\mathcal K$, the variance-absorption identity yields $W(y)(f(y^{-1}\cdot)-f)=0$. Thus, if $0\ne v\in\Ran W(y)$, the scalar component $\langle f,v\rangle$ is invariant under left translation by $y$ for all $f\in\mathcal K$. Proposition~\ref{prop:common-eigenvector}(a) implies $y\in G_v=H$. Because $W(y)\ne0$ for $\lambda$-almost every $y$, $\mu$ is carried by $H$. Each vector measure $\mu u$ is
left $H$-invariant by definition of $G_u$ and has total mass $u$. Coordinatewise uniqueness of Haar measure on the compact group $H$ gives $\mu(D)u=m_H(D)u$ for every Borel $D$ and every $u$, which proves the same conclusion without invoking the general classification.
\end{proof}

\begin{remark}\label{rem:kawada-remaining}
For any subgroup $K\le G$, define the subspace
\[
 V_K=\{0\}\cup\{u\in\C^n\setminus\{0\}:K\subseteq G_u\}.
\]

It is linear, since $K\subseteq G_u$ says that the vector measure $\mu u$ is fixed by each left translation $\delta_k$, $k\in K$. The canonical decomposition of  Theorem~\ref{thm:general-matrix-kawada} gives the precise intrinsic formula 
\begin{equation}\label{eq:stabilizer-subspace-channels}
 V_K=\bigoplus_{j:K\subseteq K_j}\Ran E_j,
 \qquad
 Q_{V_K}\mu(D)=\mu(D)Q_{V_K}
 \quad\text{for every Borel }D.
\end{equation}

Indeed, by \eqref{eq:common-eigenvector-stabilizer}, $K\subseteq G_u$ exactly when $E_ju=0$ for all indices with $K\not\subseteq K_j$. In particular, \emph{every} stabilizer subspace $V_K$ reduces all matrix values of $\mu$, and its corner measure is 
\[
 Q_{V_K}\mu(\cdot)Q_{V_K}
   =\sum_{j:K\subseteq K_j}m_{K_j}\otimes E_j.
\]

Suppose that $V=V_K\ne\{0\}$ has minimal positive dimension among such stabilizer subspaces. For $0\ne u\in V$, the subgroup $H=G_u$ contains $K$, so $V_H\subseteq V_K$ and $u\in V_H$. Minimality  gives $V_H=V$. If $0\ne v\in V$, then similarly $V_{G_v}=V$, whence $H=G_u\supseteq G_v$; the inclusion $H\subseteq G_v$ follows from $v\in V_H$. Thus $G_v=H$ for every $0\ne v\in V$. Moreover, \eqref{eq:stabilizer-subspace-channels} and minimality imply that $V=\Ran E_j$ for a unique channel $j$. To see this, if two distinct channels $i,j$ occurred in $V$, then $0\ne V_{K_i},V_{K_j}\subseteq V$; minimality would force $V_{K_i}=V_{K_j}=V$, and hence both $K_i\subseteq K_j$ and $K_j\subseteq K_i$, a contradiction. Thus $H=K_j$ and
\[
 Q_V\mu|_H Q_V=m_H\otimes I_V,
 \qquad C_V:=\bigl.Q_V\mu(H)Q_V\bigr|_V=I_V.
\]

This also proves directly that $V$ is reducing. It is useful to distinguish this conclusion from what the {variance argument alone} proves, without the general classification. In the second-countable normalized setting, the homogeneous-stabilizer property on a minimal $V$ implies $Q_V\mu|_HQ_V=m_H\otimes C_V$ with $0\le C_V\le I_V$. Writing $Q=1_{L^2(G)}\otimes Q_V$ and  $T=T_\mu$, the fixed-space characterization implies that $QTQ$ has range in the $H$-invariant subspace $\Ran(T_{m_H}\otimes Q_V)$. As $0\le QTQ\le Q$, it follows that
\[
 0\le T_V:=QTQ\le T_{m_H}\otimes Q_V.
\]

The variance absorption also shows $\Ran W(y)\cap V=\{0\}$ for $\lambda$-almost every $y\notin H$. By itself this last assertion does not imply $\Ran W(y)\perp V$ there; for example, $\operatorname{span}(e_1+e_2)$ intersects $\operatorname{span}(e_1)$ trivially but is not orthogonal to it. Hence these variance estimates alone do not prove $C_V=I_V$. The general classification,  or equivalently the intrinsic formula \eqref{eq:stabilizer-subspace-channels}, resolves this point completely. The density and transverse-subgroup methods in Theorems~\ref{thm:ac-matrix-idempotents} and \ref{thm:transverse-matrix-idempotents} provide independent proofs under additional hypotheses; finite unions of distinct maximal tori of $SU(2)$ are  treated by Corollary~\ref{cor:su2-tori-idempotents}.
\end{remark}

\subsection{Atomic splitting and countably transverse subgroups}
\label{subsec:atomic-splitting}

The next two results remove the Haar-domination hypothesis in two
directions. The first one splits off the atoms on an arbitrary locally
compact group. The second one replaces the Haar transversality
condition~\eqref{eq:transverse-subgroup-condition} by a purely
set-theoretic condition on the carrying subgroups. The common tool is
the following positivity of traces.

\begin{lemma}\label{lem:trace-coupling}
Let $G$ be a locally compact group and let $\mu_1,\mu_2\in M(G,M_n)$ be finite positive matrix-valued Radon measures. Set $\lambda_i=\operatorname{Tr}\mu_i$ and choose Borel representatives  $W_i=d\mu_i/d\lambda_i$, with $W_i\geq0$ and $\operatorname{Tr}W_i=1$ $\lambda_i$-almost everywhere. Write $\Lambda=\lambda_1\widehat\otimes\lambda_2$ for the finite {Radon product measure}  on $G\times G$. Then for every Borel $D\subseteq G$,
\begin{equation}\label{eq:trace-coupling-formula}
 \operatorname{Tr}\bigl((\mu_1*\mu_2)(D)\bigr)
 =\int_{G\times G}\mathbf1_D(yz)
 \bigl\|W_1(y)^{1/2}W_2(z)^{1/2}\bigr\|_{\mathrm{HS}}^2
 \,d\Lambda(y,z)\geq0.
\end{equation}

Moreover, the following are equivalent 
\begin{enumerate}
\item[(i)] $\operatorname{Tr}((\mu_1*\mu_2)(D))=0$;
\item[(ii)] $W_1(y)W_2(z)=0$ for $\Lambda$-almost every $(y,z)$ satisfying $yz\in D$;
\item[(iii)] $(\mu_1*\mu_2)(D\cap E)=0$ for every Borel $E\subseteq G$.
\end{enumerate}

In particular, the scalar measure $D\mapsto\operatorname{Tr}((\mu_1*\mu_2)(D))$ is positive, and
\begin{equation}\label{eq:trace-coupling-total}
 \mu_1*\mu_2=0
 \quad\Longleftrightarrow\quad
 \mu_1(G)\mu_2(G)=0
 \quad\Longleftrightarrow\quad
 \mu_2*\mu_1=0.
\end{equation}
%No second-countability assumption is needed.
\end{lemma}

\begin{proof}
By positivity, every entry of $\mu_i$ is absolutely continuous with respect to $\lambda_i$, so the Radon--Nikodym theorem gives the matrix density $W_i$. It is positive with trace one almost  everywhere; in particular $\|W_i\|\leq1$ almost everywhere. The zero measure is covered by the same argument with an arbitrary version of its density.

We use the Radon product of Lemma~\ref{lem:radon-prerequisites}(b); it agrees with the ordinary product measure on the product $\sigma$-algebra, but is defined on the whole Borel $\sigma$-algebra  of $G\times G$. Entrywise, the product measure $\mu_{1,ab}\widehat\otimes\mu_{2,cd}$ has density $(y,z)\mapsto (W_1(y))_{ab}(W_2(z))_{cd}$ with respect to $\Lambda$. Indeed, the identity holds  first on Borel rectangles by the Radon--Nikodym formulas and then on all Borel sets by uniqueness of finite Radon measures from their tensor tests. Pushing forward under the continuous  multiplication map $(y,z)\mapsto yz$, hence, for each Borel $D$,
\[
 (\mu_1*\mu_2)(D)
 =\int_{G\times G}\mathbf1_D(yz)W_1(y)W_2(z)\,d\Lambda(y,z).
\]

The integrand is bounded and Borel measurable up to a $\Lambda$-null set. Taking the trace and using, for positive matrices $A,B$, 
\[
 \operatorname{Tr}(AB)
 =\operatorname{Tr}(B^{1/2}AB^{1/2})
 =\|A^{1/2}B^{1/2}\|_{\mathrm{HS}}^2\geq0
\]
proves \eqref{eq:trace-coupling-formula}. A nonnegative integrable function has integral zero over a measurable set if and only if it vanishes almost everywhere there. Since $\operatorname{Tr}(AB)=0$ if and only if $A^{1/2}B^{1/2}=0$  if and only if $AB=0$, for positive $A,B$, conditions (i) and (ii) are equivalent. Under (ii), the matrix integrand in the convolution formula vanishes $\Lambda$-almost everywhere on the inverse image of $D$, and therefore $(\mu_1*\mu_2)(D\cap E)=0$ for every Borel $E$. Conversely (iii) with $E=G$ implies (i).

Finally, put $A_i=\mu_i(G)\geq0$. Multiplicativity of total mass gives $(\mu_1*\mu_2)(G)=A_1A_2$. If $A_1A_2=0$, positivity implies that their support projections $Q_1,Q_2$ are orthogonal, and $\mu_i(E)=Q_i\mu_i(E)Q_i$ for every Borel $E$. Thus all products $\mu_1(E)\mu_2(F)$ and $\mu_2(F)\mu_1(E)$ vanish, giving both zero convolutions. Conversely, if either convolution is zero, $\operatorname{Tr}(A_1A_2)=0$, and positivity of $A_i$ implies $A_1A_2=0$. This proves \eqref{eq:trace-coupling-total}.
\end{proof}

\begin{remark}
The scalar measure $\operatorname{Tr}(\mu_1*\mu_2)$ is positive, but $\mu_1*\mu_2$ need not be a positive matrix measure. For instance, take $\mu_1=\delta_e\otimes A$ and $\mu_2=\delta_e\otimes  B$, where $A,B\in M_2$ are positive and $AB\neq BA$; then $\mu_1*\mu_2=\delta_e\otimes AB$ is not self-adjoint. The equivalence between zero trace on $D$ and vanishing of the {restriction} of  the convolution to $D$, rather than only its value at $D$, is the feature used in support-splitting arguments.
\end{remark}

\begin{theorem}\label{thm:atomic-splitting}
Let $G$ be any locally compact group, and let $\mu\in M(G,M_n)$ be positive and convolution-idempotent. Write $\mu=\mu_a+\mu_c$, where $\mu_a$ is purely atomic and $\mu_c$ has no point masses.  Then 
\begin{enumerate}
\item[(a)] $\mu_a$ and $\mu_c$ are positive convolution idempotents.
\item[(b)] The matrices $P_a=\mu_a(G)$, and $P_c=\mu_c(G)$ are orthogonal projections satisfying $P_aP_c=P_cP_a=0$, and $P_a+P_c=\mu(G)$. Moreover, for every Borel $D\subseteq G$,
\[
 \mu_a(D)=P_a\mu(D)P_a,\qquad
 \mu_c(D)=P_c\mu(D)P_c.
\]

In particular, we have $\mu_a*\mu_c=\mu_c*\mu_a=0$.  
\item[(c)] There exist pairwise distinct finite subgroups $F_1,\ldots,F_s$ of $G$ and mutually orthogonal nonzero projections $Q_1,\ldots,Q_s\in M_n$ such that 
\[
 \mu_a=\sum_{i=1}^s m_{F_i}\otimes Q_i,
 \qquad \sum_{i=1}^sQ_i=P_a.
\]

The empty sum is allowed. In particular, $\operatorname{supp}\mu_a$ is finite.

\item[(d)] Writing $e$ for the identity of $G$ and $s(A)$ for the support projection of a positive matrix $A$,
\[
 \mu(\{e\})
   =\sum_{i=1}^s\frac{1}{|F_i|}Q_i,
 \qquad
 P_a=s\bigl(\mu(\{e\})\bigr).
\]

Thus, $\mu$ is atomless if and only if
$\mu(\{e\})=0$; equivalently, any point mass anywhere in $G$
forces a nonzero point mass at the identity. With $P=\mu(G)$,
\[
 \mu_a(D)=s(\mu(\{e\}))\,
          \mu(D)\,s(\mu(\{e\})),
\]
and
\[
 \mu_c(D)=
 (P-s(\mu(\{e\})))\,
 \mu(D)\,
 (P-s(\mu(\{e\}))).
\]
\end{enumerate}
%No second-countability or discreteness assumption on $G$ is needed.
\end{theorem}

\begin{proof}
Let $A=\{x\in G:\mu(\{x\})\ne0\}$. Since $\operatorname{Tr}\mu$ is a finite positive scalar measure, $A$ is countable (each $\{x:\operatorname{Tr}\mu(\{x\}) \ge1/k\}$ is finite). In particular $A$ is Borel, whether or not $G$ is second countable. The atomic and atomless parts are therefore given by
\[
 \mu_a(D)=\mu(D\cap A),\qquad
 \mu_c(D)=\mu(D\setminus A).
\]

In particular, both are positive finite matrix measures. We first establish idempotence of $\mu_a$. For finite complex measures $\nu,\rho$, if $\nu$ has no point masses, then
\[
 (\nu*\rho)(\{x\})
   =\int_G\nu(\{xz^{-1}\})\,d\rho(z)=0,
\]
and similarly
\[
 (\rho*\nu)(\{x\})
   =\int_G\nu(\{z^{-1}x\})\,d\rho(z)=0.
\]

Applying this entrywise shows that convolution by a matrix measure without point masses produces a measure without point masses, on either side. On the other hand, $\mu_a*\mu_a$ is purely atomic,
because it is carried by the countable set $AA$. Thus the atomic part of 
\[
 \mu*\mu
 =\mu_a*\mu_a+\mu_a*\mu_c
  +\mu_c*\mu_a+\mu_c*\mu_c
\]
is precisely $\mu_a*\mu_a$. Since $\mu*\mu=\mu$, uniqueness of the atomic decomposition gives
\begin{equation}\label{eq:atomic-part-idempotence}
 \mu_a*\mu_a=\mu_a.
\end{equation}

Comparing the remaining parts gives
\begin{equation}\label{eq:atomless-part-identity}
 \mu_c=\mu_a*\mu_c+\mu_c*\mu_a+\mu_c*\mu_c.
\end{equation}

Let $\Gamma$ be the subgroup algebraically generated by $A$. It is countable. Regard $\Gamma$ as a discrete group. The measure $\mu_a$ defines a finite positive matrix-valued measure on $\Gamma$,  and its convolution is unchanged because convolution of atomic measures depends only on the group multiplication. By Theorem~\ref{thm:discrete-matrix-idempotents} and  \eqref{eq:atomic-part-idempotence}, there are finite subgroups $F_1,\ldots,F_s\leq\Gamma$ and mutually orthogonal nonzero projections $Q_i$ such that
\[
 \mu_a=\sum_{i=1}^s m_{F_i}\otimes Q_i.
\]

Equal subgroups may be grouped together. Thus,
\[
 P_a=\mu_a(G)=\sum_{i=1}^sQ_i
\]
is an orthogonal projection, and $\mu_a$ has finite support. This proves (c). Put $P_c=\mu_c(G)\geq0$. Taking total masses in \eqref{eq:atomless-part-identity} implies
\begin{equation}\label{eq:atomic-mass-relation}
 P_c=P_aP_c+P_cP_a+P_c^2.
\end{equation}

Relative to $\mathbb C^n=\operatorname{Ran}P_a\oplus\ker P_a$, write
\[
 P_c=
 \begin{pmatrix}
  \alpha&\beta\\
  \beta^*&\gamma
 \end{pmatrix},
 \qquad \alpha\geq0,\quad\gamma\geq0.
\]

The upper-left block of \eqref{eq:atomic-mass-relation} is $\alpha=2\alpha+\alpha^2+\beta\beta^*$, and hence $\alpha+\alpha^2+\beta\beta^*=0$. Each term on the left is positive semidefinite.
Therefore $\alpha=0$ and $\beta=0$. The lower-right block now yields $\gamma^2=\gamma$. Thus $P_c$ is an orthogonal projection supported on $\ker P_a$, and $P_aP_c=P_cP_a=0$. For every Borel $D$, positivity implies that $0\leq\mu_a(D)\leq P_a$, and $0\leq\mu_c(D)\leq P_c$. Since $P_a,P_c$ are orthogonal projections, it follows that $\mu_a(D)=P_a\mu_a(D)P_a$, and $\mu_c(D)=P_c\mu_c(D)P_c$. The orthogonality of $P_a$ and $P_c$ implies that all mixed products of values of $\mu_a$ and $\mu_c$ vanish. Thus, $\mu_a*\mu_c=\mu_c*\mu_a=0$. Equation \eqref{eq:atomless-part-identity} therefore reduces to $\mu_c*\mu_c=\mu_c$. This proves (a). The same corner identities give
\[
 \mu_a(D)=P_a\mu(D)P_a,\qquad
 \mu_c(D)=P_c\mu(D)P_c,
\]
and $P_a+P_c=\mu(G)$, establishing (b).

Finally, every finite subgroup $F_i$ contains $e$,
and its normalized Haar measure satisfies
$m_{F_i}(\{e\})=|F_i|^{-1}$. As $\mu_c$ has no atoms,
\[
 \mu(\{e\})=\mu_a(\{e\})
 =\sum_{i=1}^s\frac{1}{|F_i|}Q_i.
\]

The coefficients are strictly positive and the $Q_i$ are mutually orthogonal. Thus
\[
 s(\mu(\{e\}))=\sum_{i=1}^sQ_i=P_a.
\]

Since $P_c=\mu(G)-P_a$, the intrinsic compression formulas in (d) follow from (b).
\end{proof}

\begin{theorem}\label{thm:countably-transverse-idempotents}
Let $G$ be a second-countable locally compact group,
and let $L_1,\ldots,L_r$ be closed subgroups such that
$L_i\cap L_j$ is countable whenever $i\ne j$.
Let $\mu\in M(G,M_n)$ be positive and
convolution-idempotent, and suppose that $\mu$
is carried by
\[
 S=\bigcup_{j=1}^r L_j.
\]

No absolute-continuity hypothesis is imposed. Then there is a decomposition
\begin{equation}\label{eq:countably-transverse-split}
 \mu=\mu_a+\sum_{j=1}^r\mu_j
\end{equation}
with the following properties:
\begin{enumerate}
\item[(a)]
$\mu_a$ is the purely atomic part of $\mu$.
It is a positive convolution idempotent with
finite support and has the decomposition
\[
 \mu_a=\sum_{\ell=1}^s m_{F_\ell}\otimes Q_\ell,
\]
where the $F_\ell$ are finite subgroups of $G$
and the $Q_\ell$ are mutually orthogonal
nonzero projections. The empty sum is permitted.

\item[(b)]
Each $\mu_j$ is an atomless positive
convolution idempotent carried by $L_j$. Moreover, it is
intrinsically determined by
$\mu_j(D)=\mu_c(D\cap L_j)$ for every Borel $D$;
thus the splitting does not depend on an ordering of the $L_j$.

\item[(c)]
The total masses
\[
 P_a=\mu_a(G),\qquad P_j=\mu_j(G)
\]
are mutually orthogonal projections and
\[
 P_a+\sum_{j=1}^rP_j=\mu(G).
\]
In particular, the mixed convolutions between
different summands in
\eqref{eq:countably-transverse-split} vanish.

\item[(d)]
Each $\mu_j$ admits an orthogonal Haar
decomposition
\[
 \mu_j=\sum_{\ell=1}^{s_j}
       m_{H_{j,\ell}}\otimes E_{j,\ell},
\]
where the $H_{j,\ell}$ are infinite compact
subgroups of $L_j$ and the coefficient
projections are mutually orthogonal.
Each such atomless Haar channel belongs
to exactly one of the carrying subgroups
$L_1,\ldots,L_r$.
\end{enumerate}

Thus, all values of $\mu$ commute, and $\mu$ admits the orthogonal Haar decomposition of Theorem~\ref{thm:general-matrix-kawada}. The splitting in (a)--(c) is obtained
independently of that general classification.
\end{theorem}

\begin{proof}
By Theorem~\ref{thm:atomic-splitting}, write $\mu=\mu_a+\mu_c$, where $\mu_a$ and $\mu_c$ are positive convolution idempotents with orthogonal total-mass projections. Moreover $\mu_a$ has finite support and the finite-subgroup decomposition asserted in (a).
Both parts are carried by $S$. It remains to split the atomless measure $\mu_c$. Define the disjoint Borel sets 
\[
 S_j=L_j\setminus\bigcup_{i<j}L_i,
 \qquad 1\le j\le r,
\]
and put $\mu_j(D)=\mu_c(D\cap S_j)$. Then $\mu_j$ is positive and atomless, is carried by $L_j$, and $\mu_c=\sum_{j=1}^r\mu_j$. Because each $L_i\cap L_j$ is countable and $\mu_c$ has no
point masses, $\mu_c(D\cap L_i\cap L_j)=0$ for all Borel $D$ and $i\ne j$. It follows that $\mu_j(D)=\mu_c(D\cap L_j)$, so this decomposition is independent of the arbitrary ordering used to define the disjoint sets $S_j$. Write $A_j=\mu_j(G)\ge0$, and $\lambda_j=\operatorname{Tr}\mu_j$. Fix $i\ne j$. We claim that
\[
 (\lambda_i\otimes\lambda_j)
 \bigl(\{(y,z)\in L_i\times L_j:yz\in S\}\bigr)=0.
\]

For $y\in L_i$ and $1\le k\le r$, consider
\[
 Z_{y,k}
 =\{z\in L_j:yz\in L_k\}
 =L_j\cap y^{-1}L_k.
\]

If $k\ne j$, this set is either empty or a coset in $L_j$ of $L_j\cap L_k$. It is therefore countable and hence $\lambda_j$-null, because $\lambda_j$ has no atoms. If $k=j$, the fibre $Z_{y,j}$ is empty unless $y\in L_j$. But $\lambda_i$ is carried by $L_i$, and $\lambda_i(L_i\cap L_j)=0$ because $L_i\cap L_j$ is countable and $\lambda_i$ is atomless. Since the union defining $S$ is finite,
Fubini's theorem gives the claim. By Lemma~\ref{lem:trace-coupling}, the restriction of $\mu_i*\mu_j$ to $S$ vanishes for $i\ne j$. Thus $(\mu_i*\mu_j)(S)=0$. Multiplicativity of total mass gives
\[
 (\mu_i*\mu_j)(G\setminus S)=A_iA_j.
\]

The measure $\mu_c$ is a convolution idempotent carried by $S$, and each $\mu_j*\mu_j$ is carried by $L_j\subseteq S$. Therefore
\[
 \begin{aligned}
 0
 &=\operatorname{Tr}
       ((\mu_c*\mu_c)(G\setminus S))\\
 &=\sum_{i\ne j}\operatorname{Tr}(A_iA_j).
 \end{aligned}
\]

For positive matrices $A_i,A_j$, we have $\operatorname{Tr}(A_iA_j) =\|A_i^{1/2}A_j^{1/2}\|_{\mathrm{HS}}^2 \ge0$. Hence each term vanishes, and
\begin{equation}\label{eq:countable-corner-orthogonality}
 A_iA_j=0\qquad(i\ne j).
\end{equation}

Let $\Pi_j=s(A_j)$ be the support projection of $A_j$. The preceding identity implies $\Pi_i\Pi_j=0$, for $i\ne j$. Since $0\le\mu_j(D)\le A_j$, we have $\mu_j(D)=\Pi_j\mu_j(D)\Pi_j$. It follows that
\[
 \mu_j(D)=\Pi_j\mu_c(D)\Pi_j.
\]

All mixed products of the matrix values of different $\mu_j$ vanish. Using $\mu_c*\mu_c=\mu_c$, we obtain
\[
 \begin{aligned}
 \mu_j*\mu_j
 &=\Pi_j(\mu_c*\mu_c)\Pi_j\\
 &=\Pi_j\mu_c\Pi_j\\
 &=\mu_j.
 \end{aligned}
\]

Thus each $\mu_j$ is a positive convolution idempotent on $L_j$. Its total mass satisfies $A_j^2=A_j$, so $A_j=\Pi_j$ is an orthogonal projection.  The atomic-splitting theorem gives $P_aP_j=0$  for every $j$. Together with \eqref{eq:countable-corner-orthogonality}, this proves (b) and (c).

Apply Theorem~\ref{thm:general-matrix-kawada} to each $\mu_j$, viewed as a measure on the closed subgroup $L_j$. This implies that
\[
 \mu_j=\sum_{\ell=1}^{s_j}
       m_{H_{j,\ell}}\otimes E_{j,\ell}
\]
with compact $H_{j,\ell}\le L_j$. Since $\mu_j$ is atomless, each $H_{j,\ell}$ must be infinite.  If an infinite compact subgroup $H_{j,\ell}$ were contained in some $L_i$ with $i\ne j$, it would  be contained in the countable group $L_i\cap L_j$. This is impossible: an infinite compact Hausdorff group is uncountable. Hence the carrying subgroup is unique. 

The projections occurring in different components, including $\mu_a$, are orthogonal. Combining the Haar decompositions proves (d) and the final assertion.
\end{proof}

\begin{corollary}\label{cor:countably-transverse-examples}
Let $\mu\in M(G,M_n)$ be positive and convolution-idempotent.
\begin{enumerate}
\item[(a)] Suppose $G=SU(2)$ and $\mu$ is carried by finitely many distinct maximal tori $K_1,\ldots,K_r$. 
\item[(b)] Suppose $G=SO(3)$ and $\mu$ is carried by finitely many circle subgroups $K_{\omega_1},\ldots, K_{\omega_r}$ associated with pairwise non-proportional rotation axes. 
\item[(c)] More generally, suppose $G$ is second countable and $\mu$ is carried by finitely many closed abelian subgroups having pairwise countable intersections.
\end{enumerate}

In all three cases, $\mu$ has commuting matrix values and an orthogonal Haar decomposition, without any assumption of Haar domination. In cases (a) and (b), the stronger formula holds:
\begin{equation}\label{eq:tori-atomic-continuous}
 \mu=\mu_a+\sum_{j=1}^r m_{K_j}\otimes E_j,
\end{equation}
where $\mu_a$ is a finite orthogonal sum of Haar measures of finite subgroups, the $E_j$ are orthogonal projections (possibly zero), and their ranges are orthogonal to the coefficient spaces of  $\mu_a$. In particular, the atomless part of $\mu$ is absolutely continuous with respect to $\sum_jm_{K_j}$.
\end{corollary}

\begin{proof}
For (a), any two distinct maximal tori of $SU(2)$ intersect exactly in its center $\{\pm I\}$. Thus their intersections are finite. 

 For (b), a nonidentity rotation in $SO(3)$ has a unique rotation axis. Hence circle groups corresponding to distinct unoriented axes intersect only in the identity. The hypotheses of Theorem~\ref{thm:countably-transverse-idempotents} are therefore satisfied in all three cases.

The relevant subgroups in (a) and (b) are isomorphic to the circle group $\mathbb T$. Every proper closed subgroup of $\mathbb T$ is finite. Since each atomless component $\mu_j$ has a Haar  decomposition on infinite compact subgroups of its carrying circle, only the whole circle can occur. Therefore
\[
 \mu_j=m_{K_j}\otimes E_j
\]
for an orthogonal projection $E_j$, which may be zero. Combining this with the atomic splitting gives \eqref{eq:tori-atomic-continuous}. In (c), the abelian classification on each carrying  subgroup supplies the claimed Haar decomposition.
\end{proof}

\begin{remark}\label{rem:countably-transverse-scope}
Theorem~\ref{thm:countably-transverse-idempotents} gives an independent orthogonal support-splitting mechanism for finitely many closed subgroups whose pairwise intersections are countable.  It  differs from Theorem~\ref{thm:transverse-matrix-idempotents}, which assumes domination by subgroup Haar measures and the separate mixed-convolution condition \eqref{eq:transverse-subgroup-condition}. Neither hypothesis set implies the other in general.

In particular, the components $\mu_j$ obtained from the atomless part are canonical: restriction of $\mu_c$ to $L_j$ gives $\mu_j$ regardless of how the carrying subgroups were ordered.  The  atomic component may involve finite subgroups that are not contained in any individual carrying subgroup. For example, the quaternion subgroup $Q_8\le SU(2)$ is contained in the union of three  distinct maximal tori, but not in any single one of them.

In cases (a) and (b) of Corollary~\ref{cor:countably-transverse-examples}, the atomless components are Haar measures of the carrying circles, with orthogonal matrix coefficients. They are singular  with respect to ambient Haar measure on $SU(2)$ or $SO(3)$, but they are not singular with respect to Haar measure on their individual circles. The unrestricted classification of  Theorem~\ref{thm:general-matrix-kawada} also covers nested subgroups, nontransverse configurations, infinitely many possible carrying subgroups, and groups that are not second countable. The  significance of the present result is its independent support-localization argument.
\end{remark}

\section*{Declarations}
\noindent\textbf{Funding.} No funding was received for this work.

\smallskip\noindent\textbf{Competing interests.} The author declares no relevant competing interests.

\smallskip\noindent\textbf{Data availability.} No research datasets were generated or analyzed. 

\smallskip\noindent\textbf{Use of AI tools.} ChatGPT (OpenAI) was used to assist in English editing.


\begin{thebibliography}{99}

\bibitem{ChoquetDeny}
G.~Choquet, J.~Deny, Sur l'\'equation de convolution $\mu=\mu*\sigma$,
\emph{C. R. Acad. Sci. Paris} \textbf{250} (1960), 799--801.

\bibitem{Chu}
C.-H.~Chu, \emph{Matrix Convolution Operators on Groups},
Lecture Notes in Mathematics \textbf{1956}, Springer, Berlin, 2008,
\url{https://doi.org/10.1007/978-3-540-69798-5}.

\bibitem{ChuHilberdinkHowroyd}
C.-H.~Chu, T.~Hilberdink and J.~Howroyd,
A matrix-valued Choquet--Deny theorem,
\emph{Proc. Amer. Math. Soc.} \textbf{129} (2001), 229--235,
\url{https://doi.org/10.1090/S0002-9939-00-05694-X}.

\bibitem{ChuMatrixHarmonic}
C.-H.~Chu, Matrix-valued harmonic functions on groups,
\emph{J. Reine Angew. Math.} \textbf{552} (2002), 15--52,
\url{https://doi.org/10.1515/crll.2002.089}.

\bibitem{CohenIdempotent}
P.~J.~Cohen, On a conjecture of Littlewood and idempotent measures,
\emph{Amer. J. Math.} \textbf{82} (1960), 191--212,
\url{https://doi.org/10.2307/2372731}.

\bibitem{EbadianJabbariBochner}
A.~Ebadian and A.~Jabbari, Representations and Bochner's theorem for
matrix valued group algebras, \emph{J. Math. Anal. Appl.} \textbf{528}
(2023), 127542, \url{https://doi.org/10.1016/j.jmaa.2023.127542}.

\bibitem{EbadianJabbariMIA}
A.~Ebadian and A.~Jabbari, Matrix valued conjugate convolution operators
on matrix valued $L^p$-spaces,
\emph{Math. Inequal. Appl.} \textbf{24} (2021), no.~3, 827--844,
\url{https://doi.org/10.7153/mia-2021-24-57}.

\bibitem{EbadianJabbariOaM}
A.~Ebadian and A.~Jabbari, Matrix valued $p$-convolution operators,
\emph{Oper. Matrices} \textbf{14} (2020), no.~1, 117--128,
\url{https://doi.org/10.7153/oam-2020-14-09}.

\bibitem{EngelNagel}
K.-J.~Engel, R.~Nagel, \emph{One-Parameter Semigroups for Linear Evolution
Equations}, Graduate Texts in Mathematics \textbf{194}, Springer, 2000.

\bibitem{Folland}
G.~B.~Folland, \emph{A Course in Abstract Harmonic Analysis},
2nd ed., CRC Press, Boca Raton, 2016.

\bibitem{FollandRA}
G.~B.~Folland, \emph{Real Analysis: Modern Techniques and Their
Applications}, 2nd ed., Wiley, New York, 1999.

\bibitem{FranzSkalski}
U.~Franz and A.~Skalski, On idempotent states on quantum groups,
\emph{J. Algebra} \textbf{322} (2009), 1774--1802; arXiv:0808.1683.

\bibitem{Greenleaf}
F.~P.~Greenleaf, Norm decreasing homomorphisms of group algebras,
\emph{Pacific J. Math.} \textbf{15} (1965), no.~4, 1187--1219.

\bibitem{HewittRoss}
E.~Hewitt and K.~A.~Ross, \emph{Abstract Harmonic Analysis I}, 2nd ed.,
Grundlehren der mathematischen Wissenschaften \textbf{115}, Springer,
1979.

\bibitem{Heyer}
H.~Heyer, \emph{Probability Measures on Locally Compact Groups},
Ergebnisse der Mathematik und ihrer Grenzgebiete \textbf{94}, Springer,
Berlin, 1977.

\bibitem{HostIdempotents}
B.~Host, Le th\'eor\`eme des idempotents dans $B(G)$,
\emph{Bull. Soc. Math. France} \textbf{114} (1986), 215--223,
\url{https://doi.org/10.24033/bsmf.2055}.

\bibitem{JabbariPositiveType}
A.~Jabbari, Positive type and positive definite functions on matrix
valued group algebras, \emph{Results Math.} \textbf{75} (2020),
Paper No.~149, \url{https://doi.org/10.1007/s00025-020-01278-1}.

\bibitem{KasprzakKhosraviSoltan}
P.~Kasprzak, F.~Khosravi, and P.~M.~So\l tan,
Kawada--It\^o--Kelley theorem for quantum semigroups,
\emph{J. Math. Anal. Appl.} \textbf{483} (2020), 123594,
\url{https://doi.org/10.1016/j.jmaa.2019.123594}.

\bibitem{KawadaIto}
Y.~Kawada, K.~It\^o, On the probability distribution on a compact group I,
\emph{Proc. Phys.-Math. Soc. Japan} \textbf{22} (1940), 977--998.

\bibitem{Kelley}
J.~L.~Kelley, Averaging operators on $C_\infty(X)$,
\emph{Illinois J. Math.} \textbf{2} (1958), 214--223.

\bibitem{Krengel}
U.~Krengel, \emph{Ergodic Theorems}, de Gruyter, Berlin, 1985.

\bibitem{Loynes}
R.~M.~Loynes, Fourier transforms and probability theory on a
noncommutative locally compact topological group,
\emph{Ark. Mat.} \textbf{5} (1963), 37--42.

\bibitem{Mustafayev2022}
H.~Mustafayev, A note on the Kawada--It\^o theorem,
\emph{Statist. Probab. Lett.} \textbf{181} (2022), 109261,
\url{https://doi.org/10.1016/j.spl.2021.109261}.

\bibitem{NeufangSalmiSkalskiSpronk}
M.~Neufang, P.~Salmi, A.~Skalski and N.~Spronk, Contractive idempotents on
locally compact quantum groups, \emph{Indiana Univ. Math. J.} \textbf{62}
(2013), no.~6, 1983--2002; arXiv:1209.6508.

\bibitem{Paulsen}
V.~Paulsen, \emph{Completely Bounded Maps and Operator Algebras},
Cambridge Studies in Advanced Mathematics \textbf{78}, Cambridge
University Press, 2002.

\bibitem{Pym}
J.~S.~Pym, Idempotent measures on semigroups,
\emph{Pacific J. Math.} \textbf{12} (1962), 685--698.

\bibitem{Salmi}
P.~Salmi, Idempotent states on locally compact groups and quantum
groups, in \emph{Algebraic Methods in Functional Analysis}, Oper. Theory
Adv. Appl. \textbf{233}, Birkh\"auser, Basel, 2014;
arXiv:1209.0314.

\bibitem{SalmiSkalski}
P.~Salmi and A.~Skalski, Idempotent states on locally compact quantum
groups, \emph{Q. J. Math.} \textbf{63} (2012), 1009--1032; arXiv:1102.2051.

\bibitem{SpronkCommuting}
N.~Spronk, Commuting contractive idempotents in measure algebras,
preprint, arXiv:1503.06859.

\bibitem{Wendel}
J.~G.~Wendel, Left centralizers and isomorphisms of group algebras,
\emph{Pacific J. Math.} \textbf{2} (1952), 251--261.

\bibitem{Wendel54}
J.~G.~Wendel, Haar measure and the semigroup of measures on a compact
group, \emph{Proc. Amer. Math. Soc.} \textbf{5} (1954), 923--929.

\end{thebibliography}
\end{document}